\documentclass[a4, 11pt, leqno]{article}
\usepackage{amsmath}
\usepackage{amssymb}
\newtheorem{df}{Definition}[section]
\newtheorem{pr}[df]{Proposition}
\newtheorem{Th}[df]{Theorem}

\newtheorem{lm}[df]{Lemma}

\begin{document}

\title{\bf Motion of a Vortex Filament with Axial Flow \\in the Half Space}
\author{Masashi A{\sc iki} and Tatsuo I{\sc guchi}}
\date{}
\maketitle
\vspace*{-0.5cm}
\begin{center}
Department of Mathematics, Faculty of Science and Technology, Keio University, \\
3-14-1 Hiyoshi, Kohoku-ku, Yokohama, 223-8522, J{\sc apan}
\end{center}

\begin{abstract}
We consider a nonlinear third order dispersive equation which models the motion of a vortex filament immersed in an incompressible and inviscid fluid 
occupying the three dimensional half space. 
We prove the unique solvability of initial-boundary value problems as an attempt to analyze the motion of a tornado.

\end{abstract}

\section{Introduction}
In this paper, we prove the unique solvability locally in time of the following initial-boundary value problems. For \( \alpha <0\),
\begin{eqnarray}
\left\{
\begin{array}{ll}
\mbox{\mathversion{bold}$x$}_{t}= \mbox{\mathversion{bold}$x$}_{s}\times \mbox{\mathversion{bold}$x$}_{ss} + \alpha \big\{ \mbox{\mathversion{bold}$x$}_{sss}
+ \frac{3}{2}\mbox{\mathversion{bold}$x$}_{ss}\times \big( \mbox{\mathversion{bold}$x$}_{s}\times \mbox{\mathversion{bold}$x$}_{ss}\big) \big\}, &
s>0,t>0,\\
\mbox{\mathversion{bold}$x$}(s,0) = \mbox{\mathversion{bold}$x$}_{0}(s), & s>0, \\
\mbox{\mathversion{bold}$x$}_{ss}(0,t)=\mbox{\mathversion{bold}$0$}, & t>0.
\end{array} \right.
\label{neg}
\end{eqnarray}
For \( \alpha >0\),
\begin{eqnarray}
\left\{
\begin{array}{ll}
\mbox{\mathversion{bold}$x$}_{t}= \mbox{\mathversion{bold}$x$}_{s}\times \mbox{\mathversion{bold}$x$}_{ss} + \alpha \big\{ \mbox{\mathversion{bold}$x$}_{sss}
+ \frac{3}{2}\mbox{\mathversion{bold}$x$}_{ss}\times \big( \mbox{\mathversion{bold}$x$}_{s}\times \mbox{\mathversion{bold}$x$}_{ss}\big) \big\}, &
s>0,t>0,\\
\mbox{\mathversion{bold}$x$}(s,0)= \mbox{\mathversion{bold}$x$}_{0}(s), & s>0, \\
\mbox{\mathversion{bold}$x$}_{s}(0,t)= \mbox{\mathversion{bold}$e$}_{3}, & t>0, \\
\mbox{\mathversion{bold}$x$}_{ss}(0,t)= \mbox{\mathversion{bold}$0$}, & t>0.
\end{array} \right.
\label{pos}
\end{eqnarray}
Here, \( \mbox{\mathversion{bold}$x$}(s,t)= (x^{1}(s,t), x^{2}(s,t), x^{3}(s,t))\) is the 
position vector of the vortex filament parameterized by its arc length \( s\) at time \( t\), \( \times \) is the exterior product in the three dimensional Euclidean space,
\( \alpha \) is a non-zero constant that describes the magnitude of the effect of axial flow, \( \mbox{\mathversion{bold}$e$}_{3}=(0,0,1)\), and
subscripts denote derivatives with their respective variables. Later in this paper, we will also use \( \partial_{s}\) and \( \partial _{t}\) to denote partial derivatives as well. We will refer to the equation in (\ref{neg}) and (\ref{pos}) as the vortex filament equation.
We note here that the number of boundary conditions imposed changes depending on the sign of \( \alpha \).
This is because the number of characteristic roots with a negative real part of the linearized equation, \( \mbox{\mathversion{bold}$x$}_{t}=\alpha \mbox{\mathversion{bold}$x$}_{sss}\),
changes depending on the sign of \( \alpha \). 

Our motivation for considering (\ref{neg}) and (\ref{pos}) comes from analyzing the motion of a tornado. 
This paper is our humble attempt to model the motion of a tornado. While it is obvious that a vortex filament is not the same as a tornado
and such modeling is questionable, 
many aspects of tornadoes are still unknown and we hope that our research can serve as a small step towards the 
complete analysis of the motion of a tornado.

To this end, in an earlier paper \cite{11}, the authors considered an initial-boundary value problem for the vortex filament equation with \( \alpha =0\), which is called the 
Localized Induction Equation (LIE). The LIE is a simplified model equation describing the motion of a vortex filament without axial flow. 
Other results considering the LIE can be found in Nishiyama and Tani \cite{5} and Koiso \cite{12}.

Many results are known for the Cauchy problem for the vortex filament equation with non-zero \( \alpha \), where the filament extends to spacial infinity
or the filament is closed. For example, in Nishiyama and Tani \cite{5}, they proved the unique solvability 
globally in time in Sobolev spaces. Onodera \cite{6,7} proved the unique solvability for a geometrically generalized equation. 
Segata \cite{8} proved the unique solvability and showed the asymptotic behavior in time of the solution to the 
Hirota equation, given by
\begin{eqnarray}
{\rm i}q_{t}=q_{xx}+\frac{1}{2}|q|^{2}q + {\rm i}\alpha \big( q_{xxx} + |q|^{2}q_{x}\big),
\label{hiro}
\end{eqnarray}
which can be obtained by applying the generalized Hasimoto transformation to the vortex filament equation. 
Since there are many results regarding the Cauchy problem for the Hirota equation and other Schr\"odinger type equations, it may feel more natural to see if the available theories
from these results can be utilized to solve the initial-boundary value problem for (\ref{hiro}), instead of considering (\ref{neg}) and (\ref{pos}) directly. 
Admittedly, problem (\ref{neg}) and (\ref{pos}) can be transformed into an initial-boundary value problem for the Hirota equation. But,
in light of the possibility that a new boundary condition may be considered for the vortex filament equation in the future, we thought that
it would be helpful to develop the analysis of the vortex filament equation itself because the Hasimoto transformation may not be applicable depending on the 
new boundary condition. For example, (\ref{neg}) and (\ref{pos}) model a vortex filament moving in the three dimensional 
half space, but if we consider a boundary that is not flat, it is nontrivial as to if we can apply the Hasimoto transformation or not, 
so we decided to work with the vortex filament equation directly.

For convenience, we introduce a new variable \( \mbox{\mathversion{bold}$v$}(s,t):= \mbox{\mathversion{bold}$x$}_{s}(s,t)\) and rewrite the problems in 
terms of \( \mbox{\mathversion{bold}$v$}\). Setting \( \mbox{\mathversion{bold}$v$}_{0}(s):= \mbox{\mathversion{bold}$x$}_{0s}(s)\), we have for \( \alpha <0\),
\begin{eqnarray}
\left\{
\begin{array}{ll}
\mbox{\mathversion{bold}$v$}_{t}= \mbox{\mathversion{bold}$v$}\times \mbox{\mathversion{bold}$v$}_{ss} + \alpha \big\{ \mbox{\mathversion{bold}$v$}_{sss}
+ \frac{3}{2}\mbox{\mathversion{bold}$v$}_{ss}\times \big( \mbox{\mathversion{bold}$v$}\times \mbox{\mathversion{bold}$v$}_{s}\big) 
+ \frac{3}{2}\mbox{\mathversion{bold}$v$}_{s}\times \big( \mbox{\mathversion{bold}$v$}\times \mbox{\mathversion{bold}$v$}_{ss}\big) \big\}, & s>0,t>0,\\
\mbox{\mathversion{bold}$v$}(s,0) = \mbox{\mathversion{bold}$v$}_{0}(s), & s>0,\\
\mbox{\mathversion{bold}$v$}_{s}(0,t)=\mbox{\mathversion{bold}$0$}, & t>0.
\end{array}\right.
\label{negv}
\end{eqnarray}
For \( \alpha >0\),
\begin{eqnarray}
\left\{
\begin{array}{ll}
\mbox{\mathversion{bold}$v$}_{t}= \mbox{\mathversion{bold}$v$}\times \mbox{\mathversion{bold}$v$}_{ss} + \alpha \big\{ \mbox{\mathversion{bold}$v$}_{sss}
+ \frac{3}{2}\mbox{\mathversion{bold}$v$}_{ss}\times \big( \mbox{\mathversion{bold}$v$}\times \mbox{\mathversion{bold}$v$}_{s}\big) 
+ \frac{3}{2}\mbox{\mathversion{bold}$v$}_{s}\times \big( \mbox{\mathversion{bold}$v$}\times \mbox{\mathversion{bold}$v$}_{ss}\big) \big\}, & s>0,t>0,\\
\mbox{\mathversion{bold}$v$}(s,0) = \mbox{\mathversion{bold}$v$}_{0}(s), & s>0,\\
\mbox{\mathversion{bold}$v$}(0,t) = \mbox{\mathversion{bold}$e$}_{3}, & t>0,\\
\mbox{\mathversion{bold}$v$}_{s}(0,t) = \mbox{\mathversion{bold}$0$}, & t>0.
\end{array}\right.
\label{posv}
\end{eqnarray}
Once we obtain a solution for (\ref{negv}) and (\ref{posv}), we can reconstruct \( \mbox{\mathversion{bold}$x$}(s,t)\) from the formula
\begin{eqnarray*}
\mbox{\mathversion{bold}$x$}(s,t)= \mbox{\mathversion{bold}$x$}_{0}(s) + \int^{t}_{0} \big\{ 
\mbox{\mathversion{bold}$v$}\times \mbox{\mathversion{bold}$v$}_{s} 
+ \alpha \mbox{\mathversion{bold}$v$}_{ss}
+ \frac{3}{2}\alpha \mbox{\mathversion{bold}$v$}_{s}\times \big( \mbox{\mathversion{bold}$v$}\times \mbox{\mathversion{bold}$v$}_{s}\big) \big\}(s, \tau) {\rm d}\tau,
\end{eqnarray*}
and \( \mbox{\mathversion{bold}$x$}(s,t)\) will satisfy (\ref{neg}) and (\ref{pos}) respectively, in other words, (\ref{neg}) is equivalent to (\ref{negv}) and 
(\ref{pos}) is equivalent to (\ref{posv}). Hence, we will concentrate on the solvability of (\ref{negv}) and (\ref{posv}) from now on.
Our approach for solving (\ref{negv}) and (\ref{posv}) is to consider the associated linear problem. 
Linearizing the equation around a function \( \mbox{\mathversion{bold}$w$}\) and neglecting lower order terms yield
\begin{eqnarray*}
\mbox{\mathversion{bold}$v$}_{t}= \mbox{\mathversion{bold}$w$}\times \mbox{\mathversion{bold}$v$}_{ss}
+ \alpha \big\{ \mbox{\mathversion{bold}$v$}_{sss} + \frac{3}{2}\mbox{\mathversion{bold}$v$}_{ss}\times \big(
\mbox{\mathversion{bold}$w$}\times \mbox{\mathversion{bold}$w$}_{s}\big) 
+ \frac{3}{2}\mbox{\mathversion{bold}$w$}_{s}\times \big( \mbox{\mathversion{bold}$w$}\times 
\mbox{\mathversion{bold}$v$}_{ss}\big)\big\}.
\end{eqnarray*}
Directly considering the initial-boundary value problem for the above equation seems hard.
When we try to estimate the solution in Sobolev spaces, the 
term \( \mbox{\mathversion{bold}$w$}_{s}\times \big( \mbox{\mathversion{bold}$w$}\times 
\mbox{\mathversion{bold}$v$}_{ss}\big) \) causes a loss of regularity because of the form of the
coefficient. We were able to overcome this by using the fact that if the initial datum is parameterized by its arc length, 
i.e. \( | \mbox{\mathversion{bold}$v$}_{0}|=1\), a sufficiently smooth solution 
of (\ref{negv}) and (\ref{posv}) satisfies \( |\mbox{\mathversion{bold}$v$}|=1\), and this allows us to make the transformation
\begin{eqnarray*}
\mbox{\mathversion{bold}$v$}_{s}\times \big( \mbox{\mathversion{bold}$v$}\times 
\mbox{\mathversion{bold}$v$}_{ss}\big) = \mbox{\mathversion{bold}$v$}_{ss}\times
\big( \mbox{\mathversion{bold}$v$}\times \mbox{\mathversion{bold}$v$}_{s}\big) 
-|\mbox{\mathversion{bold}$v$}_{s}|^{2}\mbox{\mathversion{bold}$v$}_{s}.
\end{eqnarray*}
Linearizing the equation in (\ref{negv}) and (\ref{posv}) after the above transformation yields
\begin{eqnarray}
\mbox{\mathversion{bold}$v$}_{t}= \mbox{\mathversion{bold}$w$}\times \mbox{\mathversion{bold}$v$}_{ss}
+ \alpha \big\{ \mbox{\mathversion{bold}$v$}_{sss} + 3\mbox{\mathversion{bold}$v$}_{ss}\times \big(
\mbox{\mathversion{bold}$w$}\times \mbox{\mathversion{bold}$w$}_{s}\big) \big\}.
\label{lv}
\end{eqnarray}
The term that was causing the loss of regularity is gone, but still, the existence of a solution 
to the initial-boundary value problem of the above third order dispersive equation is not trivial. 

One may wonder if we could treat the second order derivative terms as a 
perturbation of the linear KdV or the KdV-Burgers equation to avoid the above difficulties all together. 
This seems impossible, because as far as the authors know, the estimates obtained for the 
linear KdV and KdV-Burgers equation is insufficient to consider a second order term as a 
regular perturbation. See, for example, Hayashi and Kaikina \cite{3}, Hayashi, Kaikina, and Ruiz Paredes \cite{4}, 
or Bona and Zhang \cite{9} for known results on the initial-boundary value problems for the 
KdV and KdV-Burgers equations. 
This was our motivation in a previous paper \cite{10} to consider initial-boundary value problems for equations of the form
\begin{eqnarray}
\mbox{\mathversion{bold}$u$}_{t}= \alpha \mbox{\mathversion{bold}$u$}_{xxx}+{\rm A}(\mbox{\mathversion{bold}$w$},\partial_{x})\mbox{\mathversion{bold}$u$} + 
\mbox{\mathversion{bold}$f$},
\label{lpd}
\end{eqnarray}
where \( \mbox{\mathversion{bold}$u$}(x,t)=(u^{1}(x,t), u^{2}(x,t), \ldots  ,u^{m}(x,t)) \) is the unknown vector valued function, 
\( \mbox{\mathversion{bold}$w$}(x,t)
= ( w^{1}(x,t),w^{2}(x,t), \ldots ,w^{k}(x,t))\) and \( \mbox{\mathversion{bold}$f$}(x,t)=( f^{1}(x,t), f^{2}(x,t), \ldots ,f^{m}(x,t))\) are known vector valued 
functions, and \( {\rm A}(\mbox{\mathversion{bold}$w$},\partial_{x})\) is 
a second order differential operator of the form \( {\rm A}({\mbox{\mathversion{bold}$w$}},\partial_{x})={\rm A}_{0}({\mbox{\mathversion{bold}$w$}})\partial ^{2}_{x} 
+ {\rm A}_{1}({\mbox{\mathversion{bold}$w$}})\partial _{x} + {\rm A}_{2}({\mbox{\mathversion{bold}$w$}}) \). \( {\rm A}_{0}, \ {\rm A}_{1}, \ {\rm A}_{2}\) are
smooth matrices and \( {\rm A}({\mbox{\mathversion{bold}$w$}},\partial_{x}) \) is strongly elliptic in the sense that for any bounded domain \( E\) in 
\( \mathbf{R}^{k}\), there is a positive constant \( \delta \) such that for any \( \mbox{\mathversion{bold}$w$}\in E\)
\begin{eqnarray*}
{\rm A}_{0}(\mbox{\mathversion{bold}$w$})
+{\rm A}_{0}(\mbox{\mathversion{bold}$w$})^{*}\geq \delta {\rm I},
\end{eqnarray*}
where \( {\rm I}\) is the unit matrix and \( *\) denotes the adjoint of a matrix.
We proved the unique solvability of an initial-boundary value problem of the above equation 
in Sobolev spaces, and the precise statement we proved will be addressed in a later section. 
This result can be applied to (\ref{lv}) after we regularize it with a second order viscosity term \( \delta \mbox{\mathversion{bold}$v$}_{ss}\) with \( \delta >0\). 

The contents of this paper are as follows. In section 2, we introduce function spaces and the associated notations.
We also state the main results in this paper.
In section 3, we consider the compatibility conditions for
regularized nonlinear problems and the necessary correction of the initial datum. In section 4, we
review the existence theorem for associated linear problems which will be applied to the nonlinear problems.
In section 5, we prove an existence theorem for the case \( \alpha <0\), and in 
section 6, we prove an existence theorem for the case \( \alpha >0\).


\section{Function Spaces, Notations, and Main Theorems}

We define some function spaces that will be used throughout this paper, and notations associated with the spaces.
For a non-negative integer \( m\), and \( 1\leq p \leq \infty \), \( W^{m,p}(\mathbf{R}_{+} )\) is the Sobolev space 
containing all real-valued functions that have derivatives in the sense of distribution up to order \( m\) 
belonging to \( L^{p}(\mathbf{R}_{+} )\).
We set \( H^{m}(\mathbf{R}_{+} ) := W^{m,2}(\mathbf{R}_{+} ) \) as the Sobolev space equipped with the usual inner product.
The norm in \( H^{m}(\mathbf{R}_{+} ) \) is denoted by \( \| \cdot \|_{m} \) and we simply write \( \| \cdot \| \) 
for \( \|\cdot \|_{0} \). Otherwise, for a Banach space \( X\), the norm in \( X\) is written as \( \| \cdot \| _{X}\).
The inner product in \( L^{2}(\mathbf{R}_{+})\) is denoted by \( (\cdot ,\cdot )\).

For \( 0<T< \infty \) and a Banach space \( X\), 
\( C^{m}([0,T];X) \) denotes the space of functions that are \( m\) times continuously differentiable 
in \( t\) with respect to the norm of \( X\), and 
\( L^{2}(0,T;X)\) is the space of functions with
\(
\int ^{T}_{0} \| u(t)\|^{2}_{X}{\rm d}t
\)
being finite.

For any function space described above, we say that a vector valued function belongs to the function space 
if each of its components does. 

Finally, we state our main existence theorems here.
\begin{Th}{\rm (The case \( \alpha >0\) )}	
For a natural number \( k\), if \( \mbox{\mathversion{bold}$x$}_{0ss}
\in H^{2+3k}(\mathbf{R}_{+})\), \( |\mbox{\mathversion{bold}$x$}_{0s}|=1\), 
and \( \mbox{\mathversion{bold}$x$}_{0s}\) satisfies the compatibility conditions
for {\rm (\ref{pos2})} up to order \( k\), then there exists \( T>0\) such that 
{\rm (\ref{pos})} has a unique solution \( \mbox{\mathversion{bold}$x$}\) satisfying
\begin{eqnarray*}
\mbox{\mathversion{bold}$x$}_{ss}\in \bigcap ^{k}_{j=0}
W^{j,\infty}\big( [0,T]; H^{2+3j}(\mathbf{R}_{+})\big)
\end{eqnarray*}
and \( |\mbox{\mathversion{bold}$x$}_{s}|=1\). Here, \( T\) depends on \( \| \mbox{\mathversion{bold}$x$}_{0ss}\|_{2}\).
\label{st1}
\end{Th}
\begin{Th}{\rm (The case \( \alpha <0\) )}
For a natural number \( k\), if \( \mbox{\mathversion{bold}$x$}_{0ss}
\in H^{1+3k}(\mathbf{R}_{+})\), \( |\mbox{\mathversion{bold}$x$}_{0s}|=1\), 
and \( \mbox{\mathversion{bold}$x$}_{0s}\) satisfies the compatibility conditions
for {\rm (\ref{dneg2})} up to order \( k\), then there exists \( T>0\) such that 
{\rm (\ref{neg})} has a unique solution \( \mbox{\mathversion{bold}$x$}\) satisfying
\begin{eqnarray*}
\mbox{\mathversion{bold}$x$}_{ss}\in \bigcap ^{k}_{j=0}
W^{j,\infty}\big( [0,T]; H^{1+3j}(\mathbf{R}_{+})\big)
\end{eqnarray*}
and \( |\mbox{\mathversion{bold}$x$}_{s}|=1\). Here, \( T\) depends on \( \|\mbox{\mathversion{bold}$x$}_{0ss}\|_{3}\). 
\label{st2}
\end{Th}
%

%
%

%
%
%
%
%
\section{Regularized Nonlinear Problem and its Compatibility Conditions}
\setcounter{equation}{0}
We construct the solution of 
\begin{eqnarray}
\left\{
\begin{array}{ll}
\mbox{\mathversion{bold}$v$}_{t}=\mbox{\mathversion{bold}$v$}\times \mbox{\mathversion{bold}$v$}_{ss} 
+ \alpha \big\{ \mbox{\mathversion{bold}$v$}_{sss} + 3\mbox{\mathversion{bold}$v$}_{ss} \times 
\big( \mbox{\mathversion{bold}$v$}\times \mbox{\mathversion{bold}$v$}_{s}\big)
- \frac{3}{2}|\mbox{\mathversion{bold}$v$}_{s}|^{2}\mbox{\mathversion{bold}$v$}_{s}\big\}, & s>0, t>0,\\
\mbox{\mathversion{bold}$v$}(s,0)=\mbox{\mathversion{bold}$v$}_{0}(s), & s>0,\\
\mbox{\mathversion{bold}$v$}_{s}(0,t)=\mbox{\mathversion{bold}$0$}, & t>0,
\end{array}\right.
\label{neg2}
\end{eqnarray}
and
\begin{eqnarray}
\left\{
\begin{array}{ll}
\mbox{\mathversion{bold}$v$}_{t}=\mbox{\mathversion{bold}$v$}\times \mbox{\mathversion{bold}$v$}_{ss} 
+ \alpha \big\{ \mbox{\mathversion{bold}$v$}_{sss} + 3\mbox{\mathversion{bold}$v$}_{ss} \times 
\big( \mbox{\mathversion{bold}$v$}\times \mbox{\mathversion{bold}$v$}_{s}\big)
- \frac{3}{2}|\mbox{\mathversion{bold}$v$}_{s}|^{2}\mbox{\mathversion{bold}$v$}_{s}\big\}, & s>0, t>0,\\
\mbox{\mathversion{bold}$v$}(s,0)=\mbox{\mathversion{bold}$v$}_{0}(s), & s>0,\\
\mbox{\mathversion{bold}$v$}(0,t)=\mbox{\mathversion{bold}$e$}_{3}, & t>0,\\
\mbox{\mathversion{bold}$v$}_{s}(0,t)=\mbox{\mathversion{bold}$0$}, & t>0,
\end{array}\right.
\label{pos2}
\end{eqnarray}
which are equivalent to (\ref{neg}) and (\ref{pos}) respectively,
by taking the limit \( \delta \rightarrow +0\) in 
the following regularized problems.
\begin{eqnarray}
\left\{
\begin{array}{ll}
\mbox{\mathversion{bold}$v$}^{\delta}_{t}=\mbox{\mathversion{bold}$v$}^{\delta}\times \mbox{\mathversion{bold}$v$}^{\delta}_{ss} 
+ \alpha \big\{ \mbox{\mathversion{bold}$v$}^{\delta}_{sss} + 3\mbox{\mathversion{bold}$v$}^{\delta}_{ss} \times 
\big( \mbox{\mathversion{bold}$v$}^{\delta}\times \mbox{\mathversion{bold}$v$}^{\delta}_{s}\big)
- \frac{3}{2}|\mbox{\mathversion{bold}$v$}^{\delta}_{s}|^{2}\mbox{\mathversion{bold}$v$}^{\delta}_{s}\big\}& \ \\
\hspace*{7cm}+ \delta \big( \mbox{\mathversion{bold}$v$}^{\delta}_{ss} + |\mbox{\mathversion{bold}$v$}^{\delta}_{s}|^{2}
\mbox{\mathversion{bold}$v$}^{\delta}\big) , & s>0, t>0,\\
\mbox{\mathversion{bold}$v$}^{\delta}(s,0)=\mbox{\mathversion{bold}$v$}^{\delta}_{0}(s), & s>0,\\
\mbox{\mathversion{bold}$v$}^{\delta}_{s}(0,t)=\mbox{\mathversion{bold}$0$}, & t>0,
\end{array}\right.
\label{dneg2}
\end{eqnarray}
and
\begin{eqnarray}
\left\{
\begin{array}{ll}
\mbox{\mathversion{bold}$v$}^{\delta}_{t}=\mbox{\mathversion{bold}$v$}^{\delta}\times \mbox{\mathversion{bold}$v$}^{\delta}_{ss} 
+ \alpha \big\{ \mbox{\mathversion{bold}$v$}^{\delta}_{sss} + 3\mbox{\mathversion{bold}$v$}^{\delta}_{ss} \times 
\big( \mbox{\mathversion{bold}$v$}^{\delta}\times \mbox{\mathversion{bold}$v$}^{\delta}_{s}\big)
- \frac{3}{2}|\mbox{\mathversion{bold}$v$}^{\delta}_{s}|^{2}\mbox{\mathversion{bold}$v$}^{\delta}_{s}\big\}& \ \\
\hspace*{7cm}+ \delta \big( \mbox{\mathversion{bold}$v$}^{\delta}_{ss} + |\mbox{\mathversion{bold}$v$}^{\delta}_{s}|^{2}
\mbox{\mathversion{bold}$v$}^{\delta}\big) , & s>0, t>0,\\
\mbox{\mathversion{bold}$v$}^{\delta}(s,0)=\mbox{\mathversion{bold}$v$}^{\delta}_{0}(s), & s>0,\\
\mbox{\mathversion{bold}$v$}^{\delta}(0,t)=\mbox{\mathversion{bold}$e$}_{3}, & t>0,\\
\mbox{\mathversion{bold}$v$}^{\delta}_{s}(0,t)=\mbox{\mathversion{bold}$0$}, & t>0.
\end{array}\right.
\label{dpos2}
\end{eqnarray}
From here on in, it is assumed that \( |\mbox{\mathversion{bold}$v$}_{0}|=1\) holds, 
i.e. the initial datum is parameterized by its arc length.
Since we modified the equation, we must make corrections to the initial datum to 
insure the compatibility conditions continue to hold for each problem.
\subsection{Compatibility Conditions for (\ref{neg2}) and (\ref{pos2})}
First, we derive the compatibility conditions for (\ref{neg2}) and (\ref{pos2}). 
We set \( \mbox{\mathversion{bold}$Q$}_{(0)}(\mbox{\mathversion{bold}$v$})=
\mbox{\mathversion{bold}$v$}\) and we denote 
the right-hand side of the equation in (\ref{neg2}) and (\ref{pos2}) as
\begin{eqnarray*}
\mbox{\mathversion{bold}$Q$}_{(1)}(\mbox{\mathversion{bold}$v$})=
\mbox{\mathversion{bold}$v$}\times \mbox{\mathversion{bold}$v$}_{ss} 
+ \alpha \big\{ \mbox{\mathversion{bold}$v$}_{sss} + 3\mbox{\mathversion{bold}$v$}_{ss} \times 
\big( \mbox{\mathversion{bold}$v$}\times \mbox{\mathversion{bold}$v$}_{s}\big)
- \frac{3}{2}|\mbox{\mathversion{bold}$v$}_{s}|^{2}\mbox{\mathversion{bold}$v$}_{s}\big\}.
\end{eqnarray*}
We will also use the notation \( \mbox{\mathversion{bold}$Q$}_{(1)}(s,t)\) and \(\mbox{\mathversion{bold}$Q$}_{(1)}\) instead of 
\( \mbox{\mathversion{bold}$Q$}_{(1)}(\mbox{\mathversion{bold}$v$})\) for convenience. 
For \( n\geq 2\), we successively define \( \mbox{\mathversion{bold}$Q$}_{(n)} \) by
\begin{eqnarray*}
\begin{aligned}
\mbox{\mathversion{bold}$Q$}_{(n)}&=
\sum ^{n-1}_{j=0} 
\left(
\begin{array}{c}
n-1\\
j
\end{array}\right)
\mbox{\mathversion{bold}$Q$}_{(j)}\times \mbox{\mathversion{bold}$Q$}_{(n-1-j)ss}
+ \alpha \mbox{\mathversion{bold}$Q$}_{(n-1)sss}\\
&+ 3\alpha  \bigg\{ \sum ^{n-1}_{j=0}
\sum ^{n-1-j}_{k=0}
\left(
\begin{array}{c}
n-1\\
j
\end{array}\right)
\left(
\begin{array}{c}
n-1-j\\
k
\end{array}\right)
\mbox{\mathversion{bold}$Q$}_{(j)ss}\times \bigg( 
\mbox{\mathversion{bold}$Q$}_{(k)}\times \mbox{\mathversion{bold}$Q$}_{(n-1-j-k)s}\bigg)\bigg\}\\
&-\frac{3}{2}\alpha \left\{\sum ^{n-1}_{j=0} \sum^{n-1-j}_{k=0}
\left(
\begin{array}{c}
n-1\\
j
\end{array}\right)
\left(
\begin{array}{c}
n-1-j\\
k
\end{array}\right)
\left( \mbox{\mathversion{bold}$Q$}_{(j)s}\cdot \mbox{\mathversion{bold}$Q$}_{(k)s}\right)
\mbox{\mathversion{bold}$Q$}_{(n-1-j-k)s}\right\}.
\end{aligned}
\end{eqnarray*}
The above definition of \( \mbox{\mathversion{bold}$Q$}_{(n)}(\mbox{\mathversion{bold}$v$})\) corresponds to giving an expression for \( \partial ^{n}_{t}\mbox{\mathversion{bold}$v$}\)
in terms of \( \mbox{\mathversion{bold}$v$}\) and its \( s\) derivatives only. 
It is obvious from the definition that the term with the highest order derivative in 
\( \mbox{\mathversion{bold}$Q$}_{(n)}\) is \( \alpha ^{n}\partial ^{3n}_{s}\mbox{\mathversion{bold}$v$}\). 
From the boundary conditions of (\ref{neg2}) and (\ref{pos2}), we arrive at the following compatibility conditions.
\begin{df}{\rm (}Compatibility conditions for {\rm (\ref{neg2}))}
For \( n\in \mathbf{N}\cup \{ 0\} \), we say that \( \mbox{\mathversion{bold}$v$}_{0}\) satisfies the \( n\)-th comaptibility condition for {\rm (\ref{neg2})}
if \( \mbox{\mathversion{bold}$v$}_{0s}\in H^{1+3n}(\mathbf{R}_{+}) \) and
\begin{eqnarray*}
\big( \partial_{s}\mbox{\mathversion{bold}$Q$}_{(n)}(\mbox{\mathversion{bold}$v$}_{0})\big)(0)=\mbox{\mathversion{bold}$0$}.
\end{eqnarray*}
We also say that \( \mbox{\mathversion{bold}$v$}_{0}\) satisfies the compatibility conditions 
for {\rm (\ref{neg2})} up to order \( n\) if it satisfies the \( k\)-th compatibility condition for all \( k\) with \( 0\leq k\leq n\).
\label{ccneg2}
\end{df}
\begin{df}{\rm (}Compatibility conditions for {\rm (\ref{pos2}))}
For \( n\in \mathbf{N}\cup \{ 0\} \), we say that \( \mbox{\mathversion{bold}$v$}_{0}\) satisfies the \( n\)-th comaptibility condition for {\rm (\ref{pos2})}
if \( \mbox{\mathversion{bold}$v$}_{0s}\in H^{2+3n}(\mathbf{R}_{+}) \) and 
\begin{eqnarray*}
\mbox{\mathversion{bold}$v$}_{0}(0)=\mbox{\mathversion{bold}$e$}_{3}, \ 
\mbox{\mathversion{bold}$v$}_{0s}(0)=\mbox{\mathversion{bold}$0$},
\end{eqnarray*}
when \( n=0\), and
\begin{eqnarray*}
(\mbox{\mathversion{bold}$Q$}_{(n)}(\mbox{\mathversion{bold}$v$}_{0}))(0)=\mbox{\mathversion{bold}$0$}, \ 
\big( \partial_{s}\mbox{\mathversion{bold}$Q$}_{(n)}(\mbox{\mathversion{bold}$v$}_{0})\big) (0)=\mbox{\mathversion{bold}$0$},
\end{eqnarray*}
when \( n\geq 1\). 
We also say that \( \mbox{\mathversion{bold}$v$}_{0}\) satisfies the compatibility conditions 
for {\rm (\ref{pos2})} up to order \( n\) if it satisfies the \( k\)-th compatibility condition for all \( k\) with \( 0\leq k\leq n\).
\label{ccpos2}
\end{df}
Note that the regularity imposed on \( \mbox{\mathversion{bold}$v$}_{0s}\) in 
Definition \ref{ccpos2} is not the minimal regularity required for the trace at 
\( s=0\) to have meaning, but we defined it as above so that it corresponds to the regularity assumption in the existence theorem 
that we obtain later. Also note that the regularity assumption is made on 
\( \mbox{\mathversion{bold}$v$}_{0s}\) instead of \( \mbox{\mathversion{bold}$v$}_{0}\)
because \( |\mbox{\mathversion{bold}$v$}_{0}|=1\) and so \( \mbox{\mathversion{bold}$v$}_{0}\) is not square integrable. 
\subsection{Compatibility Conditions for (\ref{dneg2}) and (\ref{dpos2})}
We derive the compatibility conditions for (\ref{dneg2}) and (\ref{dpos2}) in the 
same way as we derived the conditions for (\ref{neg2}) and (\ref{pos2}).
Set \( \mbox{\mathversion{bold}$P$}_{(0)}(\mbox{\mathversion{bold}$v$})=\mbox{\mathversion{bold}$v$}\) 
and define \( \mbox{\mathversion{bold}$P$}_{(1)}(\mbox{\mathversion{bold}$v$})\) by
\begin{eqnarray*}
\mbox{\mathversion{bold}$P$}_{(1)}(\mbox{\mathversion{bold}$v$})=
\mbox{\mathversion{bold}$v$}\times \mbox{\mathversion{bold}$v$}_{ss} 
+ \alpha \big\{ \mbox{\mathversion{bold}$v$}_{sss} + 3\mbox{\mathversion{bold}$v$}_{ss} \times 
\big( \mbox{\mathversion{bold}$v$}\times \mbox{\mathversion{bold}$v$}_{s}\big)
- \frac{3}{2}|\mbox{\mathversion{bold}$v$}_{s}|^{2}\mbox{\mathversion{bold}$v$}_{s}\big\}
+ \delta \big( \mbox{\mathversion{bold}$v$}_{ss} + |\mbox{\mathversion{bold}$v$}_{s}|^{2}\mbox{\mathversion{bold}$v$}\big).
\end{eqnarray*}
We successively define \( \mbox{\mathversion{bold}$P$}_{(n)}\) for \( n\geq 2\) by
\begin{eqnarray*}
\begin{aligned}
\mbox{\mathversion{bold}$P$}_{(n)}&=
\sum ^{n-1}_{j=0} 
\left(
\begin{array}{c}
n-1\\
j
\end{array}\right)
\mbox{\mathversion{bold}$P$}_{(j)}\times \mbox{\mathversion{bold}$P$}_{(n-1-j)ss}
+ \alpha \mbox{\mathversion{bold}$P$}_{(n-1)sss}\\
&+ 3\alpha  \bigg\{ \sum ^{n-1}_{j=0}
\sum ^{n-1-j}_{k=0}
\left(
\begin{array}{c}
n-1\\
j
\end{array}\right)
\left(
\begin{array}{c}
n-1-j\\
k
\end{array}\right)
\mbox{\mathversion{bold}$P$}_{(j)ss}\times \bigg( 
\mbox{\mathversion{bold}$P$}_{(k)}\times \mbox{\mathversion{bold}$P$}_{(n-1-j-k)s}\bigg)\bigg\}\\
&-\frac{3}{2}\alpha \left\{\sum ^{n-1}_{j=0} \sum^{n-1-j}_{k=0}
\left(
\begin{array}{c}
n-1\\
j
\end{array}\right)
\left(
\begin{array}{c}
n-1-j\\
k
\end{array}\right)
\left( \mbox{\mathversion{bold}$P$}_{(j)s}\cdot \mbox{\mathversion{bold}$P$}_{(k)s}\right)
\mbox{\mathversion{bold}$P$}_{(n-1-j-k)s}\right\} \\
& + \delta \left\{ \mbox{\mathversion{bold}$P$}_{(n-1)ss}
+ \sum ^{n-1}_{j=0} \sum^{n-1-j}_{k=0}
\left(
\begin{array}{c}
n-1\\
j
\end{array}\right)
\left(
\begin{array}{c}
n-1-j\\
k
\end{array}\right)
\left( \mbox{\mathversion{bold}$P$}_{(j)s}\cdot \mbox{\mathversion{bold}$P$}_{(k)s}\right)
\mbox{\mathversion{bold}$P$}_{(n-1-j-k)}\right\}.
\end{aligned}
\end{eqnarray*}
We arrive at the following compatibility conditions.
\begin{df}{\rm (}Compatibility conditions for {\rm (\ref{dneg2}))}
For \( n\in \mathbf{N}\cup \{ 0\} \), we say that \( \mbox{\mathversion{bold}$v$}^{\delta}_{0}\) satisfies the \( n\)-th comaptibility condition for {\rm (\ref{dneg2})}
if \( \mbox{\mathversion{bold}$v$}^{\delta}_{0s}\in H^{1+3n}(\mathbf{R}_{+}) \) and
\begin{eqnarray*}
\big( \partial_{s}\mbox{\mathversion{bold}$P$}_{(n)}(\mbox{\mathversion{bold}$v$}^{\delta}_{0})\big) (0)=\mbox{\mathversion{bold}$0$}.
\end{eqnarray*}
We also say that \( \mbox{\mathversion{bold}$v$}^{\delta}_{0}\) satisfies the compatibility conditions 
for {\rm (\ref{dneg2})} up to order \( n\) if it satisfies the \( k\)-th compatibility condition for all \( k\) with \( 0\leq k\leq n\).
\label{ccdneg2}
\end{df}
\begin{df}{\rm (}Compatibility conditions for {\rm (\ref{dpos2}))}
For \( n\in \mathbf{N}\cup \{ 0\} \), we say that \( \mbox{\mathversion{bold}$v$}^{\delta}_{0}\) satisfies the \( n\)-th comaptibility condition for {\rm (\ref{dpos2})}
if \( \mbox{\mathversion{bold}$v$}^{\delta}_{0s}\in H^{2+3n}(\mathbf{R}_{+}) \) and 
\begin{eqnarray*}
\mbox{\mathversion{bold}$v$}^{\delta}_{0}(0)=\mbox{\mathversion{bold}$e$}_{3}, \ 
\mbox{\mathversion{bold}$v$}^{\delta}_{0s}(0)=\mbox{\mathversion{bold}$0$},
\end{eqnarray*}
when \( n=0\), and
\begin{eqnarray*}
(\mbox{\mathversion{bold}$P$}_{(n)}(\mbox{\mathversion{bold}$v$}^{\delta}_{0}))(0)=\mbox{\mathversion{bold}$0$}, \ 
\big( \partial_{s}\mbox{\mathversion{bold}$P$}_{(n)}(\mbox{\mathversion{bold}$v$}^{\delta}_{0})\big) (0)=\mbox{\mathversion{bold}$0$},
\end{eqnarray*}
when \( n\geq 1\). 
We also say that \( \mbox{\mathversion{bold}$v$}^{\delta}_{0}\) satisfies the compatibility conditions 
for {\rm (\ref{dpos2})} up to order \( n\) if it satisfies the \( k\)-th compatibility condition for all \( k\) with \( 0\leq k\leq n\).
\label{ccdpos2}
\end{df}
\subsection{Corrections to the Initial Datum}
We construct a corrected initial datum \( \mbox{\mathversion{bold}$v$}^{\delta}_{0}\)
such that given an initial datum \( \mbox{\mathversion{bold}$v$}_{0}\) that satisfies
the compatibility conditions for (\ref{neg2}) or (\ref{pos2}), 
\( \mbox{\mathversion{bold}$v$}^{\delta}_{0}\) satisfies the compatibility conditions
of (\ref{dneg2}) and (\ref{dpos2}) respectively, and 
\( \mbox{\mathversion{bold}$v$}^{\delta}_{0}\rightarrow \mbox{\mathversion{bold}$v$}_{0}\)
in the appropriate function space. As it will be shown later, a sufficiently smooth
solution of (\ref{dneg2}) or (\ref{dpos2}) with \( \delta \geq 0\) satisfies 
\( |\mbox{\mathversion{bold}$v$}^{\delta}|=1\) if \( |\mbox{\mathversion{bold}$v$}^{\delta}_{0}|=1\).
Thus, the correction of the initial datum must be done in a way that preserves this 
property. Since the argument for the construction of \( \mbox{\mathversion{bold}$v$}^{\delta}_{0}\) is the same 
for the cases \( \alpha >0\) and \( \alpha <0\), we show the details for the case \( \alpha <0\) only.

Suppose that we have an initial datum \( \mbox{\mathversion{bold}$v$}_{0}\) such that
\( \mbox{\mathversion{bold}$v$}_{0s}\in H^{1+3m}(\mathbf{R}_{+})\) satisfying the compatibility conditions for (\ref{neg2}) up to order \( m\). 
We will construct \( \mbox{\mathversion{bold}$v$}^{\delta}_{0}\) in the form
\begin{eqnarray}
\mbox{\mathversion{bold}$v$}^{\delta}_{0} = \frac{\mbox{\mathversion{bold}$v$}_{0} + \mbox{\mathversion{bold}$h$}_{\delta } }
{|\mbox{\mathversion{bold}$v$}_{0} + \mbox{\mathversion{bold}$h$}_{\delta}|},
\label{cor}
\end{eqnarray}
where \( \mbox{\mathversion{bold}$h$}_{\delta}\) is constructed so that \( \mbox{\mathversion{bold}$h$}_{\delta}\rightarrow \mbox{\mathversion{bold}$0$}\)
as \( \delta \rightarrow +0\). The method we use to construct \( \mbox{\mathversion{bold}$h$}_{\delta}\) is standard, i.e. we substitute (\ref{cor}) into the 
compatibility conditions for (\ref{dneg2}) to determine the differential coefficients of \( \mbox{\mathversion{bold}$h$}_{\delta}\) at \( s=0\) and then extend it to \( s>0\) 
so that \( \mbox{\mathversion{bold}$h$}_{\delta}\) belongs to the appropriate Sobolev space and its differential coefficients have the desired value.

We introduce some notations. We set
\begin{eqnarray*}
\begin{aligned}
\mbox{\mathversion{bold}$g$}^{\delta}_{0}(\mbox{\mathversion{bold}$V$}) &:= \mbox{\mathversion{bold}$V$},\\
\mbox{\mathversion{bold}$g$}^{\delta}_{1}(\mbox{\mathversion{bold}$V$})&:= \mbox{\mathversion{bold}$V$}\times \mbox{\mathversion{bold}$V$}_{ss}
+ \alpha \{ \mbox{\mathversion{bold}$V$}_{sss} + 3\mbox{\mathversion{bold}$V$}_{ss}\times (\mbox{\mathversion{bold}$V$}\times \mbox{\mathversion{bold}$V$}_{s})
-\frac{3}{2}|\mbox{\mathversion{bold}$V$}_{s}|^{2}\mbox{\mathversion{bold}$V$}_{s}\} + \delta ( \mbox{\mathversion{bold}$V$}_{ss} + |\mbox{\mathversion{bold}$V$}_{s}|^{2}
\mbox{\mathversion{bold}$V$}), \\
\mbox{\mathversion{bold}$g$}^{\delta}_{m+1}(\mbox{\mathversion{bold}$V$})&:= 
D\mbox{\mathversion{bold}$g$}^{\delta }_{m}(\mbox{\mathversion{bold}$V$})[\mbox{\mathversion{bold}$g$}^{\delta }_{1}(\mbox{\mathversion{bold}$V$})],
\end{aligned}
\end{eqnarray*}
where \( m\geq 1\) and \( D\) is the derivative with respect to \( \mbox{\mathversion{bold}$V$}\), i.e. 
\( D\mbox{\mathversion{bold}$g$}^{\delta}_{m}(\mbox{\mathversion{bold}$V$})[\mbox{\mathversion{bold}$W$}]= \frac{{\rm d}}{{\rm d}\varepsilon }
\left. \mbox{\mathversion{bold}$g$}^{\delta}_{m}(\mbox{\mathversion{bold}$V$}+\varepsilon \mbox{\mathversion{bold}$W$})\right|_{\varepsilon =0}\).
Note that under these notations, the \( m\)-th order compatibility condition
for (\ref{dneg2}) can be expressed as \( \mbox{\mathversion{bold}$g$}^{\delta}_{m}(\mbox{\mathversion{bold}$v$}^{\delta}_{0})_{s}(0)=\mbox{\mathversion{bold}$0$}\),
because \( \mbox{\mathversion{bold}$P$}_{(m)}(\mbox{\mathversion{bold}$V$})=\mbox{\mathversion{bold}$g$}^{\delta}_{m}(\mbox{\mathversion{bold}$V$})\).
We gave a different notation because it is more convenient for the upcoming calculations.

First we prove that if \( |\mbox{\mathversion{bold}$V$}|\equiv 1\), then 
for any \( m\geq 1\) 
\begin{eqnarray}
\sum ^{m}_{k=0}
\left(
\begin{array}{c}
m\\
k
\end{array}\right)
\mbox{\mathversion{bold}$g$}^{\delta}_{k}(\mbox{\mathversion{bold}$V$})\cdot \mbox{\mathversion{bold}$g$}^{\delta}_{m-k}(\mbox{\mathversion{bold}$V$})\equiv 0.
\label{inner}
\end{eqnarray}
We show this by induction. From direct calculation, we can prove that
\begin{eqnarray*}
\mbox{\mathversion{bold}$g$}^{\delta}_{1}(\mbox{\mathversion{bold}$V$})\cdot \mbox{\mathversion{bold}$V$} 
&=\frac{\alpha }{2}(|\mbox{\mathversion{bold}$V$}|^{2})_{sss} - 3\alpha (\mbox{\mathversion{bold}$V$}\cdot \mbox{\mathversion{bold}$V$}_{ss})
(|\mbox{\mathversion{bold}$V$}|^{2})_{s}-\frac{3}{2}|\mbox{\mathversion{bold}$V$}_{s}|^{2}(|\mbox{\mathversion{bold}$V$}|^{2})_{s}
+ \frac{\delta }{2}(|\mbox{\mathversion{bold}$V$}|^{2})_{ss}=0,
\end{eqnarray*}
which proves (\ref{inner}) with \( m=1\). Suppose that (\ref{inner}) holds up to some \( m\) with \( m\geq 1\). From the assumption of 
induction we have for any vector \( \mbox{\mathversion{bold}$W$}\) and \( t\in \mathbf{R}\),
\begin{eqnarray*}
\sum ^{m}_{k=0}
\left(
\begin{array}{c}
m\\
k
\end{array}\right)
\mbox{\mathversion{bold}$g$}^{\delta}_{k}( \frac{\mbox{\mathversion{bold}$V$}+ t\mbox{\mathversion{bold}$W$}}{|\mbox{\mathversion{bold}$V$}+ t\mbox{\mathversion{bold}$W$}|})
\cdot \mbox{\mathversion{bold}$g$}^{\delta}_{m-k}(\frac{\mbox{\mathversion{bold}$V$}+ t\mbox{\mathversion{bold}$W$}}{|\mbox{\mathversion{bold}$V$}+ t\mbox{\mathversion{bold}$W$}|})
\equiv 0.
\end{eqnarray*}
Taking the \( t\) derivative and setting \( t=0\) yield
\begin{eqnarray*}
\sum ^{m}_{k=0}
\left(
\begin{array}{c}
m\\
k
\end{array}\right)
\left\{
D\mbox{\mathversion{bold}$g$}^{\delta}_{k}(\mbox{\mathversion{bold}$V$})[\mbox{\mathversion{bold}$W$}-(\mbox{\mathversion{bold}$V$}\cdot \mbox{\mathversion{bold}$W$})
\mbox{\mathversion{bold}$V$}]\cdot \mbox{\mathversion{bold}$g$}^{\delta}_{m-k}(\mbox{\mathversion{bold}$V$})
+\mbox{\mathversion{bold}$g$}^{\delta}_{k}(\mbox{\mathversion{bold}$V$})\cdot D \mbox{\mathversion{bold}$g$}^{\delta}_{m-k}(\mbox{\mathversion{bold}$V$})
[\mbox{\mathversion{bold}$W$}-(\mbox{\mathversion{bold}$V$}\cdot \mbox{\mathversion{bold}$W$})\mbox{\mathversion{bold}$V$}]\right\} \equiv 0.
\end{eqnarray*}
By choosing \( \mbox{\mathversion{bold}$W$}= \mbox{\mathversion{bold}$g$}^{\delta}_{1}(\mbox{\mathversion{bold}$V$}) \) we have
\begin{eqnarray*}
\begin{aligned}
0 &\equiv \sum ^{m}_{k=0}
\left(
\begin{array}{c}
m\\
k
\end{array}\right)\left\{ 
\mbox{\mathversion{bold}$g$}^{\delta}_{k+1}(\mbox{\mathversion{bold}$V$})\cdot \mbox{\mathversion{bold}$g$}^{\delta}_{m-k}(\mbox{\mathversion{bold}$V$})
+ \mbox{\mathversion{bold}$g$}^{\delta}_{k}(\mbox{\mathversion{bold}$V$})\cdot \mbox{\mathversion{bold}$g$}^{\delta}_{m-k+1}(\mbox{\mathversion{bold}$V$})\right\} \\
&= \sum ^{m+1}_{k=0}
\left( 
\begin{array}{c}
m+1\\
k
\end{array}\right) \mbox{\mathversion{bold}$g$}^{\delta}_{k}(\mbox{\mathversion{bold}$V$})\cdot \mbox{\mathversion{bold}$g$}^{\delta}_{m+1-k}(\mbox{\mathversion{bold}$V$}),
\end{aligned}
\end{eqnarray*}
which proves (\ref{inner}) for the case \( m+1\), and this finishes the proof.

Next we make the following notations.
\begin{eqnarray*}
\begin{aligned}
\mbox{\mathversion{bold}$f$}_{0}(\mbox{\mathversion{bold}$V$})&:=\mbox{\mathversion{bold}$V$},\\
\mbox{\mathversion{bold}$f$}_{1}(\mbox{\mathversion{bold}$V$})&:= \mbox{\mathversion{bold}$V$}\times \mbox{\mathversion{bold}$V$}_{ss}
+ \alpha \{ \mbox{\mathversion{bold}$V$}_{sss} + 3\mbox{\mathversion{bold}$V$}_{ss}\times (\mbox{\mathversion{bold}$V$}\times \mbox{\mathversion{bold}$V$}_{s})
-\frac{3}{2}|\mbox{\mathversion{bold}$V$}_{s}|^{2}\mbox{\mathversion{bold}$V$}_{s}\}, \\
\mbox{\mathversion{bold}$f$}_{m+1}(\mbox{\mathversion{bold}$V$})&:= D\mbox{\mathversion{bold}$f$}_{m}(\mbox{\mathversion{bold}$V$})
[\mbox{\mathversion{bold}$f$}_{1}(\mbox{\mathversion{bold}$V$})],
\end{aligned}
\end{eqnarray*}
which is equivalent to taking \( \delta =0\) in \( \mbox{\mathversion{bold}$g$}^{\delta}_{m}(\mbox{\mathversion{bold}$V$})\), so that 
\( \sum ^{m}_{k=0}
\left(
\begin{array}{c}
m\\
k
\end{array}\right)
\mbox{\mathversion{bold}$f$}_{k}(\mbox{\mathversion{bold}$V$})\cdot \mbox{\mathversion{bold}$f$}_{m-k}(\mbox{\mathversion{bold}$V$})\equiv 0\) if \( |\mbox{\mathversion{bold}$V$}|=1\), 
and the \( m\)-th order
compatibility condition for (\ref{neg2}) can be expressed as \( \mbox{\mathversion{bold}$f$}_{m}(\mbox{\mathversion{bold}$v$}_{0})_{s}(0) = \mbox{\mathversion{bold}$0$}\) 
because \( \mbox{\mathversion{bold}$Q$}_{(m)}(\mbox{\mathversion{bold}$v$}_{0}) = \mbox{\mathversion{bold}$f$}_{m}(\mbox{\mathversion{bold}$v$}_{0})\).

Next, we show that 
\begin{eqnarray}
\mbox{\mathversion{bold}$g$}^{\delta}_{m}(\mbox{\mathversion{bold}$V$}) = \mbox{\mathversion{bold}$f$}_{m}(\mbox{\mathversion{bold}$V$})
+ \delta \mbox{\mathversion{bold}$r$}^{\delta}_{m}(\mbox{\mathversion{bold}$V$}),
\label{remain}
\end{eqnarray}
where \( \mbox{\mathversion{bold}$r$}^{\delta}_{1}(\mbox{\mathversion{bold}$V$}) := \mbox{\mathversion{bold}$V$}_{ss} + |\mbox{\mathversion{bold}$V$}_{s}|^{2}
\mbox{\mathversion{bold}$V$}\) and \( \mbox{\mathversion{bold}$r$}^{\delta}_{m}(\mbox{\mathversion{bold}$V$}):=
D\mbox{\mathversion{bold}$r$}^{\delta}_{m-1}(\mbox{\mathversion{bold}$V$})[\mbox{\mathversion{bold}$g$}^{\delta}_{1}(\mbox{\mathversion{bold}$V$})]
+ D\mbox{\mathversion{bold}$f$}_{m-1}(\mbox{\mathversion{bold}$V$})[\mbox{\mathversion{bold}$r$}^{\delta}_{1}(\mbox{\mathversion{bold}$V$})]\) for 
\( m\geq 2\).
From the definition, \( \mbox{\mathversion{bold}$r$}^{\delta}_{m}(\mbox{\mathversion{bold}$V$})\) contains derivatives up to order \( 3m-1\).

It is obvious that (\ref{remain}) holds for \( m=1\) from the definition of \( \mbox{\mathversion{bold}$g$}^{\delta}_{1}\) and \( \mbox{\mathversion{bold}$f$}_{1}\).
Suppose that it holds up to \( m-1\) for some \( m\geq 2\). Thus, for any vector \( \mbox{\mathversion{bold}$W$}\) and \( t\in \mathbf{R}\), we have
\begin{eqnarray*}
\mbox{\mathversion{bold}$g$}^{\delta}_{m-1}(\mbox{\mathversion{bold}$V$}+t\mbox{\mathversion{bold}$W$})
= \mbox{\mathversion{bold}$f$}_{m-1}(\mbox{\mathversion{bold}$V$}+t\mbox{\mathversion{bold}$W$}) + \delta \mbox{\mathversion{bold}$r$}^{\delta}_{m-1}
(\mbox{\mathversion{bold}$V$} + t\mbox{\mathversion{bold}$W$}).
\end{eqnarray*}
Taking the \( t\) derivative of the above equation and setting \( t=0\) yields
\begin{eqnarray*}
D\mbox{\mathversion{bold}$g$}^{\delta}_{m-1}(\mbox{\mathversion{bold}$V$})[\mbox{\mathversion{bold}$W$}]
=D\mbox{\mathversion{bold}$f$}_{m-1}(\mbox{\mathversion{bold}$V$})[\mbox{\mathversion{bold}$W$}]
+ \delta D\mbox{\mathversion{bold}$r$}^{\delta}_{m-1}(\mbox{\mathversion{bold}$V$})[\mbox{\mathversion{bold}$W$}].
\end{eqnarray*}
Finally, choosing \( \mbox{\mathversion{bold}$W$}= \mbox{\mathversion{bold}$g$}^{\delta}_{1}(\mbox{\mathversion{bold}$V$})\) yields
\begin{eqnarray*}
\begin{aligned}
\mbox{\mathversion{bold}$g$}^{\delta}_{m}(\mbox{\mathversion{bold}$V$})&=
D\mbox{\mathversion{bold}$f$}_{m-1}(\mbox{\mathversion{bold}$V$})[\mbox{\mathversion{bold}$g$}^{\delta}_{1}(\mbox{\mathversion{bold}$V$})]
+ \delta D\mbox{\mathversion{bold}$r$}^{\delta}_{m-1}(\mbox{\mathversion{bold}$V$})[\mbox{\mathversion{bold}$g$}^{\delta}_{1}(\mbox{\mathversion{bold}$V$})]\\
&= D\mbox{\mathversion{bold}$f$}_{m-1}(\mbox{\mathversion{bold}$V$})[\mbox{\mathversion{bold}$f$}_{1}(\mbox{\mathversion{bold}$V$})]
+ \delta D\mbox{\mathversion{bold}$f$}_{m-1}(\mbox{\mathversion{bold}$V$})[\mbox{\mathversion{bold}$r$}^{\delta}_{1}(\mbox{\mathversion{bold}$V$})]
+\delta D\mbox{\mathversion{bold}$r$}^{\delta}_{m-1}(\mbox{\mathversion{bold}$V$})[\mbox{\mathversion{bold}$g$}^{\delta}_{1}(\mbox{\mathversion{bold}$V$})]\\
&= \mbox{\mathversion{bold}$f$}_{m}(\mbox{\mathversion{bold}$V$}) + \delta \mbox{\mathversion{bold}$r$}^{\delta}_{m}(\mbox{\mathversion{bold}$V$}),
\end{aligned}
\end{eqnarray*}
which shows that (\ref{remain}) holds.

Next we prove that if we choose \( \mbox{\mathversion{bold}$h$}_{\delta}(0)=\mbox{\mathversion{bold}$0$}\) and \( |\mbox{\mathversion{bold}$v$}_{0}|=1\),
\begin{eqnarray}
\left. \left. \mbox{\mathversion{bold}$f$}_{m}(\mbox{\mathversion{bold}$v$}^{\delta}_{0})\right| _{s=0} = \mbox{\mathversion{bold}$f$}_{m}(\mbox{\mathversion{bold}$v$}_{0})
+ \alpha ^{m}\partial ^{3m}_{s}\mbox{\mathversion{bold}$h$}_{\delta}
-\alpha ^{m}(\mbox{\mathversion{bold}$v$}_{0}\cdot \partial ^{3m}_{s}\mbox{\mathversion{bold}$h$}_{\delta})\mbox{\mathversion{bold}$v$}_{0}
+ \mbox{\mathversion{bold}$F$}_{m}(\mbox{\mathversion{bold}$v$}_{0},\mbox{\mathversion{bold}$h$}_{\delta})\right| _{s=0},
\label{vd}
\end{eqnarray}
where \( \mbox{\mathversion{bold}$F$}_{m}(\mbox{\mathversion{bold}$v$}_{0},\mbox{\mathversion{bold}$h$}_{\delta}) \) satisfies
\begin{eqnarray*}
|\mbox{\mathversion{bold}$F$}_{m}(\mbox{\mathversion{bold}$v$}_{0},\mbox{\mathversion{bold}$h$}_{\delta})|\leq 
C( |\mbox{\mathversion{bold}$h$}_{\delta s}| + |\mbox{\mathversion{bold}$h$}_{\delta ss}| + \cdots + |\partial ^{3m-1}_{s}\mbox{\mathversion{bold}$h$}_{\delta}|),
\end{eqnarray*}
if \( |\mbox{\mathversion{bold}$h$}_{\delta s}| + |\mbox{\mathversion{bold}$h$}_{\delta ss}| + \cdots + |\partial ^{3m-1}_{s}\mbox{\mathversion{bold}$h$}_{\delta}|\leq M\), 
where \( C\) depends on \( M\) and \( \mbox{\mathversion{bold}$v$}_{0}\). 
We see from the explicit form (\ref{cor}) of \( \mbox{\mathversion{bold}$v$}^{\delta}_{0}\) that for a natural number \( n\), 
\( \partial^{n}_{s}\mbox{\mathversion{bold}$v$}^{\delta}_{0}\) has the form
\begin{eqnarray}
\left. \partial _{s}^{n}\mbox{\mathversion{bold}$v$}^{\delta}_{0}\right|_{s=0} = \left. 
\partial ^{n}_{s}\mbox{\mathversion{bold}$v$}_{0} + \partial ^{n}_{s}\mbox{\mathversion{bold}$h$}_{\delta}
-(\mbox{\mathversion{bold}$v$}_{0}\cdot \partial ^{n}_{s}\mbox{\mathversion{bold}$h$}_{\delta})\mbox{\mathversion{bold}$v$}_{0}  
+ \mbox{\mathversion{bold}$q$}_{n}(\mbox{\mathversion{bold}$v$}_{0}, \mbox{\mathversion{bold}$h$}_{\delta})\right|_{s=0}.
\label{vdd}
\end{eqnarray}
Here, we have used \( \mbox{\mathversion{bold}$h$}_{\delta}(0)= \mbox{\mathversion{bold}$0$}\) and 
\( \mbox{\mathversion{bold}$q$}_{n}(\mbox{\mathversion{bold}$v$}_{0},\mbox{\mathversion{bold}$h$}_{\delta})\) are terms containing
derivatives of \( \mbox{\mathversion{bold}$v$}_{0}\) and \( \mbox{\mathversion{bold}$h$}_{\delta}\) up to order \( n-1\), and satisfies
\begin{eqnarray*}
|\mbox{\mathversion{bold}$q$}_{n}(\mbox{\mathversion{bold}$v$}_{0},\mbox{\mathversion{bold}$h$}_{\delta})|\leq 
C( |\mbox{\mathversion{bold}$h$}_{\delta s}| + |\mbox{\mathversion{bold}$h$}_{\delta ss}| + \cdots + |\partial ^{n-1}_{s}\mbox{\mathversion{bold}$h$}_{\delta}|),
\end{eqnarray*}
if \( |\mbox{\mathversion{bold}$h$}_{\delta s}| + |\mbox{\mathversion{bold}$h$}_{\delta ss}| + \cdots + |\partial ^{n-1}_{s}\mbox{\mathversion{bold}$h$}_{\delta}|\leq M\), 
where \( C\) depends on \( M\) and \( \mbox{\mathversion{bold}$v$}_{0}\). 
From the definition of \( \mbox{\mathversion{bold}$f$}_{m}(\mbox{\mathversion{bold}$v$}^{\delta}_{0})\), we see that the term with the highest order of 
derivative is \( \alpha ^{m}\partial ^{3m}_{s}\mbox{\mathversion{bold}$v$}^{\delta}_{0}\), so combining this with (\ref{vdd}) yields (\ref{vd}).

Finally, we prove by induction that the differential coefficients of \( \mbox{\mathversion{bold}$h$}_{\delta}\) can be chosen so that
\( \mbox{\mathversion{bold}$v$}^{\delta}_{0}\) satisfies the compatibility conditions for (\ref{dneg2}), and all the coefficients
are \( O(\delta) \).
First, we choose \( \mbox{\mathversion{bold}$h$}_{\delta}(0) =
\partial _{s}\mbox{\mathversion{bold}$h$}_{\delta}(0)=\mbox{\mathversion{bold}$0$}\). This insures that \( \mbox{\mathversion{bold}$v$}^{\delta}_{0}\) 
satisfies the \(0\)-th order compatibility condition. Suppose that the differential coefficients of \( \mbox{\mathversion{bold}$h$}_{\delta}\) up to order 
\( 1+3(m-1)\) are chosen so that they are \( O(\delta) \) and the compatibility conditions for (\ref{dneg2}) up to order \( m-1\) is satisfied, i.e. 
\( \mbox{\mathversion{bold}$g$}^{\delta}_{k}(\mbox{\mathversion{bold}$v$}^{\delta}_{0})_{s}(0)=\mbox{\mathversion{bold}$0$}\) for all \( 0\leq k\leq m-1\).
By choosing \( \mbox{\mathversion{bold}$V$}=\mbox{\mathversion{bold}$v$}^{\delta}_{0}\), we have from (\ref{inner})
\begin{eqnarray*}
\sum ^{m}_{k=0}
\left(
\begin{array}{c}
m\\
k
\end{array}\right)
\mbox{\mathversion{bold}$g$}^{\delta}_{k}(\mbox{\mathversion{bold}$v$}^{\delta}_{0})\cdot \mbox{\mathversion{bold}$g$}^{\delta}_{m-k}(\mbox{\mathversion{bold}$v$}^{\delta}_{0})\equiv 0.
\end{eqnarray*}
Taking the \( s\) derivative of the above and using the assumption of induction yield 
\begin{eqnarray}
\mbox{\mathversion{bold}$v$}^{\delta}_{0}(0)\cdot \mbox{\mathversion{bold}$g$}^{\delta}_{m}(\mbox{\mathversion{bold}$v$}^{\delta}_{0})_{s}(0)
=\mbox{\mathversion{bold}$v$}_{0}(0)\cdot \mbox{\mathversion{bold}$g$}^{\delta}_{m}(\mbox{\mathversion{bold}$v$}^{\delta}_{0})_{s}(0)
= \mbox{\mathversion{bold}$0$}.
\label{vddd}
\end{eqnarray}
Now, from (\ref{remain}) and (\ref{vd}) we have at \( s=0\)
\begin{eqnarray*}
\begin{aligned}
\mbox{\mathversion{bold}$g$}^{\delta}_{m}(\mbox{\mathversion{bold}$v$}^{\delta}_{0})_{s} &= \mbox{\mathversion{bold}$f$}_{m}(\mbox{\mathversion{bold}$v$}^{\delta}_{0})_{s}
+ \delta \mbox{\mathversion{bold}$r$}^{\delta}_{m}(\mbox{\mathversion{bold}$v$}^{\delta}_{0})_{s} \\
&= \mbox{\mathversion{bold}$f$}_{m}(\mbox{\mathversion{bold}$v$}_{0})_{s} + \alpha ^{m}\partial ^{3m+1}_{s}\mbox{\mathversion{bold}$h$}_{\delta}
-\alpha ^{m}(\mbox{\mathversion{bold}$v$}_{0}\cdot \partial ^{3m+1}_{s}\mbox{\mathversion{bold}$h$}_{\delta})\mbox{\mathversion{bold}$v$}_{0} + 
\mbox{\mathversion{bold}$F$}_{m}(\mbox{\mathversion{bold}$v$}_{0},\mbox{\mathversion{bold}$h$}_{\delta})_{s} 
+ \delta\mbox{\mathversion{bold}$r$}^{\delta}_{m}(\mbox{\mathversion{bold}$v$}^{\delta}_{0})_{s}\\
&= \alpha ^{m}\partial ^{3m+1}_{s}\mbox{\mathversion{bold}$h$}_{\delta}
-\alpha ^{m}(\mbox{\mathversion{bold}$v$}_{0}\cdot \partial ^{3m+1}_{s}\mbox{\mathversion{bold}$h$}_{\delta})\mbox{\mathversion{bold}$v$}_{0} + 
\mbox{\mathversion{bold}$F$}_{m}(\mbox{\mathversion{bold}$v$}_{0},\mbox{\mathversion{bold}$h$}_{\delta})_{s} 
+ \delta\mbox{\mathversion{bold}$r$}^{\delta}_{m}(\mbox{\mathversion{bold}$v$}^{\delta}_{0})_{s}.
\end{aligned}
\end{eqnarray*}
From (\ref{vddd}), we see that 
\begin{eqnarray*}
\left. \big( \mbox{\mathversion{bold}$F$}_{m}(\mbox{\mathversion{bold}$v$}_{0},\mbox{\mathversion{bold}$h$}_{\delta})_{s} 
+ \delta\mbox{\mathversion{bold}$r$}^{\delta}_{m}(\mbox{\mathversion{bold}$v$}^{\delta}_{0})_{s}\big) \cdot \mbox{\mathversion{bold}$v$}_{0}\right|_{s=0}
= \mbox{\mathversion{bold}$0$}.
\end{eqnarray*}
From the assumption of induction we have
\begin{eqnarray*}
\left. \mbox{\mathversion{bold}$F$}_{m}(\mbox{\mathversion{bold}$v$}_{0},\mbox{\mathversion{bold}$h$}_{\delta})_{s} 
+ \delta\mbox{\mathversion{bold}$r$}^{\delta}_{m}(\mbox{\mathversion{bold}$v$}^{\delta}_{0})_{s}\right|_{s=0} = O(\delta ).
\end{eqnarray*}
So if we choose \( \partial ^{3m-1}_{s}\mbox{\mathversion{bold}$h$}_{\delta}(0) = \partial ^{3m}_{s}\mbox{\mathversion{bold}$h$}_{\delta}(0)
= \mbox{\mathversion{bold}$0$}\) and \( \partial ^{3m+1}_{s}\mbox{\mathversion{bold}$h$}_{\delta}(0)= 
-\frac{1}{\alpha ^{m}}\big( \left. \mbox{\mathversion{bold}$F$}_{m}(\mbox{\mathversion{bold}$v$}_{0},\mbox{\mathversion{bold}$h$}_{\delta})_{s} 
+ \delta\mbox{\mathversion{bold}$r$}^{\delta}_{m}(\mbox{\mathversion{bold}$v$}^{\delta}_{0})_{s}\big) \right|_{s=0} \), they are all
\( O(\delta) \) and 
\( \mbox{\mathversion{bold}$g$}^{\delta}(\mbox{\mathversion{bold}$v$}^{\delta}_{0})_{s}(0) = \mbox{\mathversion{bold}$0$} \), i.e. the \( m\)-th order 
compatibility condition is satisfied. The differential coefficients are then used to define \( \mbox{\mathversion{bold}$h$}_{\delta}(s)\) as
\begin{eqnarray*}
\mbox{\mathversion{bold}$h$}_{\delta}(s)= \phi (s)\left( \sum ^{m}_{j=0} \frac{\partial ^{3j+1}_{s}\mbox{\mathversion{bold}$h$}_{\delta}(0)}{(3j+1)!}s^{3j+1}\right),
\end{eqnarray*}
where \( \phi (s)\) is a smooth cut-off function that is \( 1\) near \( s=0\). 
We summarize the arguments so far in the following statement.
\begin{lm}
For initial datum \( \mbox{\mathversion{bold}$v$}_{0}\) with \( |\mbox{\mathversion{bold}$v$}_{0}|=1\), \( \mbox{\mathversion{bold}$v$}_{0s}\in H^{1+3m}(\mathbf{R}_{+})\), and
satisfying the compatibility conditions for {\rm (\ref{neg2})} up to order \( m\), we can construct a corrected initial datum \( \mbox{\mathversion{bold}$v$}^{\delta}_{0}\)
such that \( |\mbox{\mathversion{bold}$v$}^{\delta}_{0}|=1\), \( \mbox{\mathversion{bold}$v$}^{\delta}_{0s}\in H^{1+3m}(\mathbf{R}_{+})\), satisfies 
the compatibility conditions of {\rm (\ref{dneg2})} up to order \( m\), and 
\begin{eqnarray*}
\mbox{\mathversion{bold}$v$}^{\delta}_{0}\rightarrow \mbox{\mathversion{bold}$v$}_{0} \ {\rm in} \ L^{\infty}(\mathbf{R}_{+}), \ 
\mbox{\mathversion{bold}$v$}^{\delta}_{0s}\rightarrow \mbox{\mathversion{bold}$v$}_{0s} \ {\rm in} \ H^{1+3m}(\mathbf{R}_{+})
\end{eqnarray*}
as \( \delta \rightarrow +0\).
\label{idc}
\end{lm}
Similar arguments can be used to prove that we can approximate \( \mbox{\mathversion{bold}$v$}_{0}\) by a smoother function while satisfying the 
necessary compatibility conditions by following the method used in Rauch and Massey \cite{1}.

%
%
%
%
%
%
\section{Existence Theorems for Associated Linear Problems}
\setcounter{equation}{0}
We consider the linear problem associated to the regularized nonlinear problem.
If we linearize the nonlinear problem around a function \( \mbox{\mathversion{bold}$w$}\)
and neglect lower order terms, we obtain the following problems. For \( \alpha <0\),
\begin{eqnarray}
\left\{
\begin{array}{ll}
\mbox{\mathversion{bold}$v$}_{t}= \alpha \mbox{\mathversion{bold}$v$}_{sss} + \delta \mbox{\mathversion{bold}$v$}_{ss}
+ \mbox{\mathversion{bold}$w$}\times \mbox{\mathversion{bold}$v$}_{ss} +
3\alpha \mbox{\mathversion{bold}$v$}_{ss}\times (\mbox{\mathversion{bold}$w$}\times \mbox{\mathversion{bold}$w$}_{s}) + \mbox{\mathversion{bold}$f$}, & s>0,t>0,\\
\mbox{\mathversion{bold}$v$}(s,0)=\mbox{\mathversion{bold}$v$}_{0}(s), & s>0,\\
\mbox{\mathversion{bold}$v$}_{s}(0,t)=\mbox{\mathversion{bold}$0$}, & t>0,
\end{array}\right.
\label{lneg}
\end{eqnarray}
and for \( \alpha >0\),
\begin{eqnarray}
\left\{
\begin{array}{ll}
\mbox{\mathversion{bold}$v$}_{t}= \alpha \mbox{\mathversion{bold}$v$}_{sss} + \delta \mbox{\mathversion{bold}$v$}_{ss}
+ \mbox{\mathversion{bold}$w$}\times \mbox{\mathversion{bold}$v$}_{ss} +
3\alpha \mbox{\mathversion{bold}$v$}_{ss}\times (\mbox{\mathversion{bold}$w$}\times \mbox{\mathversion{bold}$w$}_{s}) + \mbox{\mathversion{bold}$f$}, & s>0,t>0,\\
\mbox{\mathversion{bold}$v$}(s,0)= \mbox{\mathversion{bold}$v$}_{0}(s), & s>0,\\
\mbox{\mathversion{bold}$v$}(0,t)=\mbox{\mathversion{bold}$e$}_{3}, & t>0,\\
\mbox{\mathversion{bold}$v$}_{s}(0,t)=\mbox{\mathversion{bold}$0$}, & t>0.
\end{array}\right.
\label{lpos}
\end{eqnarray}
The existence and uniqueness of solution to (\ref{lneg}) and (\ref{lpos}) can be shown as an application 
of existence theorems for a more general equation obtained in  Aiki and Iguchi \cite{10}.
In \cite{10}, we obtained existence theorems for a linear second order parabolic system with a third order 
dispersive term. The problems considered there are as follows. 
For \( \alpha <0\),
\begin{eqnarray}
\left\{
\begin{array}{ll}
\mbox{\mathversion{bold}$u$}_{t}= \alpha \mbox{\mathversion{bold}$u$}_{xxx}+{\rm A}(\mbox{\mathversion{bold}$w$},\partial_{x})\mbox{\mathversion{bold}$u$} 
+ \mbox{\mathversion{bold}$f$}, &
x>0,t>0,\\
\mbox{\mathversion{bold}$u$}(x,0)=\mbox{\mathversion{bold}$u$}_{0}(x), & x>0, \\
\mbox{\mathversion{bold}$u$}_{x}(0,t)=\mbox{\mathversion{bold}$0$}, & t>0.
\end{array} \right.
\label{lpd1}
\end{eqnarray}
For \( \alpha >0\),
\begin{eqnarray}
\left\{
\begin{array}{ll}
\mbox{\mathversion{bold}$u$}_{t}= \alpha \mbox{\mathversion{bold}$u$}_{xxx}+{\rm A}(\mbox{\mathversion{bold}$w$},\partial_{x})\mbox{\mathversion{bold}$u$}
 + \mbox{\mathversion{bold}$f$}, &
x>0,t>0,\\
\mbox{\mathversion{bold}$u$}(x,0)=\mbox{\mathversion{bold}$u$}_{0}(x), & x>0, \\
\mbox{\mathversion{bold}$u$}(0,t)=\mbox{\mathversion{bold}$e$}, & t>0, \\
\mbox{\mathversion{bold}$u$}_{x}(0,t)=\mbox{\mathversion{bold}$0$}, & t>0.
\end{array} \right.
\label{lpd2}
\end{eqnarray}
Here, \( \mbox{\mathversion{bold}$u$}(x,t)=(u^{1}(x,t), u^{2}(x,t), \ldots  ,u^{m}(x,t)) \) is the unknown vector valued function, 
\( \mbox{\mathversion{bold}$u$}_{0}(x)\), \( \mbox{\mathversion{bold}$w$}(x,t)
= ( w^{1}(x,t),w^{2}(x,t), \ldots ,w^{k}(x,t))\), and \( \mbox{\mathversion{bold}$f$}(x,t)=( f^{1}(x,t), f^{2}(x,t), \ldots ,f^{m}(x,t))\) are known vector valued 
functions, \( \mbox{\mathversion{bold}$e$}\) is an arbitrary constant vector, subscripts denote derivatives with the respective variables, 
\( {\rm A}(\mbox{\mathversion{bold}$w$},\partial_{x})\) is 
a second order differential operator of the form \( {\rm A}({\mbox{\mathversion{bold}$w$}},\partial_{x})={\rm A}_{0}({\mbox{\mathversion{bold}$w$}})\partial ^{2}_{x} 
+ {\rm A}_{1}({\mbox{\mathversion{bold}$w$}})\partial _{x} + {\rm A}_{2}({\mbox{\mathversion{bold}$w$}}) \). \( {\rm A}_{0}, \ {\rm A}_{1}, \ {\rm A}_{2}\) are
smooth matrices and \( {\rm A}({\mbox{\mathversion{bold}$w$}},\partial_{x}) \) is strongly elliptic in the sense that for any bounded domain \( E\) in 
\( \mathbf{R}^{k}\), there is a positive constant \( \delta \) such that for any \( \mbox{\mathversion{bold}$w$}\in E\)
\begin{eqnarray*}
{\rm A}_{0}(\mbox{\mathversion{bold}$w$})
+{\rm A}_{0}(\mbox{\mathversion{bold}$w$})^{*}\geq \delta {\rm I},
\end{eqnarray*}
where \( {\rm I}\) is the unit matrix and \( *\) denotes the adjoint of a matrix.
For the above problems we obtained the following.
\begin{Th}{\rm (Aiki and Iguchi \cite{10}) } Let \( \alpha <0\). 
For any \( T>0\) and an arbitrary non-negative integer \( l\), if \( {\mbox{\mathversion{bold}$u$}}_{0}\in H^{2+3l}(\mathbf{R}_{+})\), 
\( {\mbox{\mathversion{bold}$f$}}\in Y^{l}_{T}\), and \( {\mbox{\mathversion{bold}$w$}}\in Z^{l}_{T}\) satisfy the compatibility conditions up to order \( l\), 
a unique solution \( {\mbox{\mathversion{bold}$u$}}\) of {\rm (\ref{lpd1})} exists such that \( {\mbox{\mathversion{bold}$u$}}\in X^{l}_{T}\).
Furthermore, \( \mbox{\mathversion{bold}$u$}\) satisfies
\begin{eqnarray*}
\| \mbox{\mathversion{bold}$u$}\|_{X^{l}_{T}}\leq C\big( \|\mbox{\mathversion{bold}$u$}_{0}\|_{2+3l} + \|\mbox{\mathversion{bold}$f$}\|_{Y^{l}_{T}}\big),
\end{eqnarray*}
where \( C\) depends on \( \alpha \), \( T\), and \( \|\mbox{\mathversion{bold}$w$}\|_{Z^{l}_{T}}\).
\label{thneg}
\end{Th}
\begin{Th}{\rm (Aiki and Iguchi \cite{10}) } Let \( \alpha >0\). 
For any \( T>0\) and an arbitrary non-negative integer \( l\), if \( {\mbox{\mathversion{bold}$u$}}_{0}\in H^{2+3l}(\mathbf{R}_{+})\), 
\( {\mbox{\mathversion{bold}$f$}}\in Y^{l}_{T}\), and \( {\mbox{\mathversion{bold}$w$}}\in Z^{l}_{T}\) satisfy the compatibility conditions up to order \( l\), 
a unique solution \( {\mbox{\mathversion{bold}$u$}}\) of {\rm (\ref{lpd2})} exists such that \( {\mbox{\mathversion{bold}$u$}}\in X^{l}_{T}\).
Furthermore, \( \mbox{\mathversion{bold}$u$}\) satisfies
\begin{eqnarray*}
\| \mbox{\mathversion{bold}$u$}\|_{X^{l}_{T}}\leq C\big( \|\mbox{\mathversion{bold}$u$}_{0}\|_{2+3l} + \|\mbox{\mathversion{bold}$f$}\|_{Y^{l}_{T}}\big),
\end{eqnarray*}
where \( C\) depends on \( \alpha \), \( T\), and \( \|\mbox{\mathversion{bold}$w$}\|_{Z^{l}_{T}}\).
\label{thpos}
\end{Th}
Here, 
\begin{eqnarray*}
X^{l}_{T}:= \bigcap ^{l}_{j=0} \bigg( C^{j}\big( [0,T];H^{2+3(l-j)}(\mathbf{R}_{+})\big) \cap H^{j}\big( 0,T;H^{3+3(l-j)}(\mathbf{R}_{+})\big) \bigg),
\label{X}
\end{eqnarray*}
\begin{eqnarray*}
Y^{l}_{T}:= \bigg\{ f; \  f\in \bigcap ^{l-1}_{j=0} C^{j}\big( [0,T];H^{2+3(l-1-j)}(\mathbf{R}_{+})\big) , \ \partial ^{l}_{t}f \in L^{2}\big( 0,T; H^{1}(\mathbf{R}_{+})\big) \bigg\},
\end{eqnarray*}
\begin{eqnarray*}
Z^{l}_{T}:= \bigg\{ w; \ w\in \bigcap ^{l-1}_{j=0}C^{j}\big( [0,T];H^{2+3(l-1-j)}(\mathbf{R}_{+})\big), \ \partial ^{l}_{t}w\in L^{\infty}\big (0,T;H^{1}(\mathbf{R}_{+})\big) \bigg\}.
\end{eqnarray*}
We apply these theorems with 
\begin{eqnarray}
{\rm A}(\mbox{\mathversion{bold}$w$},\partial_{x})\mbox{\mathversion{bold}$v$}= \delta \mbox{\mathversion{bold}$v$}_{xx} + \mbox{\mathversion{bold}$w$}\times \mbox{\mathversion{bold}$v$}_{xx} 
+ 3\alpha \mbox{\mathversion{bold}$v$}_{xx}\times ( \mbox{\mathversion{bold}$w$}\times \mbox{\mathversion{bold}$w$}_{x}),
\label{oper}
\end{eqnarray}
which obviously satisfies the assumptions on the elliptic operator. Thus we have existence and uniqueness of the solution to (\ref{lneg}) and (\ref{lpos}).
Based on this linear existence theorem, we construct the solution to (\ref{dneg2}) and (\ref{dpos2}). 
%
%
%
%

%
%
%
\section{Construction of the Solution in the Case \( \alpha <0\)}
\setcounter{equation}{0}

\subsection{Existence of Solution}
We construct the solution by the following iteration scheme. For \( n\geq 2\) and \( R\geq 1\), we define \( \mbox{\mathversion{bold}$v$}^{(n),R}\) as the solution of
\begin{eqnarray*}
\left\{
\begin{array}{ll}
\mbox{\mathversion{bold}$v$}^{(n),R}_{t}=\alpha \mbox{\mathversion{bold}$v$}^{(n),R}_{sss} 
+ {\rm A}(\mbox{\mathversion{bold}$v$}^{(n-1),R},\partial_{s})\mbox{\mathversion{bold}$v$}^{(n),R}
- \frac{3}{2}\alpha |\mbox{\mathversion{bold}$v$}^{(n-1),R}_{s}|^{2}\mbox{\mathversion{bold}$v$}^{(n-1),R}_{s} & \ \\
\hspace*{8cm}+ \delta |\mbox{\mathversion{bold}$v$}^{(n-1),R}_{s}|^{2}\mbox{\mathversion{bold}$v$}^{(n-1),R}, & s>0, t>0,\\
\mbox{\mathversion{bold}$v$}^{(n),R}(s,0)=\mbox{\mathversion{bold}$v$}_{0}^{\delta , R}(s), & s>0,\\
\mbox{\mathversion{bold}$v$}^{(n),R}_{s}(0,t)=\mbox{\mathversion{bold}$0$}, & t>0,
\end{array}\right.
\end{eqnarray*}
where \( {\rm A}(\mbox{\mathversion{bold}$v$}^{(n-1),R},\partial_{s}) \) is the operator (\ref{oper}) in the last section and 
\( \mbox{\mathversion{bold}$v$}^{\delta ,R}_{0}(s)= \phi (\frac{s}{R})\mbox{\mathversion{bold}$v$}^{\delta }_{0}(s)\). Here, 
\( \mbox{\mathversion{bold}$v$}^{\delta }_{0}\) is the modified initial datum constructed in section 3 and
\( \phi (s)\) is a smooth cut-off function satisfying \( 0\leq \phi \leq 1\), \( \phi (s)=1\) for \( 0\leq s\leq 1\), and
\( \phi (s)=0\) for \( s>2\). 
Now, we must choose \( \mbox{\mathversion{bold}$v$}^{(1),R}\) appropriately so that the necessary compatibility conditions are satisfied at each
iteration step. This is accomplished by choosing
\begin{eqnarray*}
\mbox{\mathversion{bold}$v$}^{(1),R}(s,t)= \mbox{\mathversion{bold}$v$}^{\delta ,R}_{0}(s) + \sum^{m}_{j=1} \frac{t^{j}}{j!}
\mbox{\mathversion{bold}$P$}_{(j)}(\mbox{\mathversion{bold}$v$}^{\delta ,R}_{0}(s)),
\end{eqnarray*}
where \( m\) is a fixed natural number and \( P_{(j)}\) is defined in section 3. 
Note that multiplying the initial datum by \( \phi \) does not change the fact that 
\( \mbox{\mathversion{bold}$v$}^{\delta ,R}_{0}\) satisfies the compatibility conditions for (\ref{dneg2}).
Recall that we are assuming that \( \mbox{\mathversion{bold}$v$}^{\delta}_{0} \)
is smooth, satisfies the compatibility conditions up to an arbitrary fixed order, and \( \mbox{\mathversion{bold}$v$}^{\delta}_{0}\rightarrow \mbox{\mathversion{bold}$v$}_{0}\)
in the appropriate function space. More specifically, we assume that \( \mbox{\mathversion{bold}$v$}^{\delta}_{0}\) is smooth enough so that 
\( \mbox{\mathversion{bold}$v$}^{(1),R}\in X^{N}_{T}\) for a large \( N>m\) which will be determined later. For each \( R\geq 1\) and natural number \( n\), 
\( \mbox{\mathversion{bold}$v$}^{(n),R}\) is well-defined by Theorem \ref{thneg} and \( \mbox{\mathversion{bold}$v$}^{(n),R}\in X^{m}_{T}\).
We define the function space \( \tilde{X}^{m}_{T}\) as
\begin{eqnarray*}
\tilde{X}^{m}_{T}:= \big\{ \mbox{\mathversion{bold}$v$}; \mbox{\mathversion{bold}$v$}_{s}\in C\big( [0,T]; H^{1+3m}(\mathbf{R}_{+})\big) \big\}
\cap \bigg\{ \bigcap ^{m}_{j=1}C^{j}\big( [0,T]; H^{2+3j}(\mathbf{R}_{+})\big) \bigg\}
\cap C\big( [0,T]; L^{\infty}(\mathbf{R}_{+})\big).
\end{eqnarray*}
The above function space is the space we construct the solution to the nonlinear problem.
Note that from the definition, we have
\begin{eqnarray*}
\| \mbox{\mathversion{bold}$v$}^{(1),R}\|_{\tilde{X}^{m}_{T}}\leq 1 + C\| \mbox{\mathversion{bold}$v$}^{\delta ,R}_{0s}\| _{2+6m}
(1+\| \mbox{\mathversion{bold}$v$}^{\delta ,R}_{0s}\| _{2+6m})^{1+2m} =: M_{0}.
\end{eqnarray*}
Here, \( C\) depends on \( \alpha \) and \( T\), but not on \( \mbox{\mathversion{bold}$v$}^{\delta ,R}_{0}\).
Note that \( \mbox{\mathversion{bold}$v$}^{\delta ,R}_{0}\rightarrow \mbox{\mathversion{bold}$v$}^{\delta }_{0}\)
in \( \tilde{X}^{m}_{T}\)
 as \( R\rightarrow +\infty\), and
there is a positive constant \( C\) independent of \( R\geq 1\) such that
\begin{eqnarray}
\| \mbox{\mathversion{bold}$v$}^{\delta ,R}_{0s}\|_{1+3m} \leq C\|\mbox{\mathversion{bold}$v$}^{\delta }_{0s}\|_{1+3m}.
\label{iniuni}
\end{eqnarray}
This uniform estimate does not hold for \( \|\mbox{\mathversion{bold}$v$}^{\delta ,R}_{0}\| \) because \( \mbox{\mathversion{bold}$v$}^{\delta }_{0}\) does not
belong to \( L^{2}(\mathbf{R}_{+})\).
We show the uniform boundedness of \( \{\mbox{\mathversion{bold}$v$}^{(n),R}\} \) with respect to \( n\) and \( R\) on some time interval \( [0,T_{0}]\) by induction. 
Suppose that for any \( j\) with \( 1\leq j\leq n-1\), \( \|\mbox{\mathversion{bold}$v$}^{(j),R}\|_{\tilde{X}^{m}_{T}}\leq M\).
Then, by a standard energy estimate, we have
\begin{eqnarray*}
\begin{aligned}
\frac{1}{2}\frac{{\rm d}}{{\rm d}t}\| \mbox{\mathversion{bold}$v$}^{(n),R}_{s}\|^{2} &\leq -\frac{| \alpha |}{2}|\mbox{\mathversion{bold}$v$}^{(n),R}_{ss}(0)|^{2}
+ C\|\mbox{\mathversion{bold}$v$}^{(n),R}_{ss}\|^{2} - \delta \|\mbox{\mathversion{bold}$v$}^{(n),R}_{ss}\|^{2} + CM^{3},\\
\frac{1}{2}\frac{{\rm d}}{{\rm d}t}\| \mbox{\mathversion{bold}$v$}^{(n),R}_{ss}\|^{2} &\leq -\frac{| \alpha |}{2}|\mbox{\mathversion{bold}$v$}^{(n),R}_{sss}(0)|^{2}
-\frac{\delta }{4}\| \mbox{\mathversion{bold}$v$}^{(n),R}_{sss}\|^{2} + CM^{2}\|\mbox{\mathversion{bold}$v$}^{(n),R}_{ss}\|^{2} + CM^{3},
\end{aligned}
\end{eqnarray*}
where \( C\) is independent of \( M\) and \( n\). Combining the above estimates yields for any \( 0\leq t\leq T\), 
\begin{eqnarray*}
\| \mbox{\mathversion{bold}$v$}^{(n),R}_{s}(t)\|^{2}_{1} + \int ^{t}_{0}\|\mbox{\mathversion{bold}$v$}^{(n),R}_{s}(\tau )\|^{2}_{2}{\rm d}\tau 
\leq C{\rm e}^{M^{2}T}\big( \|\mbox{\mathversion{bold}$v$}^{\delta ,R}_{0s}\| ^{2}_{1} + M^{3}T\big),
\end{eqnarray*}
where \( C\) is independent of \( T\), \( M\), and \( n\). For a natural number \( k\) with \( 1\leq k\leq m\), we set 
\( \mbox{\mathversion{bold}$v$}^{(n),k}:= \partial ^{k}_{t}\mbox{\mathversion{bold}$v$}^{(n),R}\), and then \( \mbox{\mathversion{bold}$v$}^{(n),k}\) satisfies
\begin{eqnarray*}
\begin{aligned}
\mbox{\mathversion{bold}$v$}^{(n),k}_{t} &= \alpha \mbox{\mathversion{bold}$v$}^{(n),k}_{sss} 
+ \sum ^{k}_{j=0} 
\left(
\begin{array}{c}
k\\
j
\end{array}\right)
(\mbox{\mathversion{bold}$v$}^{(n-1),j}\times \mbox{\mathversion{bold}$v$}^{(n),k-j}_{ss})\\
&\hspace*{1cm}+ 3\alpha \left\{ \sum ^{k}_{j=0}\sum ^{j}_{i=0}
\left(
\begin{array}{c}
k\\
j
\end{array}\right)
\left(
\begin{array}{c}
j\\
i
\end{array}\right)
\mbox{\mathversion{bold}$v$}^{(n),k-j}_{ss}\times ( \mbox{\mathversion{bold}$v$}^{(n-1),i}\times \mbox{\mathversion{bold}$v$}^{(n-1),j-i}_{s})\right\} \\
&\hspace*{1cm}-\frac{3}{2}\alpha \left\{ \sum ^{k}_{j=0}\sum ^{j}_{i=0}
\left(
\begin{array}{c}
k\\
j
\end{array}\right)
\left(
\begin{array}{c}
j\\
i
\end{array}\right)(\mbox{\mathversion{bold}$v$}^{(n-1),i}_{s}\cdot \mbox{\mathversion{bold}$v$}^{(n-1),j-i}_{s})\mbox{\mathversion{bold}$v$}^{(n-1),k-j}_{s}\right\} \\
&\hspace*{1cm} + \delta \mbox{\mathversion{bold}$v$}^{(n),k}_{ss} + \delta \left\{ 
\sum ^{k}_{j=0}\sum ^{j}_{i=0}
\left(
\begin{array}{c}
k\\
j
\end{array}\right)
\left(
\begin{array}{c}
j\\
i
\end{array}\right)
(\mbox{\mathversion{bold}$v$}^{(n-1),i}_{s}\cdot \mbox{\mathversion{bold}$v$}^{(n-1),j-i}_{s})\mbox{\mathversion{bold}$v$}^{(n-1),k-j}\right\} \\
&=: \alpha \mbox{\mathversion{bold}$v$}^{(n),k}_{sss} + \mbox{\mathversion{bold}$v$}^{(n-1)}\times \mbox{\mathversion{bold}$v$}^{(n),k}_{ss}
+ 3\alpha \mbox{\mathversion{bold}$v$}^{(n),k}_{ss}\times (\mbox{\mathversion{bold}$v$}^{(n-1)}\times \mbox{\mathversion{bold}$v$}^{(n-1)}_{s})
+ \delta \mbox{\mathversion{bold}$v$}^{(n),k}_{ss} + \mbox{\mathversion{bold}$F$}^{k}.
\end{aligned}
\end{eqnarray*}
By a similar energy estimate, we have
\begin{eqnarray*}
\begin{aligned}
\frac{1}{2} \frac{{\rm d}}{{\rm d}t}\| \mbox{\mathversion{bold}$v$}^{(n),k}\|_{2}^{2} &\leq CM^{2}(1+M^{2})\big\{
\|\mbox{\mathversion{bold}$v$}^{(n),k}\|^{2}_{2} + (1+M^{2})^{5}	+ \|  \mbox{\mathversion{bold}$F$}^{k}\|^{2}_{1}\big\}\\
&\leq CM^{2}(1+M^{2})\big\{
\|\mbox{\mathversion{bold}$v$}^{(n),k}\|^{2}_{2} + (1+M^{2})^{5}\big\},
\end{aligned}
\end{eqnarray*}
where we have used \( \|\mbox{\mathversion{bold}$v$}^{(j)}\|_{\tilde{X}^{m}_{T}}\leq M\) for \( 1\leq j\leq n-1\) to estimate
\( \mbox{\mathversion{bold}$F$}^{k}\).
Thus we have 
\begin{eqnarray*}
\| \mbox{\mathversion{bold}$v$}^{(n),k}\|^{2}_{2}\leq C{\rm e}^{CM^{2}(1+M^{2})T}\big\{ \| \mbox{\mathversion{bold}$v$}^{(n),k}(\cdot, 0)\|^{2}_{2}
+ (1+M^{2})^{5}T\big\}.
\end{eqnarray*}
By using the equation, we obtain
\begin{eqnarray*}
\| \mbox{\mathversion{bold}$v$}^{(n),k}(\cdot, 0)\|^{2}_{2}\leq C \| \mbox{\mathversion{bold}$v$}^{\delta ,R}_{0s}\|_{1+3k}^{2}
(1+\| \mbox{\mathversion{bold}$v$}^{\delta ,R}_{0s}\|_{1+3k})^{2+4m},
\end{eqnarray*}
and we have 
\begin{eqnarray*}
\| \mbox{\mathversion{bold}$v$}^{(n),k}\|_{2} ^{2}\leq C{\rm e}^{CM^{2}(1+M^{2})T}\big\{ \| \mbox{\mathversion{bold}$v$}^{\delta ,R}_{0s}\|_{1+3k}^{2}
(1+\| \mbox{\mathversion{bold}$v$}^{\delta ,R}_{0s}\|_{1+3k})^{2+4m}
+ (1+M^{2})^{5}T\big\}.
\end{eqnarray*}
Finally, by using the equation and the above estimates, we can convert the regularity in \( t\) into the regularity in \( s\) and obtain for
\( 1\leq j\leq m\),
\begin{eqnarray*}
\| \mbox{\mathversion{bold}$v$}^{(n),j}\|_{1+3(m-j)} ^{2}\leq C{\rm e}^{CM^{2}(1+M^{2})T}\big\{ \| \mbox{\mathversion{bold}$v$}^{\delta ,R}_{0s}\|_{1+3m}^{2}
(1+\| \mbox{\mathversion{bold}$v$}^{\delta ,R}_{0s}\|_{1+3m})^{2+4m}
+ (1+M^{2})^{5}T\big\}.
\end{eqnarray*}
Thus, by choosing \( M:= C_{0}M_{0}\), with a sufficiently large \( C_{0}>0\) independent of \( n\) and \( R\), there is a 
\( T_{0}>0\) such that 
\begin{eqnarray*}
\| \mbox{\mathversion{bold}$v$}^{(n),R}_{s}(t)\|^{2}_{1+3m} + \sum ^{m}_{j=1}\| \partial ^{j}_{t}\mbox{\mathversion{bold}$v$}^{(n),R}(t)\|^{2}_{2+3(m-j)}
\leq \frac{C_{0}M_{0}}{2}
\end{eqnarray*}

Next, we estimate the solution in \( C\big( [0,T];L^{\infty}(\mathbf{R}_{+})\big) \). To do this, we introduce a new variable
\( \mbox{\mathversion{bold}$W$}^{(n),R}:= \mbox{\mathversion{bold}$v$}^{(n),R}-\mbox{\mathversion{bold}$v$}^{\delta ,R}_{0} \). Then, 
\( \mbox{\mathversion{bold}$W$}^{(n),R}\) satisfies
\begin{eqnarray*}
\left\{ 
\begin{array}{ll}
\mbox{\mathversion{bold}$W$}^{(n),R}_{t} = \alpha \mbox{\mathversion{bold}$W$}^{(n),R}_{sss} + \mbox{\mathversion{bold}$v$}^{(n-1),R}\times 
\mbox{\mathversion{bold}$W$}^{(n),R}_{ss} + 3\alpha \mbox{\mathversion{bold}$W$}^{(n),R}_{ss}\times ( \mbox{\mathversion{bold}$v$}^{(n-1),R}\times
\mbox{\mathversion{bold}$v$}^{(n-1),R})
+ \delta \mbox{\mathversion{bold}$W$}^{(n),R}_{ss} & \ \ \\
\hspace*{2cm} -\frac{3}{2}\alpha |\mbox{\mathversion{bold}$v$}^{(n-1),R}_{s}|^{2}\mbox{\mathversion{bold}$v$}^{(n-1),R}_{s}
+ \delta |\mbox{\mathversion{bold}$v$}^{(n-1),R}_{s}|^{2}\mbox{\mathversion{bold}$v$}^{(n-1),R}
+ \alpha \mbox{\mathversion{bold}$v$}^{\delta ,R}_{0sss}  +\mbox{\mathversion{bold}$v$}^{(n-1),R}\times \mbox{\mathversion{bold}$v$}^{\delta ,R}_{0ss} & \ \ \\
\hspace*{2cm}+ 3\alpha \mbox{\mathversion{bold}$v$}^{\delta ,R}_{0ss}\times ( \mbox{\mathversion{bold}$v$}^{(n-1),R}\times \mbox{\mathversion{bold}$v$}^{(n-1),R}_{s})
+ \delta \mbox{\mathversion{bold}$v$}^{\delta ,R}_{0ss}, & \hspace*{-2cm}s>0,t>0,\\
\mbox{\mathversion{bold}$W$}^{(n),R}(s,0) = \mbox{\mathversion{bold}$0$}, & \hspace*{-2cm}s>0,\\
\mbox{\mathversion{bold}$W$}^{(n),R}_{s}(0,t) = \mbox{\mathversion{bold}$0$},& \hspace*{-2cm}t>0.
\end{array}\right.
\end{eqnarray*}
We have by a direct calculation,
\begin{eqnarray*}
\frac{1}{2}\frac{{\rm d}}{{\rm d}t}\| \mbox{\mathversion{bold}$W$}^{(n),R}\|^{2} \leq 
-\frac{\delta }{2} \|\mbox{\mathversion{bold}$W$}^{(n),R}_{s}\|^{2} + C\big( \| \mbox{\mathversion{bold}$W$}^{(n),R}\|^{2} + (1+M)^{3}\big).
\end{eqnarray*}
Thus we have 
\begin{eqnarray*}
\begin{aligned}
\| \mbox{\mathversion{bold}$v$}^{(n),R}\|_{L^{\infty}(\mathbf{R}_{+})} &\leq \| \mbox{\mathversion{bold}$W$}^{(n),R}\|_{L^{\infty}(\mathbf{R}_{+})}
+ \| \mbox{\mathversion{bold}$v$}^{\delta ,R}_{0}\|_{L^{\infty}(\mathbf{R}_{+})}\\
&\leq C\|\mbox{\mathversion{bold}$W$}^{(n),R}\|_{1} + 1\\
&\leq C(1 + M^{2})^{5}T + 1.
\end{aligned}
\end{eqnarray*}
Thus, by choosing \( T_{0}\) smaller if necessary, we have a uniform estimate of the form
\( \| \mbox{\mathversion{bold}$v$}^{(n),R}\|^{2}_{\tilde{X}^{m}_{T_{0}}} \leq C_{0}M^{2}_{0}\). 

Next we show that \( \{ \mbox{\mathversion{bold}$v$}^{(n),R} \} _{n\geq 1}\) converges. Set \( \mbox{\mathversion{bold}$V$}^{(n),R}:=
\mbox{\mathversion{bold}$v$}^{(n),R}-\mbox{\mathversion{bold}$v$}^{(n-1),R}\) for \( n\geq 2\).
Then, \( \mbox{\mathversion{bold}$V$}^{(n),R}\) satisfies
\begin{eqnarray*}
\left\{
\begin{array}{ll}
\mbox{\mathversion{bold}$V$}^{(n),R}_{t} = \alpha \mbox{\mathversion{bold}$V$}^{(n),R}_{sss} + \mbox{\mathversion{bold}$v$}^{(n-1),R}\times \mbox{\mathversion{bold}$V$}^{(n),R}_{ss}
+ 3\alpha \mbox{\mathversion{bold}$V$}^{(n),R}_{ss}\times (\mbox{\mathversion{bold}$v$}^{(n-1),R}\times \mbox{\mathversion{bold}$v$}^{(n-1),R}_{s}) & \\
\hspace*{10.5cm}+ \delta \mbox{\mathversion{bold}$V$}^{(n),R}_{ss} + \mbox{\mathversion{bold}$G$}^{R}_{n}, & s>0,t>0, \\
\mbox{\mathversion{bold}$V$}^{(n),R}(s,0)=\mbox{\mathversion{bold}$0$}, & s>0,\\
\mbox{\mathversion{bold}$V$}^{(n),R}_{s}(0,t)=\mbox{\mathversion{bold}$0$}, & t>0,
\end{array}\right.
\end{eqnarray*}
where \( \mbox{\mathversion{bold}$G$}^{R}_{n}\) are terms depending linearly on \( \mbox{\mathversion{bold}$V$}^{(n-1),R}\).
In the same way we estimated \( \mbox{\mathversion{bold}$v$}^{(n),R}\), we have
\begin{eqnarray*}
\begin{aligned}
\| \mbox{\mathversion{bold}$V$}^{(n),R}(t)\|^{2}_{2} + \int ^{t}_{0}\|\mbox{\mathversion{bold}$V$}^{(n),R}_{s}(\tau )\|^{2}_{2}{\rm d}\tau
&\leq C\int ^{t}_{0} \| \mbox{\mathversion{bold}$V$}^{(n-1),R}(\tau )\|^{2}_{2}{\rm d}\tau \\
&\leq \frac{(CT_{0})^{n-1}}{(n-1)!}M_{0},
\end{aligned}
\end{eqnarray*}
which implies that \( \mbox{\mathversion{bold}$v$}^{(n),R}\) converges to some \( \mbox{\mathversion{bold}$v$}^{R}\) in \( \tilde{X}^{0}_{T_{0}}\).
Combining this convergence, uniform estimate, and the interpolation inequality, we see that \( \mbox{\mathversion{bold}$v$}^{(n),R}\) converges
to \( \mbox{\mathversion{bold}$v$}^{R}\) in \( \tilde{X}^{m-1}_{T_{0}}\). 
Since 
we have approximated the initial datum as smooth as we desire, the above argument implies that for any natural number \( m\), we can construct a solution \( \mbox{\mathversion{bold}$v$}^{R}\) to 
(\ref{dneg2}) in \( \tilde{X}^{m}_{T_{0}}\) with the initial datum \( \mbox{\mathversion{bold}$v$}^{\delta ,R}_{0}\).

Finally, we take the limit \( R\rightarrow +\infty \). First, from the estimate uniform in \( R\), we have 
\( \sum^{2}_{j=0}\sum^{2(2-j)}_{k=0}\sup _{t}\| \partial ^{j}_{t}\partial ^{k}_{s}\mbox{\mathversion{bold}$v$}^{R}(t)\|_{L^{\infty}(\mathbf{R}_{+})}
\leq C\) with \( C>0\) independent of \( R\). 
Therefore, by a standard compactness argument, we see that there is a subsequence \( \{ \mbox{\mathversion{bold}$v$}^{R_{j}}\}_{j}\) and \(\mbox{\mathversion{bold}$v$}\) such that
for \( l=0,1\) and \( 0\leq k\leq 2(1-l)\), 
\( \partial ^{l}_{t}\partial ^{k}_{s}\mbox{\mathversion{bold}$v$}^{R_{j}}\rightarrow 
\partial ^{l}_{t}\partial ^{k}_{s}\mbox{\mathversion{bold}$v$}\) uniformly in any compact subset of \( [0,T_{0}]\times \mathbf{R}_{+}\).

On the other hand, from the uniform estimate, we see that there is a subsequence of \( \{ \mbox{\mathversion{bold}$v$}^{R_{j}}\}_{j}\), which we also denote by
\( \{ \mbox{\mathversion{bold}$v$}^{R_{j}}\}_{j}\), such that
\( \mbox{\mathversion{bold}$v$}^{R_{j}}\) converges to \( \mbox{\mathversion{bold}$v$}\) weakly\( \ast \)  in 
\( \tilde{\tilde{X}}^{m}_{T_{0}}:= \big\{  \mbox{\mathversion{bold}$v$}; \mbox{\mathversion{bold}$v$}_{s}\in L^{\infty}\big( [0,T_{0}]; H^{1+3m}(\mathbf{R}_{+})\big) \big\}
\cap \big\{ \bigcap ^{m}_{j=1}W^{j,\infty}\big( [0,T_{0}]; H^{2+3j}(\mathbf{R}_{+})\big) \big\}\cap L^{\infty}\big( [0,T_{0}]\times \mathbf{R}_{+}\big) \). 
Combining these convergence, we see that we have a solution \( \mbox{\mathversion{bold}$v$}\) of (\ref{dneg2}) with
\( \mbox{\mathversion{bold}$v$}\in \tilde{\tilde{X}}^{m}_{T_{0}}\). 
By taking \( N\), which was mentioned in the beginning of this subsection, large enough, we can construct a solution 
\( \mbox{\mathversion{bold}$v$}\in \tilde{\tilde{X}}^{l}_{T_{0}} \) for any fixed \( l\). By taking \( l>m+1\), this solution belongs to 
\( \tilde{X}^{m}_{T_{0}}\), which follows from Sobolev's embedding with respect to \( t\). 

We summarize the conclusion of this subsection.
\begin{pr}
For a natural number \( m\) and \( \delta >0 \), there exists a \( T_{0}>0\) such that a unique solution \( \mbox{\mathversion{bold}$v$}^{\delta} 
\in \tilde{X}^{m}_{T_{0}}\) to {\rm (\ref{dneg2})}
exists with smooth initial datum \( \mbox{\mathversion{bold}$v$}^{\delta }_{0}\). 
\end{pr}
%
%

%
%
%

\subsection{Uniform Estimate of the Solution with respect to \( \delta \)}

To take the limit \( \delta \rightarrow +0\), we need to obtain uniform estimate of the solution.
To do this, we first show a property of the solution to (\ref{dneg2}) that is very important in the upcoming analysis. 
In the following, we omit the superscript \( \delta \) on the solution to shorten the notation.

\begin{lm}
If \( \mbox{\mathversion{bold}$v$}\) is a solution of {\rm (\ref{dneg2})}
with \( \mbox{\mathversion{bold}$v$}_{s}\in C\big([0,T], H^{2}(\mathbf{R}_{+})\big) \), \( \mbox{\mathversion{bold}$v$}\in C\big( [0,T]; L^{\infty}(\mathbf{R}_{+})\big) \), 
and \( |\mbox{\mathversion{bold}$v$}^{\delta}_{0}|=1\), then \( |\mbox{\mathversion{bold}$v$}|=1\) 
in \( \mathbf{R}_{+}\times [0,T]\).
\label{unit}
\end{lm}
{\it Proof.} Following Nishiyama and Tani \cite{5}, we set \( h(s,t)= |\mbox{\mathversion{bold}$v$}(s,t)|^{2} -1\). From direct calculation and from the fact that 
\( \mbox{\mathversion{bold}$v$}\) is a solution of (\ref{dneg2}), we have
\begin{eqnarray*}
\begin{aligned}
h_{t}&= 2\mbox{\mathversion{bold}$v$}\cdot \mbox{\mathversion{bold}$v$}_{t}\\
&= 2\big\{ \alpha \mbox{\mathversion{bold}$v$}\cdot \mbox{\mathversion{bold}$v$}_{sss} + 3\alpha \mbox{\mathversion{bold}$v$}\cdot
( \mbox{\mathversion{bold}$v$}_{ss}\times (\mbox{\mathversion{bold}$v$}\times \mbox{\mathversion{bold}$v$}_{s}))
-\frac{3}{2} \alpha |\mbox{\mathversion{bold}$v$}_{s}|^{2}(\mbox{\mathversion{bold}$v$}\cdot \mbox{\mathversion{bold}$v$}_{s})
+ \delta (\mbox{\mathversion{bold}$v$}\cdot \mbox{\mathversion{bold}$v$}_{ss}) + \delta |\mbox{\mathversion{bold}$v$}_{s}|^{2}|\mbox{\mathversion{bold}$v$}|^{2}\big\} \\
&=\alpha h_{sss} + \delta h_{ss} + (2\delta |\mbox{\mathversion{bold}$v$}_{s}|^{2} + 3\alpha (\mbox{\mathversion{bold}$v$}_{s}
\cdot \mbox{\mathversion{bold}$v$}_{ss}))h.
\end{aligned}
\end{eqnarray*}
Thus, \( h\) satisfies 
\begin{eqnarray*}
\left\{
\begin{array}{ll}
h_{t}=\alpha h_{sss} + \delta h_{ss} + [2\delta |\mbox{\mathversion{bold}$v$}_{s}|^{2} 
+ 3\alpha (\mbox{\mathversion{bold}$v$}_{s}\cdot \mbox{\mathversion{bold}$v$}_{ss})]h, & s>0, t>0,\\
h(s,0)= 0, & s>0, \\
h_{s}(0,t) = 0, & t>0.
\end{array}\right.
\end{eqnarray*}
By a standard energy method, we have
\begin{eqnarray*}
\begin{aligned}
&\frac{1}{2}\frac{{\rm d}}{{\rm d}t} \| h\|^{2}\leq -\alpha h(0)\cdot h_{ss}(0)
-\delta \|h_{s}\|^{2} + C\| h\|^{2}_{1}, \\
&\frac{1}{2}\frac{{\rm d}}{{\rm d}t} \|h_{s}\|^{2}=-(h_{ss}, h_{t})\\
&\hspace*{1.65cm} \leq -\frac{|\alpha |}{2}|h_{ss}(0)|^{2} -\frac{\delta }{2} \|h_{ss}\|^{2} 
+ C\|h\|^{2}_{1},
\end{aligned}
\end{eqnarray*}
where \( C\) depends on \( \sup_{0\leq t\leq T}\| \mbox{\mathversion{bold}$v$}_{s}(t)\|_{2} \). Combining the two estimates and applying Gronwall's inequality, we obtain
\( h\equiv 0\), which finishes the proof of the lemma. \hfill \( \Box\)
\bigskip

Now that we have established that \( |\mbox{\mathversion{bold}$v$}|=1\), we rewrite the nonlinear terms in (\ref{dneg2})
into its original form.
\begin{eqnarray}
\left\{
\begin{array}{ll}
\mbox{\mathversion{bold}$v$}_{t}=\mbox{\mathversion{bold}$v$}\times \mbox{\mathversion{bold}$v$}_{ss} 
+ \alpha \big\{ \mbox{\mathversion{bold}$v$}_{sss} + \frac{3}{2}\mbox{\mathversion{bold}$v$}_{ss} \times 
\big( \mbox{\mathversion{bold}$v$}\times \mbox{\mathversion{bold}$v$}_{s}\big)
+ \frac{3}{2}\mbox{\mathversion{bold}$v$}_{s}\times (\mbox{\mathversion{bold}$v$}\times \mbox{\mathversion{bold}$v$}_{ss}) \big\}& \ \\
\hspace*{7cm}+ \delta \big( \mbox{\mathversion{bold}$v$}_{ss} + |\mbox{\mathversion{bold}$v$}_{s}|^{2}
\mbox{\mathversion{bold}$v$}\big) , & s>0, t>0,\\
\mbox{\mathversion{bold}$v$}(s,0)=\mbox{\mathversion{bold}$v$}^{\delta}_{0}(s), & s>0,\\
\mbox{\mathversion{bold}$v$}_{s}(0,t)=\mbox{\mathversion{bold}$0$}, & t>0.
\label{dneg222}	
\end{array}\right.
\end{eqnarray}
We will refer to this form of the problem when estimating the solution.

The following two equalities were derived from the property \( |\mbox{\mathversion{bold}$v$}|=1\) in Nishiyama and Tani \cite{5}
which will be used to derive the uniform estimate.
\begin{eqnarray}
\mbox{\mathversion{bold}$v$}\cdot \partial ^{n}_{s}\mbox{\mathversion{bold}$v$} 
= -\frac{1}{2}\sum ^{n-1}_{j=1}
\left(
\begin{array}{c}
n\\
j
\end{array}\right)
\partial ^{j}_{s}\mbox{\mathversion{bold}$v$}\cdot \partial ^{n-j}_{s}\mbox{\mathversion{bold}$v$}.
\label{par}
\end{eqnarray}
\begin{eqnarray}
\mbox{\mathversion{bold}$v$}_{s}\times \partial ^{n}_{s}\mbox{\mathversion{bold}$v$}
= -[\mbox{\mathversion{bold}$v$}\cdot \partial ^{n}_{s}\mbox{\mathversion{bold}$v$}](\mbox{\mathversion{bold}$v$}\times \mbox{\mathversion{bold}$v$}_{s})
+ [(\mbox{\mathversion{bold}$v$}\times \mbox{\mathversion{bold}$v$}_{s})\cdot \partial ^{n}_{s}\mbox{\mathversion{bold}$v$}]\mbox{\mathversion{bold}$v$}
\ \ \ \ \ \ {\rm for} \ n\geq 2.
\label{par2}
\end{eqnarray}
(\ref{par}) is derived by differentiating the equality \( |\mbox{\mathversion{bold}$v$}|^{2}=1\) with respect to \( s\). 
We show (\ref{par2}) in a little more detail for the convenience of the reader. Suppose \( \mbox{\mathversion{bold}$v$}_{s}\neq \mbox{\mathversion{bold}$0$}\).
Then, since \( |\mbox{\mathversion{bold}$v$}|=1\) and \( \mbox{\mathversion{bold}$v$}\cdot \mbox{\mathversion{bold}$v$}_{s} =0\),
\( \{ \mbox{\mathversion{bold}$v$},  \frac{\mbox{\mathversion{bold}$v$}_{s}}{|\mbox{\mathversion{bold}$v$}_{s}|}, 
\frac{\mbox{\mathversion{bold}$v$}\times \mbox{\mathversion{bold}$v$}_{s}}{|\mbox{\mathversion{bold}$v$}_{s}|}\} \) form a 
orthonormal frame of \( \mathbf{R}^{3}\). Thus for \( n\geq 2\), we have
\begin{eqnarray*}
 \partial ^{n}_{s}\mbox{\mathversion{bold}$v$}
=[ \mbox{\mathversion{bold}$v$}\cdot \partial ^{n}_{s}\mbox{\mathversion{bold}$v$}]
\mbox{\mathversion{bold}$v$} 
+[\frac{\mbox{\mathversion{bold}$v$}_{s}}{|\mbox{\mathversion{bold}$v$}_{s}|}\cdot \partial ^{n}_{s}\mbox{\mathversion{bold}$v$}]\frac{\mbox{\mathversion{bold}$v$}_{s}}{|\mbox{\mathversion{bold}$v$}_{s}|}
+ [ \frac{(\mbox{\mathversion{bold}$v$}\times \mbox{\mathversion{bold}$v$}_{s})}{|\mbox{\mathversion{bold}$v$}_{s}|}\cdot \partial ^{n}_{s}\mbox{\mathversion{bold}$v$}]
\frac{\mbox{\mathversion{bold}$v$}\times \mbox{\mathversion{bold}$v$}_{s}}{|\mbox{\mathversion{bold}$v$}_{s}|}.
\end{eqnarray*}
Taking the exterior product with \( \mbox{\mathversion{bold}$v$}_{s}\) from the left yields
\begin{eqnarray*}
\begin{aligned}
\mbox{\mathversion{bold}$v$}_{s}\times \partial ^{n}_{s}\mbox{\mathversion{bold}$v$}&=
-[\mbox{\mathversion{bold}$v$}\cdot \partial ^{n}_{s}\mbox{\mathversion{bold}$v$}](\mbox{\mathversion{bold}$v$}\times \mbox{\mathversion{bold}$v$}_{s})	
+ [ \frac{(\mbox{\mathversion{bold}$v$}\times \mbox{\mathversion{bold}$v$}_{s})}{|\mbox{\mathversion{bold}$v$}_{s}|}\cdot \partial ^{n}_{s}\mbox{\mathversion{bold}$v$}]
\frac{[\mbox{\mathversion{bold}$v$}_{s}\times (\mbox{\mathversion{bold}$v$}\times \mbox{\mathversion{bold}$v$}_{s})]}{|\mbox{\mathversion{bold}$v$}_{s}|}\\
&=-[\mbox{\mathversion{bold}$v$}\cdot \partial ^{n}_{s}\mbox{\mathversion{bold}$v$}](\mbox{\mathversion{bold}$v$}\times \mbox{\mathversion{bold}$v$}_{s})
+ [(\mbox{\mathversion{bold}$v$}\times \mbox{\mathversion{bold}$v$}_{s})\cdot \partial ^{n}_{s}\mbox{\mathversion{bold}$v$}]\mbox{\mathversion{bold}$v$}.
\end{aligned}
\end{eqnarray*}
When \( \mbox{\mathversion{bold}$v$}_{s}=\mbox{\mathversion{bold}$0$}\), each term in (\ref{par2}) is zero, so that (\ref{par2}) holds in either case.
\bigskip

Now we estimate the solution.
We first derive the basic estimate.
\begin{pr}
Let \( M,T>0\). Suppose that \( \mbox{\mathversion{bold}$v$}\) is a solution of {\rm (\ref{dneg222})} with \( \mbox{\mathversion{bold}$v$}^{\delta}_{0s}\in H^{4}(\mathbf{R}_{+})\),
\( |\mbox{\mathversion{bold}$v$}^{\delta }_{0}|=1\), and \( \|\mbox{\mathversion{bold}$v$}^{\delta }_{0s}\|_{1}\leq M\)
satisfying \( \mbox{\mathversion{bold}$v$}_{s}\in 
C\big( [0,T]; H^{4}(\mathbf{R}_{+})\big) \) and \( \mbox{\mathversion{bold}$v$}\in C\big( [0,T];L^{\infty}(\mathbf{R}_{+})\big) \). 
Then, there exist \( \delta_{*}>0\) and \( C_{*}>0\) such that
for \( \delta \in (0,\delta _{*}]\), the following estimate holds.
\begin{eqnarray*}
\sup _{0\leq t\leq T}\| \mbox{\mathversion{bold}$v$}_{s}(t)\|_{1}\leq C_{*}.
\end{eqnarray*}
We emphasize that \( C_{*}\) depends on \( M \) and \( T\), but not on \( \delta \in (0,\delta _{*}] \).
\label{basest}
\end{pr}
{\it Proof.} From Lemma \ref{unit}, we have \( |\mbox{\mathversion{bold}$v$}|=1\).
We make use of quantities which are conserved for the initial value problem in \( \mathbf{R}\) with \( \delta =0\). First we estimate
\begin{eqnarray*}
\begin{aligned}
\frac{{\rm d}}{{\rm d}t}\| \mbox{\mathversion{bold}$v$}_{s}\|^{2} &= -(\mbox{\mathversion{bold}$v$}_{t},\mbox{\mathversion{bold}$v$}_{ss})\\
&= -\alpha \big\{ (\mbox{\mathversion{bold}$v$}_{sss},\mbox{\mathversion{bold}$v$}_{ss}) 
+ \frac{3}{2}(\mbox{\mathversion{bold}$v$}_{s}\times (\mbox{\mathversion{bold}$v$}\times \mbox{\mathversion{bold}$v$}_{ss}),\mbox{\mathversion{bold}$v$}_{ss})\big\}
- \delta \big\{ \|\mbox{\mathversion{bold}$v$}_{ss}\|^{2} + (|\mbox{\mathversion{bold}$v$}_{s}|^{2}\mbox{\mathversion{bold}$v$},\mbox{\mathversion{bold}$v$}_{ss})\big\} \\
&= -\frac{|\alpha |}{2}|\mbox{\mathversion{bold}$v$}_{ss}(0)|^{2}-\delta \|\mbox{\mathversion{bold}$v$}_{ss}\|^{2}+ \delta \|\mbox{\mathversion{bold}$v$}_{s}\|^{4}_{L^{4}(\mathbf{R}_{+})}\\
&\leq -\frac{|\alpha |}{2}|\mbox{\mathversion{bold}$v$}_{ss}(0)|^{2} - \frac{\delta }{2}\| \mbox{\mathversion{bold}$v$}_{ss}\|^{2} + C\delta \|\mbox{\mathversion{bold}$v$}_{s}\|^{6}.
\end{aligned}
\end{eqnarray*}
Here, \( C\) is independent of \( \delta \) and is determined from the interpolation inequality \( \|\mbox{\mathversion{bold}$v$}_{s}\|_{L^{4}(\mathbf{R}_{+})}
\leq C\|\mbox{\mathversion{bold}$v$}_{s}\|^{3/4}\|\mbox{\mathversion{bold}$v$}_{ss}\|^{1/4}\). Thus, we have 
\( \frac{{\rm d}}{{\rm d}t}\|\mbox{\mathversion{bold}$v$}_{s}\|^{2} \leq C\delta \|\mbox{\mathversion{bold}$v$}_{s}\|^{6}\). On the other hand, the ordinary differential
equation
\begin{eqnarray*}
\left\{
\begin{array}{ll}
r_{t} = C\delta r^{3}, & t>0,\\
r(0)= \| \mbox{\mathversion{bold}$v$}^{\delta }_{0s}\|^{2}
\end{array}\right.
\end{eqnarray*}
has the explicit solution \( r(t) = \big( \|\mbox{\mathversion{bold}$v$}^{\delta}_{0s}\|^{-4}-C\delta t\big) ^{-1/2} \) as long as 
\( \|\mbox{\mathversion{bold}$v$}^{\delta}_{0s}\|^{-4}>C\delta t\). So, if we choose \( \delta _{*}>0\) such that 
\( M^{-4}>C\delta_{*} T\) holds, \( r(t)\) is well-defined on \( [0,T]\) and from the comparison principle,
\begin{eqnarray*}
\begin{aligned}
\| \mbox{\mathversion{bold}$v$}_{s}(t)\| \leq r(t)^{1/2} &=  \big( \|\mbox{\mathversion{bold}$v$}^{\delta}_{0s}\|^{-4}-C\delta t\big) ^{-1/4}\\
&\leq \big( M^{-4}-C\delta_{*} T\big) ^{-1/4}=: C_{1},
\end{aligned}
\end{eqnarray*}
which is a uniform estimate for \( \|\mbox{\mathversion{bold}$v$}_{s}\|\). 
Next we derive a uniform estimate for \( \| \mbox{\mathversion{bold}$v$}_{ss}\|\).
For the initial value problem with \( \delta =0\), this was achieved by fully utilizing the conserved quantity 
\( \|  \mbox{\mathversion{bold}$v$}_{ss}\|^{2} - \frac{5}{4}\| | \mbox{\mathversion{bold}$v$}_{s}|\|^{2}\).
In our case, we also use this quantity but we have to take care of boundary terms. 
\begin{eqnarray*}
\begin{aligned}
\frac{{\rm d}}{{\rm d}t}\left\{ \|\mbox{\mathversion{bold}$v$}_{ss}\|^{2}-\frac{5}{4}\| |\mbox{\mathversion{bold}$v$}_{s}|^{2}\|^{2}\right\}
= 2(\mbox{\mathversion{bold}$v$}_{ss},\mbox{\mathversion{bold}$v$}_{sst})-5(|\mbox{\mathversion{bold}$v$}_{s}|^{2}\mbox{\mathversion{bold}$v$}_{s},\mbox{\mathversion{bold}$v$}_{st})
&= -2(\mbox{\mathversion{bold}$v$}_{sss},\mbox{\mathversion{bold}$v$}_{st}) - 5(|\mbox{\mathversion{bold}$v$}_{s}|^{2}\mbox{\mathversion{bold}$v$}_{s},\mbox{\mathversion{bold}$v$}_{st})\\
&=: I_{1} + \alpha I_{2} + \delta I_{3}.
\end{aligned}
\end{eqnarray*}
We estimate each term separately.
\begin{eqnarray*}
\begin{aligned}
I_{1}&= -2 (\mbox{\mathversion{bold}$v$}_{sss},\mbox{\mathversion{bold}$v$}_{s}\times \mbox{\mathversion{bold}$v$}_{ss}) 
-5(|\mbox{\mathversion{bold}$v$}_{s}|^{2}\mbox{\mathversion{bold}$v$}_{s},\mbox{\mathversion{bold}$v$}\times \mbox{\mathversion{bold}$v$}_{sss}) \\
&=-2 \int_{\mathbf{R}_{+}} (\mbox{\mathversion{bold}$v$}_{s}\cdot \mbox{\mathversion{bold}$v$}_{ss})
\big[ \mbox{\mathversion{bold}$v$}_{s}\cdot (\mbox{\mathversion{bold}$v$}\times \mbox{\mathversion{bold}$v$}_{ss})\big] {\rm d}s
+ 4\int_{\mathbf{R}_{+}}(\mbox{\mathversion{bold}$v$}\cdot \mbox{\mathversion{bold}$v$}_{ss})
\big[ \mbox{\mathversion{bold}$v$}_{sss}\cdot (\mbox{\mathversion{bold}$v$}\times \mbox{\mathversion{bold}$v$}_{s})\big] {\rm d}s\\
&\hspace*{3cm}-5\int_{\mathbf{R}_{+}}|\mbox{\mathversion{bold}$v$}_{s}|^{2}\mbox{\mathversion{bold}$v$}_{s}\cdot (\mbox{\mathversion{bold}$v$}\times \mbox{\mathversion{bold}$v$}_{sss})
 {\rm d}s\\
&= -2\int_{\mathbf{R}_{+}}(\mbox{\mathversion{bold}$v$}_{s}\cdot \mbox{\mathversion{bold}$v$}_{ss})\big[ \mbox{\mathversion{bold}$v$}_{s}\cdot 
(\mbox{\mathversion{bold}$v$}\times \mbox{\mathversion{bold}$v$}_{ss})\big] {\rm d}s
-\int_{\mathbf{R}_{+}}|\mbox{\mathversion{bold}$v$}_{s}|^{2}\big[ \mbox{\mathversion{bold}$v$}_{s}\times (\mbox{\mathversion{bold}$v$}\times \mbox{\mathversion{bold}$v$}_{sss})\big]
 {\rm d}s \\
&= -\int_{\mathbf{R}_{+}}\big\{ |\mbox{\mathversion{bold}$v$}_{s}|^{2}\mbox{\mathversion{bold}$v$}_{s}\cdot 
(\mbox{\mathversion{bold}$v$}\times \mbox{\mathversion{bold}$v$}_{ss})\big\}_{s} {\rm d}s = 0.
\end{aligned}
\end{eqnarray*}
Here, we have used integration by parts, (\ref{par}), and (\ref{par2}). From here on, integration with respect to \( s\) is assumed to be taken over \( \mathbf{R}_{+}\).
Next we have
\begin{eqnarray*}
\begin{aligned}
I_{2} &= -2\int \mbox{\mathversion{bold}$v$}_{sss}\cdot \mbox{\mathversion{bold}$v$}_{ssss} {\rm d}s
-6\int \mbox{\mathversion{bold}$v$}_{sss}\cdot \big[ \mbox{\mathversion{bold}$v$}_{ss}\times ( \mbox{\mathversion{bold}$v$}\times \mbox{\mathversion{bold}$v$}_{ss})\big] {\rm d}s
-3\int \mbox{\mathversion{bold}$v$}_{sss}\cdot \big[ \mbox{\mathversion{bold}$v$}_{s}\times (\mbox{\mathversion{bold}$v$}_{s}\times \mbox{\mathversion{bold}$v$}_{ss})\big] {\rm d}s \\
& \ \ \ -3\int \mbox{\mathversion{bold}$v$}_{sss}\cdot 
\big[ \mbox{\mathversion{bold}$v$}_{s}\times (\mbox{\mathversion{bold}$v$}\times \mbox{\mathversion{bold}$v$}_{sss})\big] {\rm d}s 
-5\int |\mbox{\mathversion{bold}$v$}_{s}|^{2}\mbox{\mathversion{bold}$v$}_{s}\cdot \mbox{\mathversion{bold}$v$}_{ssss} {\rm d}s \\
& \ \ \ -\frac{15}{2}\int |\mbox{\mathversion{bold}$v$}_{s}|^{2}\mbox{\mathversion{bold}$v$}_{s}\cdot 
\big[ \mbox{\mathversion{bold}$v$}_{sss}\times (\mbox{\mathversion{bold}$v$}\times \mbox{\mathversion{bold}$v$}_{s})\big] {\rm d}s 
-15\int |\mbox{\mathversion{bold}$v$}_{s}|^{2}\mbox{\mathversion{bold}$v$}_{s}\cdot 
\big[ \mbox{\mathversion{bold}$v$}_{ss}\times (\mbox{\mathversion{bold}$v$}\times \mbox{\mathversion{bold}$v$}_{ss})\big] {\rm d}s \\
&= |\mbox{\mathversion{bold}$v$}_{sss}(0)|^{2} + 9\int (|\mbox{\mathversion{bold}$v$}_{s}|^{2})_{s}|\mbox{\mathversion{bold}$v$}_{ss}|^{2} {\rm d}s
-3\int |\mbox{\mathversion{bold}$v$}_{s}|^{2}(|\mbox{\mathversion{bold}$v$}_{ss}|^{2})_{s} {\rm d}s
-\frac{3}{2}\int (\mbox{\mathversion{bold}$v$}_{s}\cdot \mbox{\mathversion{bold}$v$}_{sss})(|\mbox{\mathversion{bold}$v$}_{s}|^{2})_{s} {\rm d}s \\
& \ \ \ +\frac{3}{2}\int |\mbox{\mathversion{bold}$v$}_{s}|^{2}(|\mbox{\mathversion{bold}$v$}_{ss}|^{2})_{s} {\rm d}s
-\frac{9}{2}\int (|\mbox{\mathversion{bold}$v$}_{s}|^{2})_{s}|\mbox{\mathversion{bold}$v$}_{ss}|^{2} {\rm d}s
-5\int |\mbox{\mathversion{bold}$v$}_{s}|^{2}\mbox{\mathversion{bold}$v$}_{s}\cdot \mbox{\mathversion{bold}$v$}_{ssss} {\rm d}s \\
& \ \ \ -\frac{45}{4}\int |\mbox{\mathversion{bold}$v$}_{s}|^{4}(|\mbox{\mathversion{bold}$v$}_{s}|^{2})_{s} {\rm d}s 
 -\frac{15}{2}\int |\mbox{\mathversion{bold}$v$}_{s}|^{4}(|\mbox{\mathversion{bold}$v$}_{s}|^{2})_{s} {\rm d}s \\
&=  |\mbox{\mathversion{bold}$v$}_{sss}(0)|^{2} + \frac{9}{2}\int (|\mbox{\mathversion{bold}$v$}_{s}|^{2})_{s}|\mbox{\mathversion{bold}$v$}_{ss}|^{2} {\rm d}s 
-\frac{3}{2}\int |\mbox{\mathversion{bold}$v$}_{s}|^{2}(|\mbox{\mathversion{bold}$v$}_{ss}|^{2})_{s} {\rm d}s
-\frac{25}{4}\int \big\{ |\mbox{\mathversion{bold}$v$}_{s}|^{6}\big\}_{s} {\rm d}s \\
& \ \ \ -\frac{3}{2}\int (\mbox{\mathversion{bold}$v$}_{s}\cdot \mbox{\mathversion{bold}$v$}_{sss})(|\mbox{\mathversion{bold}$v$}_{s}|^{2})_{s} {\rm d}s
-5\int |\mbox{\mathversion{bold}$v$}_{s}|^{2}\mbox{\mathversion{bold}$v$}_{s}\cdot \mbox{\mathversion{bold}$v$}_{ssss} {\rm d}s \\
&=  |\mbox{\mathversion{bold}$v$}_{sss}(0)|^{2} + \int \big\{ |\mbox{\mathversion{bold}$v$}_{s}|^{2}|\mbox{\mathversion{bold}$v$}_{ss}|^{2}\big\}_{s} {\rm d}s
+ \frac{7}{2}\int (|\mbox{\mathversion{bold}$v$}_{s}|^{2})_{s}|\mbox{\mathversion{bold}$v$}_{ss}|^{2} {\rm d}s
+ \frac{7}{2}\int (\mbox{\mathversion{bold}$v$}_{s}\cdot \mbox{\mathversion{bold}$v$}_{sss})(|\mbox{\mathversion{bold}$v$}_{s}|^{2})_{s} {\rm d}s \\
&=  |\mbox{\mathversion{bold}$v$}_{sss}(0)|^{2}.
\end{aligned}
\end{eqnarray*}
Again, we have used integration by parts, (\ref{par}), and (\ref{par2}). Finally, we calculate
\begin{eqnarray*}
\begin{aligned}
I_{3}&= - \| \mbox{\mathversion{bold}$v$}_{sss}\|^{2}- 2( \mbox{\mathversion{bold}$v$}_{sss}, |\mbox{\mathversion{bold}$v$}_{s}|^{2}\mbox{\mathversion{bold}$v$}_{s})
-4( \mbox{\mathversion{bold}$v$}_{sss}, (\mbox{\mathversion{bold}$v$}_{s}\cdot \mbox{\mathversion{bold}$v$}_{ss})\mbox{\mathversion{bold}$v$})
-5( |\mbox{\mathversion{bold}$v$}_{s}|^{2}\mbox{\mathversion{bold}$v$}_{s}, \mbox{\mathversion{bold}$v$}_{sss}) \\
& \ \ \ -5 (|\mbox{\mathversion{bold}$v$}_{s}|^{2}\mbox{\mathversion{bold}$v$}_{s}, |\mbox{\mathversion{bold}$v$}_{s}|^{2}\mbox{\mathversion{bold}$v$}_{s})
-10(|\mbox{\mathversion{bold}$v$}_{s}|^{2}\mbox{\mathversion{bold}$v$}_{s}, (\mbox{\mathversion{bold}$v$}_{s}\cdot \mbox{\mathversion{bold}$v$}_{ss})\mbox{\mathversion{bold}$v$})\\
&\leq - \frac{1}{2}\|\mbox{\mathversion{bold}$v$}_{sss}\|^{2} + C\big( \|\mbox{\mathversion{bold}$v$}_{s}\|_{L^{6}(\mathbf{R}_{+})}^{6} 
+ \| \mbox{\mathversion{bold}$v$}_{s}\cdot \mbox{\mathversion{bold}$v$}_{ss}\|^{2}\big) \\
&\leq -\frac{1}{4}\| \mbox{\mathversion{bold}$v$}_{sss}\|^{2} + C_{2}.
\end{aligned}
\end{eqnarray*}
Here, \( C_{2}\) is a constant depending on \( C_{1}\). We also used the interpolation inequalities
\( \| \mbox{\mathversion{bold}$v$}_{s}\|_{L^{6}(\mathbf{R}_{+})} \leq C\|\mbox{\mathversion{bold}$v$}_{s}\|^{2/3}\|\mbox{\mathversion{bold}$v$}_{ss}\|^{1/3}\),
\( \| \mbox{\mathversion{bold}$v$}_{s}\|_{L^{\infty}(\mathbf{R}_{+})}\leq C \|\mbox{\mathversion{bold}$v$}_{s}\|^{1/2}\|\mbox{\mathversion{bold}$v$}_{ss}\|^{1/2}\), 
and \( \|\mbox{\mathversion{bold}$v$}_{ss}\|\leq C\|\mbox{\mathversion{bold}$v$}_{s}\|^{1/2}\|\mbox{\mathversion{bold}$v$}_{sss}\|^{1/2}\).
By combining the three estimates, we arrive at
\begin{eqnarray*}
\frac{{\rm d}}{{\rm d}t}\left\{ \|\mbox{\mathversion{bold}$v$}_{ss}\|^{2}-\frac{5}{4}\| |\mbox{\mathversion{bold}$v$}_{s}|^{2}\|^{2}\right\}
\leq -|\alpha | |\mbox{\mathversion{bold}$v$}_{sss}(0)|^{2}-\frac{\delta }{4}\|\mbox{\mathversion{bold}$v$}_{sss}\|^{2} + C_{2}.
\end{eqnarray*}
Integrating over \( [0,t]\) yields
\begin{eqnarray*}
\begin{aligned}
\|\mbox{\mathversion{bold}$v$}_{ss}(t)\|^{2} &+ \int^{t}_{0} \big( |\alpha ||\mbox{\mathversion{bold}$v$}_{sss}(0,\tau )|^{2} + \frac{\delta }{4}\|\mbox{\mathversion{bold}$v$}_{sss}(\tau )\|^{2}\big)
{\rm d}\tau \\
&\leq \| \mbox{\mathversion{bold}$v$}^{\delta }_{0ss}\|^{2} + \frac{5}{4}\| |\mbox{\mathversion{bold}$v$}_{s}(t)|^{2}\|^{2}
+ C_{2}t\\
&\leq C\| \mbox{\mathversion{bold}$v$}^{\delta }_{0s}\|^{2}_{1} + \frac{1}{2}\|\mbox{\mathversion{bold}$v$}_{ss}(t)\|^{2} + C\|\mbox{\mathversion{bold}$v$}_{s}(t)\|^{6} + C_{2}t,
\end{aligned}
\end{eqnarray*}
where we have used \( \| \mbox{\mathversion{bold}$v$}_{s}\|_{L^{4}(\mathbf{R}_{+})}\leq C\|\mbox{\mathversion{bold}$v$}_{s}\|^{3/4}\|\mbox{\mathversion{bold}$v$}_{ss}\|^{1/4}\) again.
Thus we have 
\begin{eqnarray*}
\sup _{0\leq t\leq T}\| \mbox{\mathversion{bold}$v$}_{ss}(t)\|^{2} + \int^{T}_{0} \big( |\alpha ||\mbox{\mathversion{bold}$v$}_{sss}(0,t)|^{2} + \delta \|\mbox{\mathversion{bold}$v$}_{sss}(t)\|^{2}\big)
{\rm d}t 
\leq C\| \mbox{\mathversion{bold}$v$}^{\delta }_{0s}\|^{2}_{1} + C_{3} + C_{2}T,
\end{eqnarray*}
where \( C_{3}\) is a constant depending on \( C_{1}\). Thus if we choose 
\( C_{*}^{2}:= CM^{2} + C^{2}_{1} + C_{3} + C_{2}T\), we see that the proposition holds. \hfill \( \Box\)
\bigskip

Based on the estimate derived in Proposition \ref{basest}, we derive the higher order estimate. 
\begin{pr}
For a natural number \( k\) and \( M>0\), let \( \mbox{\mathversion{bold}$v$}\) be a solution of { \rm (\ref{dneg222})} with \( | \mbox{\mathversion{bold}$v$}^{\delta}_{0} |=1\),
\( \mbox{\mathversion{bold}$v$}^{\delta}_{0s}\in H^{1+3k}(\mathbf{R}_{+})\), and \( \| \mbox{\mathversion{bold}$v$}^{\delta }_{0s}\|_{H^{1+3k}(\mathbf{R}_{+})}\leq M\)
satisfying \( \mbox{\mathversion{bold}$v$}_{s}\in 
C\big( [0,T];H^{1+3k}(\mathbf{R}_{+})\big) \) and \( \mbox{\mathversion{bold}$v$}\in C\big( [0,T]; L^{\infty}(\mathbf{R}_{+})\big) \).
Then, there is a \( C_{**}>0\) and \( T_{1}\in (0,T]\) such that for \( 0< \delta \leq \delta_{*}\), \( \mbox{\mathversion{bold}$v$}\) satisfies
\begin{eqnarray*}
\sup_{0\leq t\leq T_{1}}\| \mbox{\mathversion{bold}$v$}_{s}(t)\|_{1+3k}\leq C_{**}.
\end{eqnarray*}
Here, \( T_{1}\) depends on \( \|\mbox{\mathversion{bold}$v$}_{0s}\|_{3}\) and 
\( C_{**}\) depends on \( C_{*}\) and \( \delta _{*}\), but not on \( \delta \in (0,\delta _{*}]\). \( C_{*}\) and \( \delta _{*}\) are
defined in Proposition {\rm \ref{basest}}.
\end{pr}
{\it Proof}. From Proposition \ref{basest}, we have a \( C_{*}>0 \) and \( \delta _{*}>0\) such that 
\begin{eqnarray*}
\sup_{0\leq t\leq T}\| \mbox{\mathversion{bold}$v$}_{s}(t)\|_{1}\leq C_{*}
\end{eqnarray*}
holds for \( \delta \in (0,\delta _{*}]\). We also know from Lemma \ref{unit} that \( |\mbox{\mathversion{bold}$v$}|=1\).

Now, we take the derivative with respect to \( s\) of the equation \( m\) times with \( 4\leq m\leq 1+3k \) to obtain
\begin{eqnarray*}
\begin{aligned}
\partial ^{m}_{s}\mbox{\mathversion{bold}$v$}_{t} &= \mbox{\mathversion{bold}$v$}\times \partial ^{m+2}_{s}\mbox{\mathversion{bold}$v$}
+ m\mbox{\mathversion{bold}$v$}_{s}\times \partial ^{m+1}_{s}\mbox{\mathversion{bold}$v$}
+ \alpha \bigg\{ \partial ^{m+3}_{s}\mbox{\mathversion{bold}$v$} + \frac{3}{2}\big( \partial ^{m+2}_{s}\mbox{\mathversion{bold}$v$}\big) \times 
(\mbox{\mathversion{bold}$v$}\times \mbox{\mathversion{bold}$v$}_{s}) \\
&\hspace*{0.5cm}+ \frac{3}{2}(m+1)\big( \partial^{m+1}_{s}\mbox{\mathversion{bold}$v$}\big) \times
(\mbox{\mathversion{bold}$v$}\times \mbox{\mathversion{bold}$v$}_{ss}) 
+ \frac{3}{2}(m+1)\mbox{\mathversion{bold}$v$}_{ss}\times
(\mbox{\mathversion{bold}$v$}\times \partial ^{m+1}_{s}\mbox{\mathversion{bold}$v$}) \\
&\hspace*{0.5cm} + \frac{3}{2}\mbox{\mathversion{bold}$v$}_{s}\times
(\mbox{\mathversion{bold}$v$}\times \partial ^{m+2}_{s}\mbox{\mathversion{bold}$v$}) + \frac{3m}{2}\mbox{\mathversion{bold}$v$}_{s}\times
(\mbox{\mathversion{bold}$v$}_{s}\times \partial ^{m+1}_{s}\mbox{\mathversion{bold}$v$})\bigg\} \\
&\hspace*{0.5cm} +\delta \bigg\{ \partial ^{m+2}_{s}\mbox{\mathversion{bold}$v$}  
+ 2(\mbox{\mathversion{bold}$v$}_{s}\cdot \partial ^{m+1}_{s}\mbox{\mathversion{bold}$v$})\mbox{\mathversion{bold}$v$} + \mbox{\mathversion{bold}$z$}_{m} \bigg\}
+ \mbox{\mathversion{bold}$w$}_{m},
\end{aligned}
\end{eqnarray*}
where \( \mbox{\mathversion{bold}$z$}_{m}\) and \( \mbox{\mathversion{bold}$w$}_{m}\) are terms that contain derivatives of \( \mbox{\mathversion{bold}$v$}\) up to order \( m\)
and are independent of \( \delta \). 
We estimate the solution in the following way.
\begin{eqnarray*}
\begin{aligned}
\frac{1}{2}\frac{{\rm d}}{{\rm d}t}\| \partial ^{m+1}_{s}\mbox{\mathversion{bold}$v$}\| ^{2} 
&=-(\partial ^{m}_{s}\mbox{\mathversion{bold}$v$}_{t},\partial ^{m+2}_{s}\mbox{\mathversion{bold}$v$})- \big( \partial ^{m}_{s}\mbox{\mathversion{bold}$v$}_{t}\cdot
\partial ^{m+1}_{s}\mbox{\mathversion{bold}$v$}\big)(0) \\
&=-m(\mbox{\mathversion{bold}$v$}_{s}\times \partial ^{m+1}_{s}\mbox{\mathversion{bold}$v$},\partial ^{m+2}_{s}\mbox{\mathversion{bold}$v$})
- \big( \partial ^{m}_{s}\mbox{\mathversion{bold}$v$}_{t}\cdot
\partial ^{m+1}_{s}\mbox{\mathversion{bold}$v$}\big)(0)
- \alpha \bigg\{ (\partial ^{m+3}_{s}\mbox{\mathversion{bold}$v$},\partial ^{m+2}_{s}\mbox{\mathversion{bold}$v$}) \\
&+ \frac{3}{2}(m+1)(\partial^{m+1}_{s}\mbox{\mathversion{bold}$v$}\times (\mbox{\mathversion{bold}$v$}\times \mbox{\mathversion{bold}$v$}_{ss}), 
\partial ^{m+2}_{s}\mbox{\mathversion{bold}$v$}) + \frac{3}{2}(m+1)(\mbox{\mathversion{bold}$v$}_{ss}\times (\mbox{\mathversion{bold}$v$}\times 
\partial ^{m+1}_{s}\mbox{\mathversion{bold}$v$}), \partial ^{m+2}_{s}\mbox{\mathversion{bold}$v$}) \\
&+ \frac{3}{2}(\mbox{\mathversion{bold}$v$}_{s}\times (\mbox{\mathversion{bold}$v$}\times \partial ^{m+2}_{s}\mbox{\mathversion{bold}$v$}),
\partial ^{m+2}_{s}\mbox{\mathversion{bold}$v$}) + \frac{3m}{2}(\mbox{\mathversion{bold}$v$}_{s}\times 
(\mbox{\mathversion{bold}$v$}_{s}\times \partial ^{m+1}_{s}\mbox{\mathversion{bold}$v$}), \partial ^{m+2}_{s}\mbox{\mathversion{bold}$v$}) \bigg\} \\
&- \delta \bigg\{ (\partial ^{m+2}_{s}\mbox{\mathversion{bold}$v$}, \partial ^{m+2}_{s}\mbox{\mathversion{bold}$v$})
+ 2((\mbox{\mathversion{bold}$v$}_{s}\cdot \partial ^{m+1}_{s}\mbox{\mathversion{bold}$v$})\mbox{\mathversion{bold}$v$},\partial ^{m+2}_{s}\mbox{\mathversion{bold}$v$})
+(\mbox{\mathversion{bold}$z$}_{m}, \partial ^{m+2}_{s}\mbox{\mathversion{bold}$v$} )\bigg\} 
- ( \mbox{\mathversion{bold}$w$}_{m}, \partial ^{m+2}_{s}\mbox{\mathversion{bold}$v$}).
\end{aligned}
\end{eqnarray*}
Each term is estimated by using the fact that \( |\mbox{\mathversion{bold}$v$}|=1\),  (\ref{par}), and (\ref{par2}). The usage of these properties is sometimes hard to notice and 
somewhat complicated, so we
give a detailed calculation for such term even though the calculus itself is elementary. Set \( m_{*}:=\max \{3, m-3\} \). First we have
\begin{eqnarray*}
\begin{aligned}
-m(\mbox{\mathversion{bold}$v$}_{s}\times \partial ^{m+1}_{s}\mbox{\mathversion{bold}$v$},\partial ^{m+2}_{s}\mbox{\mathversion{bold}$v$})
&=m (\mbox{\mathversion{bold}$v$}_{s}\times \partial ^{m+2}_{s}\mbox{\mathversion{bold}$v$},\partial ^{m+1}_{s}\mbox{\mathversion{bold}$v$})\\
&= -m((\mbox{\mathversion{bold}$v$}\cdot \partial ^{m+2}_{s}\mbox{\mathversion{bold}$v$})\mbox{\mathversion{bold}$v$}\times \mbox{\mathversion{bold}$v$}_{s},
\partial ^{m+1}_{s}\mbox{\mathversion{bold}$v$})+m([(\mbox{\mathversion{bold}$v$}\times \mbox{\mathversion{bold}$v$}_{s})\cdot \partial ^{m+2}_{s}\mbox{\mathversion{bold}$v$}]\mbox{\mathversion{bold}$v$},
\partial ^{m+1}_{s}\mbox{\mathversion{bold}$v$}) \\
&=\frac{1}{2}m\sum ^{m+1}_{j=1}
\left(
\begin{array}{c}
m+2\\
j
\end{array}\right)
((\partial ^{j}_{s}\mbox{\mathversion{bold}$v$}\cdot \partial ^{m+2-j}_{s}\mbox{\mathversion{bold}$v$})\mbox{\mathversion{bold}$v$}\times
\mbox{\mathversion{bold}$v$}_{s}, \partial ^{m+1}_{s}\mbox{\mathversion{bold}$v$}) \\
& \ \ \ - m([(\mbox{\mathversion{bold}$v$}\times \mbox{\mathversion{bold}$v$}_{ss})\cdot \partial ^{m+1}_{s}\mbox{\mathversion{bold}$v$}]\mbox{\mathversion{bold}$v$},
\partial ^{m+1}_{s}\mbox{\mathversion{bold}$v$})
- m ([(\mbox{\mathversion{bold}$v$}\times \mbox{\mathversion{bold}$v$}_{s})\cdot \partial ^{m+1}_{s}\mbox{\mathversion{bold}$v$}]\mbox{\mathversion{bold}$v$}_{s},
\partial ^{m+1}_{s}\mbox{\mathversion{bold}$v$}) \\
& \ \ \ +\frac{m}{2}\sum^{m+1}_{j=1}
\left(
\begin{array}{c}
m+2\\
j
\end{array}\right)
([(\mbox{\mathversion{bold}$v$}\times \mbox{\mathversion{bold}$v$}_{s})\cdot \partial ^{m+1}_{s}\mbox{\mathversion{bold}$v$}]\partial ^{j}_{s}\mbox{\mathversion{bold}$v$},
\partial ^{m+2-j}_{s}\mbox{\mathversion{bold}$v$}) \\
&\leq C\| \mbox{\mathversion{bold}$v$}_{s}\|^{2}_{m},
\end{aligned}
\end{eqnarray*}
where \( C\) depends on \( \|\mbox{\mathversion{bold}$v$}_{s}\|_{m_{*}}\).
Next we have
\begin{eqnarray*}
\begin{aligned}
\frac{3}{2}(m+1)&(\partial ^{m+1}_{s}\mbox{\mathversion{bold}$v$}\times (\mbox{\mathversion{bold}$v$}\times \mbox{\mathversion{bold}$v$}_{ss}),\partial ^{m+2}_{s}\mbox{\mathversion{bold}$v$})\\
&= \frac{3}{2}(m+1)\bigg\{ ((\mbox{\mathversion{bold}$v$}_{ss}\cdot \partial ^{m+1}_{s}\mbox{\mathversion{bold}$v$})\mbox{\mathversion{bold}$v$},\partial ^{m+2}_{s}\mbox{\mathversion{bold}$v$})
-((\mbox{\mathversion{bold}$v$}\cdot \partial ^{m+1}_{s}\mbox{\mathversion{bold}$v$})\mbox{\mathversion{bold}$v$}_{ss},\partial ^{m+2}_{s}\mbox{\mathversion{bold}$v$})\bigg\} \\
&=\frac{3}{2}(m+1)\bigg\{ -\frac{1}{2}\sum^{m+1}_{j=1}
\left( 
\begin{array}{c}
m+2\\
j
\end{array}\right)
((\mbox{\mathversion{bold}$v$}_{ss}\cdot \partial ^{m+1}_{s}\mbox{\mathversion{bold}$v$})\partial ^{j}_{s}\mbox{\mathversion{bold}$v$},\partial ^{m+2-j}_{s}\mbox{\mathversion{bold}$v$})\\
& \ \ \ \ +( (\mbox{\mathversion{bold}$v$}_{s}\cdot \partial ^{m+1}_{s}\mbox{\mathversion{bold}$v$})\mbox{\mathversion{bold}$v$}_{ss},\partial^{m+1}_{s}\mbox{\mathversion{bold}$v$})
+((\mbox{\mathversion{bold}$v$}\cdot \partial ^{m+2}_{s}\mbox{\mathversion{bold}$v$})\mbox{\mathversion{bold}$v$}_{ss},\partial ^{m+1}_{s}\mbox{\mathversion{bold}$v$}) \\
& \ \ \ \ +((\mbox{\mathversion{bold}$v$}\cdot \partial ^{m+1}_{s}\mbox{\mathversion{bold}$v$})\mbox{\mathversion{bold}$v$}_{sss},\partial ^{m+1}_{s}\mbox{\mathversion{bold}$v$})
+\big( (\mbox{\mathversion{bold}$v$}\cdot\partial ^{m+1}_{s}\mbox{\mathversion{bold}$v$})\mbox{\mathversion{bold}$v$}_{ss}\cdot \partial ^{m+1}_{s}\mbox{\mathversion{bold}$v$}
\big) (0) \bigg\} \\
&\leq C\big( \|\mbox{\mathversion{bold}$v$}_{s}\|^{2}_{m} + |\partial ^{m+1}_{s}\mbox{\mathversion{bold}$v$}(0)|^{2}\big),
\end{aligned}
\end{eqnarray*}
where \( C\) depends on \( \| \mbox{\mathversion{bold}$v$}_{s}\|_{m_{*}}\).
We continue with
\begin{eqnarray*}
\begin{aligned}
\frac{3}{2}(m+1)&(\mbox{\mathversion{bold}$v$}_{ss}\times (\mbox{\mathversion{bold}$v$}\times \partial ^{m+1}_{s}),\partial ^{m+2}_{s}\mbox{\mathversion{bold}$v$})\\
&= \frac{3}{2}(m+1) \bigg\{ ((\mbox{\mathversion{bold}$v$}_{ss}\cdot \partial ^{m+1}_{s}\mbox{\mathversion{bold}$v$})\mbox{\mathversion{bold}$v$},\partial ^{m+2}_{s}\mbox{\mathversion{bold}$v$})
-((\mbox{\mathversion{bold}$v$}\cdot \mbox{\mathversion{bold}$v$}_{ss})\partial ^{m+1}_{s}\mbox{\mathversion{bold}$v$},\partial ^{m+2}_{s}\mbox{\mathversion{bold}$v$})\bigg\} \\
&= \frac{3}{2}(m+1)\bigg\{ -\frac{1}{2}\sum ^{m+1}_{j=0}
\left( 
\begin{array}{c}
m+2\\
j
\end{array}\right)
((\mbox{\mathversion{bold}$v$}_{ss}\cdot \partial ^{m+1}_{s}\mbox{\mathversion{bold}$v$})\partial ^{j}_{s}\mbox{\mathversion{bold}$v$},\partial ^{m+2-j}_{s}\mbox{\mathversion{bold}$v$})\\
& \ \ \ \ + ((\mbox{\mathversion{bold}$v$}\cdot \mbox{\mathversion{bold}$v$}_{ss})_{s}\partial ^{m+1}_{s}\mbox{\mathversion{bold}$v$},\partial ^{m+1}_{s}\mbox{\mathversion{bold}$v$})
-\big( (\mbox{\mathversion{bold}$v$}\cdot \mbox{\mathversion{bold}$v$}_{ss})|\partial ^{m+1}_{s}\mbox{\mathversion{bold}$v$}|^{2}\big)(0)\bigg\} \\
&\leq C\|\mbox{\mathversion{bold}$v$}_{s}\|^{2}_{m},
\end{aligned}
\end{eqnarray*}
where, again, \( C\) depends on \( \|\mbox{\mathversion{bold}$v$}_{s}\|_{m_{*}}\). From here on, it will be assumed that generic constants \( C\) depend on 
\( \|\mbox{\mathversion{bold}$v$}_{s}\|_{m_{*}}\) unless explicitly mentioned otherwise.
We calculate furthermore
\begin{eqnarray*}
\begin{aligned}
\frac{3}{2}&(\mbox{\mathversion{bold}$v$}_{s}\times (\mbox{\mathversion{bold}$v$}\times \partial ^{m+2}_{s}\mbox{\mathversion{bold}$v$}),\partial ^{m+2}_{s}\mbox{\mathversion{bold}$v$})
= \frac{3}{2}((\mbox{\mathversion{bold}$v$}_{s}\cdot \partial ^{m+2}_{s})\mbox{\mathversion{bold}$v$}, \partial ^{m+2}_{s}\mbox{\mathversion{bold}$v$}) \\
&= \frac{3}{2}\bigg\{ -\frac{1}{2}((\mbox{\mathversion{bold}$v$}_{s}\cdot \partial ^{m+2}_{s}\mbox{\mathversion{bold}$v$})\mbox{\mathversion{bold}$v$}_{s},
\partial ^{m+1}_{s}\mbox{\mathversion{bold}$v$})-\frac{1}{2}\sum ^{m}_{j=2}
\left(
\begin{array}{c}
m+2\\
j
\end{array}\right)
((\mbox{\mathversion{bold}$v$}_{s}\cdot \partial ^{m+2}_{s}\mbox{\mathversion{bold}$v$})\partial ^{j}_{s}\mbox{\mathversion{bold}$v$},\partial ^{m+2-j}_{s}\mbox{\mathversion{bold}$v$})
\bigg\} \\
&= \frac{3}{4}\bigg\{ ((\mbox{\mathversion{bold}$v$}_{ss}\cdot \partial ^{m+1}_{s}\mbox{\mathversion{bold}$v$})\mbox{\mathversion{bold}$v$}_{s},
\partial ^{m+1}_{s}\mbox{\mathversion{bold}$v$}) + 
\sum ^{m}_{j=2}
\left(
\begin{array}{c}
m+2\\
j
\end{array}\right)
\big[ ((\mbox{\mathversion{bold}$v$}_{ss}\cdot \partial ^{m+1}_{s}\mbox{\mathversion{bold}$v$})\partial ^{j}_{s}\mbox{\mathversion{bold}$v$},
\partial ^{m+2-j}_{s}\mbox{\mathversion{bold}$v$}) \\
&\ \ \ + ((\mbox{\mathversion{bold}$v$}_{s}\cdot \partial ^{m+1}_{s}\mbox{\mathversion{bold}$v$})\partial ^{j+1}_{s}\mbox{\mathversion{bold}$v$},
\partial ^{m+2-j}_{s}\mbox{\mathversion{bold}$v$}) + ((\mbox{\mathversion{bold}$v$}_{s}\cdot \partial ^{m+1}_{s}\mbox{\mathversion{bold}$v$})\partial ^{j}_{s}\mbox{\mathversion{bold}$v$},
\partial ^{m+3-j}_{s}\mbox{\mathversion{bold}$v$}) \big] \bigg\}\\
&\leq C\|\mbox{\mathversion{bold}$v$}_{s}\|^{2}_{m},
\end{aligned}
\end{eqnarray*}
\begin{eqnarray*}
\begin{aligned}
\frac{3m}{2}(\mbox{\mathversion{bold}$v$}_{s}\times (\mbox{\mathversion{bold}$v$}_{s}\times \partial ^{m+1}_{s}\mbox{\mathversion{bold}$v$}),
\partial ^{m+2}_{s}\mbox{\mathversion{bold}$v$}) &= 
\frac{3m}{2}\bigg\{ ((\mbox{\mathversion{bold}$v$}_{s}\cdot \partial ^{m+1}_{s}\mbox{\mathversion{bold}$v$})\mbox{\mathversion{bold}$v$}_{s},
\partial ^{m+2}_{s}\mbox{\mathversion{bold}$v$})-((\mbox{\mathversion{bold}$v$}_{s}\cdot \mbox{\mathversion{bold}$v$}_{s})\partial ^{m+1}_{s}\mbox{\mathversion{bold}$v$},
\partial ^{m+2}_{s}\mbox{\mathversion{bold}$v$})\bigg\} \\
&=\frac{3m}{2}\bigg\{ -((\mbox{\mathversion{bold}$v$}_{ss}\cdot \partial ^{m+1}_{s}\mbox{\mathversion{bold}$v$})\mbox{\mathversion{bold}$v$}_{s},
\partial ^{m+1}_{s}\mbox{\mathversion{bold}$v$}) + \frac{1}{2}((\mbox{\mathversion{bold}$v$}_{s}\cdot \mbox{\mathversion{bold}$v$}_{s})_{s}\partial ^{m+1}_{s}\mbox{\mathversion{bold}$v$},
\partial ^{m+1}_{s}\mbox{\mathversion{bold}$v$}) \bigg\} \\
&\leq C\|\mbox{\mathversion{bold}$v$}_{s}\|^{2}_{m}.
\end{aligned}
\end{eqnarray*}
Next, we estimate the boundary terms.
\begin{eqnarray*}
\begin{aligned}
\partial ^{m}_{s}\mbox{\mathversion{bold}$v$}_{t}\cdot \partial ^{m+1}_{s}\mbox{\mathversion{bold}$v$} &=
(\mbox{\mathversion{bold}$v$}\times \partial ^{m+2}_{s}\mbox{\mathversion{bold}$v$})\cdot \partial ^{m+1}_{s}\mbox{\mathversion{bold}$v$}
+ \alpha \bigg\{ \partial ^{m+3}_{s}\mbox{\mathversion{bold}$v$}\cdot \partial ^{m+1}_{s} \mbox{\mathversion{bold}$v$}
+ \frac{3}{2}[\partial ^{m+2}_{s}\mbox{\mathversion{bold}$v$}\times (\mbox{\mathversion{bold}$v$}\times \mbox{\mathversion{bold}$v$}_{s})]\cdot 
\partial ^{m+1}_{s}\mbox{\mathversion{bold}$v$} \\
& \ \ \ +\frac{3}{2}(m+1)[\mbox{\mathversion{bold}$v$}_{ss}\times ( \mbox{\mathversion{bold}$v$}\times \partial ^{m+1}_{s}\mbox{\mathversion{bold}$v$})]\cdot \partial ^{m+1}_{s}\mbox{\mathversion{bold}$v$} + \frac{3}{2}[\mbox{\mathversion{bold}$v$}_{s}\times (\mbox{\mathversion{bold}$v$}\times \partial ^{m+2}_{s}\mbox{\mathversion{bold}$v$})]
\cdot \partial ^{m+1}_{s}\mbox{\mathversion{bold}$v$} \\
& \ \ \ + \frac{3m}{2}[\mbox{\mathversion{bold}$v$}_{s}\times (\mbox{\mathversion{bold}$v$}_{s}\times \partial ^{m+1}_{s}\mbox{\mathversion{bold}$v$})]\cdot
\partial ^{m+1}_{s}\mbox{\mathversion{bold}$v$}\bigg\} 
+ \delta \bigg\{ \partial ^{m+2}_{s}\mbox{\mathversion{bold}$v$}\cdot \partial ^{m+1}_{s}\mbox{\mathversion{bold}$v$} \\
& \ \ \ + 2(\mbox{\mathversion{bold}$v$}_{s}\cdot \partial ^{m+1}_{s}\mbox{\mathversion{bold}$v$})(\mbox{\mathversion{bold}$v$}\cdot \partial ^{m+1}_{s}\mbox{\mathversion{bold}$v$})
+ \mbox{\mathversion{bold}$z$}_{m}\cdot \partial ^{m+1}_{s}\mbox{\mathversion{bold}$v$} \bigg\} + \mbox{\mathversion{bold}$w$}_{m}\cdot 
\partial ^{m+1}_{s}\mbox{\mathversion{bold}$v$} ,
\end{aligned}
\end{eqnarray*}
thus we have
\begin{eqnarray*}
\begin{aligned}
\big( \partial ^{m}_{s}\mbox{\mathversion{bold}$v$}_{t}\cdot \partial ^{m+1}_{s}\mbox{\mathversion{bold}$v$}\big) (0) &= 
[(\mbox{\mathversion{bold}$v$}\times \partial ^{m+2}_{s}\mbox{\mathversion{bold}$v$})\cdot \partial ^{m+1}_{s}\mbox{\mathversion{bold}$v$}](0)\\
& \ \ \ +\alpha \bigg\{ \partial ^{m+3}_{s}\mbox{\mathversion{bold}$v$}\cdot \partial ^{m+1}_{s}\mbox{\mathversion{bold}$v$} 
+ \frac{3}{2}(m+1)[\mbox{\mathversion{bold}$v$}_{ss}\times ( \mbox{\mathversion{bold}$v$}\times \partial ^{m+1}_{s}\mbox{\mathversion{bold}$v$})]
\cdot \partial ^{m+1}_{s}\mbox{\mathversion{bold}$v$}\bigg\}(0) \\
& \ \ \ + \delta \bigg\{ \partial ^{m+2}_{s}\mbox{\mathversion{bold}$v$}\cdot \partial ^{m+1}_{s}\mbox{\mathversion{bold}$v$}
+ \mbox{\mathversion{bold}$z$}_{m}\cdot \partial ^{m+1}_{s}\mbox{\mathversion{bold}$v$}\bigg\} (0)
+ (\mbox{\mathversion{bold}$w$}_{m}\cdot \partial ^{m+1}_{s}\mbox{\mathversion{bold}$v$})(0).
\end{aligned}
\end{eqnarray*}
Again, we estimate each term separately.
\begin{eqnarray*}
\begin{aligned}
\big\{ [\mbox{\mathversion{bold}$v$}_{ss}\times (\mbox{\mathversion{bold}$v$}\times \partial ^{m+1}_{s}\mbox{\mathversion{bold}$v$})]
\cdot \partial ^{m+1}_{s}\mbox{\mathversion{bold}$v$}\big\} (0) 
&= \left. \bigg\{  (\mbox{\mathversion{bold}$v$}_{ss}\cdot \partial ^{m+1}_{s}\mbox{\mathversion{bold}$v$})\mbox{\mathversion{bold}$v$}
-(\mbox{\mathversion{bold}$v$}\cdot \mbox{\mathversion{bold}$v$}_{ss})\partial ^{m+1}_{s}\mbox{\mathversion{bold}$v$}\bigg\} \cdot 
\partial ^{m+1}_{s}\mbox{\mathversion{bold}$v$} \right| _{s=0} \\
&=  -\frac{1}{2}\sum ^{m}_{j=1}
\left( 
\begin{array}{c}
m+1\\
j
\end{array}\right)
(\mbox{\mathversion{bold}$v$}_{ss}\cdot \partial ^{m+1}_{s}\mbox{\mathversion{bold}$v$})
(\partial^{j}_{s}\mbox{\mathversion{bold}$v$}\cdot \partial ^{m+1-j}_{s}\mbox{\mathversion{bold}$v$}) \bigg| _{s=0} \\
&\leq C\big( \| \mbox{\mathversion{bold}$v$}_{s}\|^{2}_{m}+ |\partial ^{m+1}_{s}\mbox{\mathversion{bold}$v$}(0)|^{2}\big)
\end{aligned}
\end{eqnarray*}
holds.
Combining the estimates yields
\begin{eqnarray*}
\begin{aligned}
\frac{1}{2}\frac{{\rm d}}{{\rm d}t}\| \partial ^{m+1}_{s}\mbox{\mathversion{bold}$v$}\|^{2} &+ \frac{1}{2}|\alpha | |\partial ^{m+2}_{s}\mbox{\mathversion{bold}$v$}(0)|^{2} \\
&\leq C\big( \|\mbox{\mathversion{bold}$v$}_{s}\|^{2}_{m} + |\partial ^{m+1}_{s}\mbox{\mathversion{bold}$v$}(0)|^{2}\big) 
+ (\mbox{\mathversion{bold}$v$}\times \partial ^{m+2}_{s}\mbox{\mathversion{bold}$v$})\cdot \partial ^{m+1}_{s}\mbox{\mathversion{bold}$v$} \bigg| _{s=0}
+ \alpha \partial ^{m+3}_{s}\mbox{\mathversion{bold}$v$}\cdot \partial ^{m+1}_{s}\mbox{\mathversion{bold}$v$} \bigg|_{s=0} \\
& \ \ \ + \delta (\partial ^{m+2}_{s}\mbox{\mathversion{bold}$v$}\cdot \partial ^{m+1}_{s}\mbox{\mathversion{bold}$v$})\bigg|_{s=0}.
\end{aligned}
\end{eqnarray*}
On the other hand, from the boundary condition we see that the solution satisfies 
\( \partial ^{j}_{t}\mbox{\mathversion{bold}$v$}_{s}(0,t) = \mbox{\mathversion{bold}$0$}\) for any \( j\) with \( 0\leq j\leq k\).
Substituting the equation to convert all \( t\) derivatives into \( s\) derivatives yields
\( \alpha ^{j} ( \partial ^{3j+1}_{s}\mbox{\mathversion{bold}$v$})(0,t) = \mbox{\mathversion{bold}$F$}(\mbox{\mathversion{bold}$v$},\mbox{\mathversion{bold}$v$}_{s},
\cdots , \partial ^{3j}_{s}\mbox{\mathversion{bold}$v$})(0,t)\), i.e. boundary terms with \( (3j+1) \)-th order derivative can be expressed 
in terms of boundary terms with derivatives up to order \( 3j\). By choosing
\( m=3j+1 \), we have
\begin{eqnarray}
\begin{aligned}
\frac{1}{2}\frac{{\rm d}}{{\rm d}t}\| \partial ^{3j+2}_{s}\mbox{\mathversion{bold}$v$}\|^{2} &+ \frac{|\alpha |}{2}|\partial ^{3j+3}_{s}\mbox{\mathversion{bold}$v$}(0)|^{2}\\
&\leq C\big( \|\mbox{\mathversion{bold}$v$}_{s}\|^{2}_{1+3j} + |\partial ^{3j+2}_{s}\mbox{\mathversion{bold}$v$}(0)|^{2}\big) 
+ (\mbox{\mathversion{bold}$v$}\times \partial ^{3j+3}_{s}\mbox{\mathversion{bold}$v$})\cdot \partial ^{3j+2}_{s}\mbox{\mathversion{bold}$v$} \bigg| _{s=0} \\
& \ \ \ + \alpha \partial ^{3(j+1)+1}_{s}\mbox{\mathversion{bold}$v$}\cdot \partial ^{3j+2}_{s}\mbox{\mathversion{bold}$v$} \bigg|_{s=0} 
+ \delta (\partial ^{3j+3}_{s}\mbox{\mathversion{bold}$v$}\cdot \partial ^{3j+2}_{s}\mbox{\mathversion{bold}$v$})\bigg|_{s=0}\\
&\leq C\big( \|\mbox{\mathversion{bold}$v$}_{s}\|^{2}_{1+3j} + |\partial ^{3j+2}_{s}\mbox{\mathversion{bold}$v$}(0)|^{2}\big) 
+ (\mbox{\mathversion{bold}$v$}\times \partial ^{3j+3}_{s}\mbox{\mathversion{bold}$v$})\cdot \partial ^{3j+2}_{s}\mbox{\mathversion{bold}$v$} \bigg| _{s=0} \\
& \ \ \ +  |\alpha | \big| \partial ^{3j+3}_{s}\mbox{\mathversion{bold}$v$}\big| \big| \partial ^{3j+2}_{s}\mbox{\mathversion{bold}$v$}\big| \bigg|_{s=0} 
+ \delta \big| \partial ^{3j+3}_{s}\mbox{\mathversion{bold}$v$}\big| \big| \partial ^{3j+2}_{s}\mbox{\mathversion{bold}$v$}\big| \bigg|_{s=0}.
\label{uest}
\end{aligned}
\end{eqnarray}
By a similar estimate, we can show that 
\begin{eqnarray*}
\begin{aligned}
\frac{1}{2}\frac{{\rm d}}{{\rm d}t}\| \partial ^{j+1}_{s}\mbox{\mathversion{bold}$v$}\|^{2}
&+\frac{|\alpha|}{2}|\partial ^{j+2}_{s}\mbox{\mathversion{bold}$v$}(0)|^{2} \\
&\leq C\big( \| \mbox{\mathversion{bold}$v$}_{s}\|^{3}_{3} + |\partial ^{j+1}_{s}\mbox{\mathversion{bold}$v$}(0)|^{2}\big)
+(\mbox{\mathversion{bold}$v$}\times \partial ^{j+2}_{s}\mbox{\mathversion{bold}$v$})\cdot \partial ^{j+1}_{s}\mbox{\mathversion{bold}$v$}\bigg|_{s=0} \\
& \ \ \ + |\alpha | \big| \partial ^{j+3}_{s}\mbox{\mathversion{bold}$v$}\big| \big| \partial ^{j+1}_{s}\mbox{\mathversion{bold}$v$}\big|\bigg|_{s=0}
+\delta \big| \partial ^{j+1}_{s}\mbox{\mathversion{bold}$v$}\big| \big| \partial ^{j+2}_{s}\mbox{\mathversion{bold}$v$}\big| \bigg|_{s=0}
\end{aligned}
\end{eqnarray*}
holds for \( j=0,1,2,3\). Thus, for \( \eta >0\) we have
\begin{eqnarray*}
\frac{1}{2}\frac{{\rm d}}{{\rm d}t}\| \partial ^{j+1}_{s}\mbox{\mathversion{bold}$v$}\|^{2}+\frac{|\alpha |}{4}\big| \partial ^{j+2}_{s}\mbox{\mathversion{bold}$v$}(0)\big|^{2}
\leq C\big( \| \mbox{\mathversion{bold}$v$}_{s}\|^{2}_{3} + \big| \partial ^{j+1}_{s}\mbox{\mathversion{bold}$v$}(0)\big|^{2}\big)
+ \eta \big| \partial ^{j+3}_{s}\mbox{\mathversion{bold}$v$}(0)\big|^{2}
\end{eqnarray*}
for \( j=1,2\), and
\begin{eqnarray*}
\begin{aligned}
\frac{1}{2}\frac{{\rm d}}{{\rm d}t}\| \partial ^{4}_{s}\mbox{\mathversion{bold}$v$}\|^{2} + \frac{|\alpha |}{4}\big| \partial^{5}_{s}\mbox{\mathversion{bold}$v$}(0)\big|^{2}
&\leq C\| \mbox{\mathversion{bold}$v$}_{s}\|^{3}_{3} + C\big| \partial ^{4}_{s}\mbox{\mathversion{bold}$v$}(0)\big|^{2}, \\
\frac{1}{2}\frac{{\rm d}}{{\rm d}t}\|\mbox{\mathversion{bold}$v$}_{s}\|^{2} + \frac{|\alpha |}{4}\big| \partial ^{2}_{s}\mbox{\mathversion{bold}$v$}(0)\big|^{2}
&\leq C\| \mbox{\mathversion{bold}$v$}_{s}\|^{3}_{3}.
\end{aligned}
\end{eqnarray*}
Here, \( C\) depends on \( \eta \) and \( C_{*}\). 
By taking a linear combination of the above estimates, we arrive at
\begin{eqnarray*}
\frac{1}{2}\frac{{\rm d}}{{\rm d}t}\|\mbox{\mathversion{bold}$v$}_{s}\|^{2}_{3} \leq C\|\mbox{\mathversion{bold}$v$}_{s}\|^{3}_{3},
\end{eqnarray*}
where \( C\) depends on \( C_{*}\) but not on \( \delta \in (0, \delta_{*}] \). Utilizing the comparison principle as before, the above estimate gives a
time-local uniform estimate in \( C\big( [0,T_{1}];H^{3}(\mathbf{R}_{+})\big) \) for some \( T_{1}\in (0,T] \). 
From the \( H^{3}\) estimate and (\ref{uest}), we can derive the uniform estimate in \( C\big( [0,T_{1}]; H^{1+3k}(\mathbf{R}_{+})\big) \)
in the same manner. Here, \( T_{1}\) is determined from the \( H^{3}\) estimate and only depends on \( \|\mbox{\mathversion{bold}$v$}_{0s}\|_{3}\).
This finishes the proof of the proposition. \hfill \( \Box \)


\subsection{Taking the limit \( \delta \rightarrow +0\)}

Now, we are ready to take the limit \( \delta \rightarrow +0\). We prove the 
existence theorem for the case \( \alpha <0\).

{\it Proof of Theorem {\rm \ref{st2}}}.
From the way we constructed \( \mbox{\mathversion{bold}$v$}^{\delta}_{0} \), the 
following holds.
\( \mbox{\mathversion{bold}$v$}^{\delta }_{0s}\rightarrow \mbox{\mathversion{bold}$v$}_{0s}\) in \( H^{1+3k}(\mathbf{R}_{+})\) and \( \mbox{\mathversion{bold}$v$}^{\delta}_{0}
\rightarrow \mbox{\mathversion{bold}$v$}_{0}\) in \( L^{\infty}(\mathbf{R}_{+})\) as \( \delta \rightarrow +0\).
Thus, by taking \( \delta _{*}>0\) smaller if necessary, we have \( \| \mbox{\mathversion{bold}$v$}^{\delta}_{0s}\|_{1+3k}
\leq 2\| \mbox{\mathversion{bold}$v$}_{0s}\|_{1+3k}\) for any \( \delta \in (0,\delta _{*}]\). 
For such \( \delta \), the solution \( \mbox{\mathversion{bold}$v$}^{\delta}\) constructed in section 5.1 with initial datum 
\( \mbox{\mathversion{bold}$v$}^{\delta }_{0}\) satisfies the assumptions of Lemma \ref{unit} and Proposition \ref{basest} with
\( M=2\| \mbox{\mathversion{bold}$v$}_{0s}\|_{1+3k}\),
i.e., the solution satisfies \( \big| \mbox{\mathversion{bold}$v$}^{\delta}\big| =1\) and a uniform estimate in 
\( C\big( [0,T]; H^{1+3k}(\mathbf{R}_{+})\big) \) for some \( T>0\) holds.
For any \( \delta , \delta ^{\prime}\in (0,\delta _{*}]\), we set
\( \mbox{\mathversion{bold}$V$}:=
\mbox{\mathversion{bold}$v$}^{\delta ^{\prime}}-\mbox{\mathversion{bold}$v$}^{\delta }
-(\mbox{\mathversion{bold}$v$}^{\delta ^{\prime}}_{0}-\mbox{\mathversion{bold}$v$}^{\delta }_{0})\), then \( \mbox{\mathversion{bold}$V$}
\) satisfies
\begin{eqnarray*}
\left\{
\begin{array}{ll}
\mbox{\mathversion{bold}$V$}_{t} =\alpha \mbox{\mathversion{bold}$V$}_{sss}
+ \mbox{\mathversion{bold}$v$}^{\delta ^{\prime}}\times
\mbox{\mathversion{bold}$V$}_{ss}+3\alpha 
\mbox{\mathversion{bold}$V$}_{ss}\times (\mbox{\mathversion{bold}$v$}^{\delta ^{\prime}}
\times \mbox{\mathversion{bold}$v$}^{\delta ^{\prime}}_{s})
 
+\delta ^{\prime}(\mbox{\mathversion{bold}$v$}^{\delta ^{\prime}}_{ss}
+ |\mbox{\mathversion{bold}$v$}^{\delta ^{\prime}}_{s}|^{2}
\mbox{\mathversion{bold}$v$}^{\delta ^{\prime}} ), & \ \\
\hspace*{8.5cm} -\delta (\mbox{\mathversion{bold}$v$}^{\delta }_{ss}
+ |\mbox{\mathversion{bold}$v$}^{\delta }_{s}|^{2}\mbox{\mathversion{bold}$v$}^{\delta}) 
+ \mbox{\mathversion{bold}$F$}, & s>0, t>0, \\
\mbox{\mathversion{bold}$V$}(s,0)=\mbox{\mathversion{bold}$0$}, & s>0, \\
\mbox{\mathversion{bold}$V$}_{s}(0,t)=\mbox{\mathversion{bold}$0$}, & t>0,
\end{array}\right.
\end{eqnarray*}
where \( \mbox{\mathversion{bold}$F$}\) is the collection of terms that are lower order in \( \mbox{\mathversion{bold}$V$}\) 
and depends linearly on \( \mbox{\mathversion{bold}$V$}_{0}:= \mbox{\mathversion{bold}$v$}^{\delta ^{\prime}}_{0}-\mbox{\mathversion{bold}$v$}^{\delta }_{0}\).
By a standard energy method, we have
\begin{eqnarray*}
\begin{aligned}
\frac{1}{2}\frac{{\rm d}}{{\rm d}t}\|\mbox{\mathversion{bold}$V$}\|^{2}
&\leq \alpha \mbox{\mathversion{bold}$V$}(0)\cdot \mbox{\mathversion{bold}$V$}_{ss}(0)
+ C\| \mbox{\mathversion{bold}$V$}\|^{2}_{1} + C\big[ (\delta + \delta ^{\prime}) + \|\mbox{\mathversion{bold}$V$}_{0}\|^{2}_{L^{\infty}(\mathbf{R}_{+})} 
+ \|\mbox{\mathversion{bold}$V$}_{0s}\|^{2}_{2}\big] ,\\
\frac{1}{2}\frac{{\rm d}}{{\rm d}t}\|\mbox{\mathversion{bold}$V$}_{s}\|^{2}
&\leq -\frac{|\alpha |}{2}|\mbox{\mathversion{bold}$V$}_{ss}(0)|^{2}
+C\|\mbox{\mathversion{bold}$V$}\|^{2}_{1} + C\big[ (\delta + \delta ^{\prime}) + \|\mbox{\mathversion{bold}$V$}_{0}\|^{2}_{L^{\infty}(\mathbf{R}_{+})} 
+ \|\mbox{\mathversion{bold}$V$}_{0s}\|^{2}_{3}\big],
\end{aligned}
\end{eqnarray*}
where \( C\) is independent of \( \delta\) and \( \delta ^{\prime}\). Here, we have used
\( \mbox{\mathversion{bold}$v$}^{\delta^{\prime}}\cdot
(\mbox{\mathversion{bold}$v$}^{\delta ^{\prime}}-\mbox{\mathversion{bold}$v$}^{\delta })_{ss}=
-\mbox{\mathversion{bold}$v$}^{\delta^{\prime}}_{s}\cdot (\mbox{\mathversion{bold}$v$}^{\delta ^{\prime}}-\mbox{\mathversion{bold}$v$}^{\delta })_{s}
-(\mbox{\mathversion{bold}$v$}^{\delta ^{\prime}}-\mbox{\mathversion{bold}$v$}^{\delta })_{s}\cdot \mbox{\mathversion{bold}$v$}^{\delta}_{s}
-(\mbox{\mathversion{bold}$v$}^{\delta ^{\prime}}-\mbox{\mathversion{bold}$v$}^{\delta })\cdot \mbox{\mathversion{bold}$v$}^{\delta}_{ss}\),
which follows from the fact that \( |\mbox{\mathversion{bold}$v$}^{\delta}| 
=|\mbox{\mathversion{bold}$v$}^{\delta^{\prime}}|=1\). The above
estimate implies
\begin{eqnarray*}
\| \mbox{\mathversion{bold}$V$}\|^{2}_{1}\leq CT\big[ (\delta + \delta^{\prime}) + \|\mbox{\mathversion{bold}$V$}_{0}\|^{2}_{L^{\infty}(\mathbf{R}_{+})} 
+ \|\mbox{\mathversion{bold}$V$}_{0s}\|^{2}_{3}\big],
\end{eqnarray*}
where \( C\) is independent of \( \delta \) and \( \delta ^{\prime}\). Thus,
there is a \( \mbox{\mathversion{bold}$v$}\) such that
\( \mbox{\mathversion{bold}$v$}^{\delta}\rightarrow \mbox{\mathversion{bold}$v$}\) in \( C\big([0,T];L^{\infty}(\mathbf{R}_{+})\big) \)
and
\( \mbox{\mathversion{bold}$v$}^{\delta}_{s} \rightarrow \mbox{\mathversion{bold}$v$}_{s}\) in
\( C\big( [0,T]; L^{2}(\mathbf{R}_{+})\big) \). Combining these convergence with 
the uniform estimate, we have a solution \( \mbox{\mathversion{bold}$v$}\)
to (\ref{neg2}) such that
\( \mbox{\mathversion{bold}$v$}_{s}\in 
\bigcap ^{k}_{j=0} W^{j,\infty}\big( [0,T]; H^{1+3j}(\mathbf{R}_{+})\big) \)
and \( |\mbox{\mathversion{bold}$v$}|=1\).
Again, since we can approximate the initial datum by a smooth function, we have a solution 
\( \mbox{\mathversion{bold}$v$}\in \tilde{X}^{k}_{T}\), i.e. the continuity with respect to \( t\) can be recovered.
The uniform estimate that we obtained in the last subsection is essentially 
the energy estimate for \( \mbox{\mathversion{bold}$v$}\), thus the uniqueness of the 
solution follows. 
Based on this estimate, by taking a sequence of smooth initial datum \( \{ \mbox{\mathversion{bold}$v$}^{n}_{0}\} \) 
such that \( \mbox{\mathversion{bold}$v$}^{n}_{0}\rightarrow \mbox{\mathversion{bold}$v$}_{0}\) in 
\( L^{\infty}(\mathbf{R}_{+})\) and \( \mbox{\mathversion{bold}$v$}^{n}_{0s}\rightarrow \mbox{\mathversion{bold}$v$}_{0s}\) in \( H^{1+3k}(\mathbf{R}_{+})\)
and 
considering the convergence of the corresponding solution as \( n\rightarrow +\infty \) in the same manner as we did with \( \delta \rightarrow +0\), 
we have a solution \( \mbox{\mathversion{bold}$v$}\) satisfying \( \mbox{\mathversion{bold}$v$}_{s}\in 
\bigcap ^{k}_{j=0} W^{j,\infty}\big( [0,T]; H^{1+3j}(\mathbf{R}_{+})\big) \)
and \( |\mbox{\mathversion{bold}$v$}|=1\) with initial datum \( \mbox{\mathversion{bold}$v$}_{0}\) satisfying \( \mbox{\mathversion{bold}$v$}_{0s}\in H^{1+3k}(\mathbf{R}_{+})\)
and \( |\mbox{\mathversion{bold}$v$}_{0}|=1\).

This is the point where we are unable to recover the continuity in \( t\). This is because our problem cannot be solved reverse in time, which is required, 
to apply the standard method of recovering the continuity. So far, we have no new idea to solve this issue. 

Finally, as we mentioned in the introduction, we can construct \( \mbox{\mathversion{bold}$x$}\)
from \( \mbox{\mathversion{bold}$v$}\). This finishes the proof of 
Theorem \ref{st2}. \hfill \( \Box\)

%
%
%
%
\section{Construction of the Solution in the Case \( \alpha >0\) }
\setcounter{equation}{0}
\subsection{Existence of Solution}
We construct the solution in a similar manner as in the case \( \alpha <0\). For \( n\geq 2\), 
we define \( \mbox{\mathversion{bold}$v$}^{(n)}\) by
\begin{eqnarray*}
\left\{
\begin{array}{ll}
\mbox{\mathversion{bold}$v$}^{(n)}_{t}=\alpha \mbox{\mathversion{bold}$v$}^{(n)}_{sss} + {\rm A}(\mbox{\mathversion{bold}$v$}^{(n-1)},\partial_{s})\mbox{\mathversion{bold}$v$}^{(n)}
- \frac{3}{2}\alpha |\mbox{\mathversion{bold}$v$}^{(n-1)}_{s}|^{2}\mbox{\mathversion{bold}$v$}^{(n-1)}_{s} + \delta |\mbox{\mathversion{bold}$v$}^{(n-1)}_{s}|^{2}
\mbox{\mathversion{bold}$v$}^{(n-1)}, & s>0, t>0,\\
\mbox{\mathversion{bold}$v$}^{(n)}(s,0)=\mbox{\mathversion{bold}$v$}_{0}^{\delta ,R}(s), & s>0,\\
\mbox{\mathversion{bold}$v$}^{(n)}(0,t)=\mbox{\mathversion{bold}$e$}_{3}, & t>0, \\
\mbox{\mathversion{bold}$v$}^{(n)}_{s}(0,t)=\mbox{\mathversion{bold}$0$}, & t>0,
\end{array}\right.
\end{eqnarray*}
where \( \mbox{\mathversion{bold}$e$}_{3}= (0,0,1)\) and \( \mbox{\mathversion{bold}$v$}_{0}^{\delta ,R}\) is the 
initial datum cut-off at spatial infinity. Again, we define \( \mbox{\mathversion{bold}$v$}^{(1)}\) by
\begin{eqnarray*}
\mbox{\mathversion{bold}$v$}^{(1)}(s,t)= \mbox{\mathversion{bold}$v$}^{\delta ,R }_{0}(s)
+ \sum ^{m}_{j=1} \frac{t^{j}}{j!}\mbox{\mathversion{bold}$Q$}_{(j)}(\mbox{\mathversion{bold}$v$}^{\delta ,R}_{0}(s))
\end{eqnarray*}
so that the compatibility conditions are satisfied at each iteration step. By Theorem \ref{thpos}, each \( \mbox{\mathversion{bold}$v$}^{(n)}\) is 
well-defined. 

Since the arguments for the uniform estimate and the convergence with respect to \( n\) and \( R\) are the same as in the case \( \alpha <0\), 
we omit most of the details and just show the basic energy estimate used to derive the uniform estimate.
For any \( \eta >0\) we have 
\begin{eqnarray*}
\begin{aligned}
\frac{1}{2}\frac{{\rm d}}{{\rm d}t}\| \mbox{\mathversion{bold}$v$}^{(n)}_{s}\|^{2}&=-( \mbox{\mathversion{bold}$v$}^{(n)}_{ss}, \mbox{\mathversion{bold}$v$}^{(n)}_{t})\\
&\leq \frac{\alpha }{2}| \mbox{\mathversion{bold}$v$}^{(n)}_{ss}(0)|^{2} - \frac{\delta }{2}\| \mbox{\mathversion{bold}$v$}^{(n)}_{ss}\|^{2}
+ C\|\mbox{\mathversion{bold}$v$}^{(n-1)}_{s}\|^{2}_{1} \\
&\leq \eta \|\mbox{\mathversion{bold}$v$}^{(n)}_{sss}\|^{2} + C_{\eta}\|\mbox{\mathversion{bold}$v$}^{(n)}_{ss}\|^{2}
-\frac{\delta }{2}\| \mbox{\mathversion{bold}$v$}^{(n)}_{ss}\|^{2} + C \| \mbox{\mathversion{bold}$v$}^{(n-1)}_{s}\|^{2}_{1}, \\
\frac{1}{2}\frac{{\rm d}}{{\rm d}t}\|\mbox{\mathversion{bold}$v$}^{(n)}_{ss}\|^{2} &= -(\mbox{\mathversion{bold}$v$}^{(n)}_{sss},\mbox{\mathversion{bold}$v$}^{(n)}_{st})\\
&\leq \frac{\alpha }{2}|\mbox{\mathversion{bold}$v$}^{(n)}_{sss}(0)|^{2} - \delta \| \mbox{\mathversion{bold}$v$}^{(n)}_{sss}\|^{2}
+\eta \| \mbox{\mathversion{bold}$v$}^{(n)}_{sss}\|^{2} + C_{\eta }( \| \mbox{\mathversion{bold}$v$}^{(n)}_{ss}\|^{2} + \|\mbox{\mathversion{bold}$v$}^{(n-1)}_{s}\|^{2}).
\end{aligned}
\end{eqnarray*}
By using the equation and Sobolev's embedding theorem, we have 
\begin{eqnarray*}
|\mbox{\mathversion{bold}$v$}^{(n)}_{sss}(0)|^{2}\leq \eta \| \mbox{\mathversion{bold}$v$}^{(n)}_{sss}\|^{2} + C_{\eta }\|\mbox{\mathversion{bold}$v$}^{(n)}_{ss}\|^{2}
+ C(1+ \|\mbox{\mathversion{bold}$v$}^{(n-1)}_{s}\|_{1})^{2}.
\end{eqnarray*}
Combining all the estimates yields
\begin{eqnarray*}
\sup_{0\leq t\leq T}\| \mbox{\mathversion{bold}$v$}_{s}^{(n)}(t)\|^{2}_{1} + \int ^{T}_{0}\| \mbox{\mathversion{bold}$v$}^{(n)}_{s}(t)\|^{2}_{2}{\rm d}t
\leq C\int ^{T}_{0} \|\mbox{\mathversion{bold}$v$}^{(n-1)}_{s}(t)\|^{2}_{1} {\rm d}t + CT,
\end{eqnarray*}
where in the above estimates, the constants \( C\) depend on \( \|\mbox{\mathversion{bold}$v$}^{(n-1)}\|_{L^{\infty}(\mathbf{R}_{+})}\). From the above, 
estimates uniform in \( n\) and \( R\) can be obtained by induction with respect to \( n\). 
%
%
%


\subsection{Uniform Estimate of the Solution with respect to \( \delta \)}

As before, we derive a uniform estimate. 
First we prove the following.
\begin{lm}
If \( \mbox{\mathversion{bold}$v$}\) is a solution of {\rm (\ref{dpos2})}
with \( \mbox{\mathversion{bold}$v$}_{s}\in C\big([0,T], H^{2}(\mathbf{R}_{+})\big) \), \( \mbox{\mathversion{bold}$v$}\in C\big( [0,T]; L^{\infty}(\mathbf{R}_{+})\big) \), 
and \( |\mbox{\mathversion{bold}$v$}^{\delta}_{0}|=1\), then \( |\mbox{\mathversion{bold}$v$}|=1\) 
in \( \mathbf{R}_{+}\times [0,T]\).
\label{unit2}
\end{lm}
{\it Proof}. As in the proof of Lemma \ref{unit}, if we set \( h(s,t):= |\mbox{\mathversion{bold}$v$}(s,t)|^{2}-1\), \( h\) satisfies
\begin{eqnarray*}
\left\{ 
\begin{array}{ll}
h_{t}=\alpha h_{sss} + \delta h_{ss} + (2\delta |\mbox{\mathversion{bold}$v$}_{s}|^{2} + 3\alpha (\mbox{\mathversion{bold}$v$}_{s}
\cdot \mbox{\mathversion{bold}$v$}_{ss}))h, & s>0,t>0, \\
h(s,0)=0, & s>0, \\
h(0,t)=0, & t>0, \\
h_{s}(0,t)=0, & t>0.
\end{array}\right.
\end{eqnarray*}
We estimate as follows. For any \( \eta >0\), 
\begin{eqnarray*}
\frac{1}{2}\frac{{\rm d}}{{\rm d}t}\| h\|^{2}\leq C_{\eta}\| h\|^{2}-\delta \| h_{s}\|^{2} + \eta \| h_{s}\|^{2}
\end{eqnarray*}
holds. Thus, after choosing \( \eta >0\) sufficiently small, \( h\equiv 0\) follows. This finishes the proof of the lemma. \hfill \( \Box\) 

\bigskip

As before, we rewrite the nonlinear terms in (\ref{dpos2}) into its original form.
\begin{eqnarray}
\left\{
\begin{array}{ll}
\mbox{\mathversion{bold}$v$}_{t} = \mbox{\mathversion{bold}$v$}\times \mbox{\mathversion{bold}$v$}_{ss}
+ \alpha \big\{ \mbox{\mathversion{bold}$v$}_{sss} + \frac{3}{2}\mbox{\mathversion{bold}$v$}_{ss}\times (\mbox{\mathversion{bold}$v$}\times \mbox{\mathversion{bold}$v$}_{s})
+ \frac{3}{2}\mbox{\mathversion{bold}$v$}_{s}\times (\mbox{\mathversion{bold}$v$}\times \mbox{\mathversion{bold}$v$}_{ss})\big\} \\
\hspace*{7.5cm}+ \delta \big( \mbox{\mathversion{bold}$v$}_{ss} + |\mbox{\mathversion{bold}$v$}_{s}|^{2}\mbox{\mathversion{bold}$v$}\big),  & s>0, t>0, \\
\mbox{\mathversion{bold}$v$}(s,0)= \mbox{\mathversion{bold}$v$}^{\delta }_{0}, & s>0, \\
\mbox{\mathversion{bold}$v$}(0,t) = \mbox{\mathversion{bold}$e$}_{3}, & t>0, \\
\mbox{\mathversion{bold}$v$}_{s}(0,t) = \mbox{\mathversion{bold}$0$}, & t>0.
\end{array}\right.
\label{dpos222}
\end{eqnarray}
Now, we derive a basic uniform estimate with respect to \( \delta \). The main method and properties used to derive the estimate is the same as in the case \( \alpha <0\), namely,
utilizing \( |\mbox{\mathversion{bold}$v$}|=1\),
(\ref{par}), and (\ref{par2}), but the energy is slightly modified and we do not need to use a higher order conserved quantity. First we have
\begin{eqnarray*}
\begin{aligned}
\frac{1}{2}\frac{{\rm d}}{{\rm d}t}\| \mbox{\mathversion{bold}$v$}_{s}\|^{2}
&= \frac{\alpha }{2}|\mbox{\mathversion{bold}$v$}_{ss}(0)|^{2} -\delta \| \mbox{\mathversion{bold}$v$}_{ss}\|^{2}, \\
\frac{1}{2}\frac{{\rm d}}{{\rm d}t}\| \mbox{\mathversion{bold}$v$}_{ss}\|^{2}
&\leq C\| \mbox{\mathversion{bold}$v$}_{s}\|^{2}_{2}(1+ \|\mbox{\mathversion{bold}$v$}_{s}\|_{2}) + \frac{\alpha }{2}|\mbox{\mathversion{bold}$v$}_{sss}(0)|^{2}
-\delta \| \mbox{\mathversion{bold}$v$}_{sss}\|^{2} \\
&\leq C\| \mbox{\mathversion{bold}$v$}_{s}\|^{2}_{2}(1+ \|\mbox{\mathversion{bold}$v$}_{s}\|_{2}) - \delta \| \mbox{\mathversion{bold}$v$}_{sss}\|^{2},
\end{aligned}
\end{eqnarray*}
where we have used \( |\mbox{\mathversion{bold}$v$}_{sss}(0)|^{2}\leq C\| \mbox{\mathversion{bold}$v$}_{s}\|^{2}_{2}(1+ \|\mbox{\mathversion{bold}$v$}_{s}\|_{2})\),
which follows by rewriting the boundary term using the equation. To close the estimate, we will derive estimates for 
\( \mbox{\mathversion{bold}$v$}_{sss}\). However, like the estimates above, the boundary terms have a bad sign unlike in the case \( \alpha <0\). 
Thus, we must modify the energy to obtain the desired estimate. Specifically, to obtain an estimate for \( \mbox{\mathversion{bold}$v$}_{sss}\),
we use the following.
\begin{eqnarray*}
\begin{aligned}
\frac{1}{2}\frac{{\rm d}}{{\rm d}t}\bigg\{ \| \mbox{\mathversion{bold}$v$}_{sss}\|^{2} + \frac{2}{\alpha }
(\mbox{\mathversion{bold}$v$}\times \mbox{\mathversion{bold}$v$}_{ss},\mbox{\mathversion{bold}$v$}_{sss})\bigg\}
\leq C\|\mbox{\mathversion{bold}$v$}_{s}\|^{2}_{2}(1+ \|\mbox{\mathversion{bold}$v$}_{s}\|^{2}_{2}).
\end{aligned}
\end{eqnarray*}
In each estimate, \( C\) is independent of \( \delta \).
Combining the three estimates, we obtain a uniform estimate for \( \| \mbox{\mathversion{bold}$v$}_{s}\|_{2}\) for sufficiently small \( \delta \). 
We denote this threshold as \( \delta _{*}\).

The reason we modified the energy from the standard Sobolev norm is to take care of the boundary term. If we directly estimate 
\( \| \mbox{\mathversion{bold}$v$}_{sss}\|^{2}\), boundary term of the form 
\( \mbox{\mathversion{bold}$v$}_{sss}(0)\cdot \partial ^{5}_{s}\mbox{\mathversion{bold}$v$}(0) \) comes out and the 
order of derivative is too high to estimate. By adding a lower order modification term in the energy, we can cancel out this term.
This kind of modification is needed every three derivatives. We use the above energy as an example to demonstrate the 
idea behind finding the correct modifying term.
Taking the trace \( s=0\) in the equation yields
\begin{eqnarray*}
\alpha \mbox{\mathversion{bold}$v$}_{sss}(0,t) + (\mbox{\mathversion{bold}$v$}\times \mbox{\mathversion{bold}$v$}_{ss})(0,t) = \mbox{\mathversion{bold}$0$}
\end{eqnarray*}
for any \( t>0\). Thus, replacing 
\( \| \mbox{\mathversion{bold}$v$}_{sss}\|^{2} \) with \( \| \mbox{\mathversion{bold}$v$}_{sss}\|^{2} + \frac{2}{\alpha }
(\mbox{\mathversion{bold}$v$}\times \mbox{\mathversion{bold}$v$}_{ss},\mbox{\mathversion{bold}$v$}_{sss}) \)  changes the boundary term from 
\( \mbox{\mathversion{bold}$v$}_{sss}(0)\cdot \partial ^{5}_{s}\mbox{\mathversion{bold}$v$}(0)\) to 
\( \big( \mbox{\mathversion{bold}$v$}_{sss}(0) + \frac{1}{\alpha }\mbox{\mathversion{bold}$v$}\times \mbox{\mathversion{bold}$v$}_{ss}(0)\big)
\cdot \partial ^{5}_{s}\mbox{\mathversion{bold}$v$}(0) \), which is zero.

We continue the estimate in this pattern. Suppose that we have a uniform estimate \\
\( \sup_{0\leq t\leq T}\|\mbox{\mathversion{bold}$v$}_{s}(t)\|_{2+3(i-1)}\leq M\) for some \( i\geq 1\). 
For \( j=1,2\), we have
\begin{eqnarray*}
\frac{1}{2}\frac{{\rm d}}{{\rm d}t}\| \partial ^{3i+j}_{s}\mbox{\mathversion{bold}$v$}\|^{2}
\leq C( 1+ \| \mbox{\mathversion{bold}$v$}_{s}\|^{2}_{2+3i}),
\end{eqnarray*}
where we have used \( |\partial ^{3(i+1)}_{s}\mbox{\mathversion{bold}$v$}(0)|^{2} \leq C\|\mbox{\mathversion{bold}$v$}_{s}\|^{2}_{2+3i}\).
Here, \( C\) depends on \( M\), but not on \( \delta \). 
Set \( \mbox{\mathversion{bold}$W$}_{(m)}(\mbox{\mathversion{bold}$v$}):=\mbox{\mathversion{bold}$P$}_{(m)}(\mbox{\mathversion{bold}$v$})
-\alpha ^{m}\partial ^{3m}_{s}\mbox{\mathversion{bold}$v$}\), which is \( \mbox{\mathversion{bold}$P$}_{(m)}(\mbox{\mathversion{bold}$v$})\) without 
the highest order derivative term.
Then, the final estimate is
\begin{eqnarray*}
\frac{1}{2}\frac{{\rm d}}{{\rm d}t}\bigg\{ \| \partial ^{3(i+1)}_{s}\mbox{\mathversion{bold}$v$}\|^{2} 
+ \frac{2}{\alpha ^{i+1}}(\mbox{\mathversion{bold}$W$}_{(i+1)}(\mbox{\mathversion{bold}$v$}), \partial ^{3(i+1)}_{s}\mbox{\mathversion{bold}$v$}) \bigg\}
\leq C\| \mbox{\mathversion{bold}$v$}_{s}\|^{2}_{2+3i} + C,
\end{eqnarray*}
where, again, \( C\) depends on \( M\), but not on \( \delta \). Thus, we have proven the following time-local uniform estimate 
by induction.
\begin{pr}
For a natural number \( k\) and \( M>0\), let \( \mbox{\mathversion{bold}$v$}\) be a solution of {\rm (\ref{dpos222})} with 
\( |\mbox{\mathversion{bold}$v$}^{\delta }_{0}|=1\), \( \mbox{\mathversion{bold}$v$}^{\delta }_{0s}\in
H^{2+3k}(\mathbf{R}_{+})\), and \( \| \mbox{\mathversion{bold}$v$}^{\delta}_{0s}\|_{H^{2+3k}(\mathbf{R}_{+})}\leq M\) satisfying 
\( \mbox{\mathversion{bold}$v$}_{s}\in C\big( [0,T];H^{2+3k}(\mathbf{R}_{+})\big)\) and 
\( \mbox{\mathversion{bold}$v$}\in C\big( [0,T];L^{\infty}(\mathbf{R}_{+})\big) \). Then, there is a \( C_{**}>0\) and
\( T_{1}\in (0, T]\) such that for \( 0\leq \delta \leq \delta_{*}\), \( \mbox{\mathversion{bold}$v$}\) satisfies
\begin{eqnarray*}
\sup_{0\leq t\leq T_{1}}\| \mbox{\mathversion{bold}$v$}_{s}(t)\|_{2+3k}\leq C_{**}.
\end{eqnarray*}
Here, \( T_{1}\) depends on \( \|\mbox{\mathversion{bold}$v$}_{0s}\|_{2}\) and \( C_{**}\) is independent of \( \delta \in (0,\delta _{*}]\). 
\label{posest}
\end{pr}
\subsection{ Taking the limit \( \delta \rightarrow +0\)}
Now we take the limit \( \delta \rightarrow +0\). For \( \delta ^{\prime},\delta \in (0,\delta_{*}]\), we set the difference of the corresponding 
solutions as \( \mbox{\mathversion{bold}$V$}:= \mbox{\mathversion{bold}$v$}^{\delta ^{\prime}}-\mbox{\mathversion{bold}$v$}^{\delta}
-(\mbox{\mathversion{bold}$v$}^{\delta ^{\prime}}_{0}-\mbox{\mathversion{bold}$v$}^{\delta}_{0})
\). 
Then, \( \mbox{\mathversion{bold}$V$} \) satisfies
\begin{eqnarray*}
\left\{
\begin{array}{ll}
\mbox{\mathversion{bold}$V$}_{t} = \mbox{\mathversion{bold}$v$}^{\delta ^{\prime}}\times \mbox{\mathversion{bold}$V$}_{ss}
+ \alpha \big\{ \mbox{\mathversion{bold}$V$}_{sss} + 3\mbox{\mathversion{bold}$V$}_{ss}\times ( \mbox{\mathversion{bold}$v$}^{\delta ^{\prime}}
\times \mbox{\mathversion{bold}$v$}^{\delta ^{\prime}}) \big\}
+ \delta ^{\prime}\mbox{\mathversion{bold}$V$}_{ss} + \mbox{\mathversion{bold}$G$}, & s>0, t>0, \\
\mbox{\mathversion{bold}$V$}(s,0)=\mbox{\mathversion{bold}$0$}, & s>0, \\
\mbox{\mathversion{bold}$V$}(0,t)=\mbox{\mathversion{bold}$0$}, & t>0, \\
\mbox{\mathversion{bold}$V$}_{s}(0,t)= \mbox{\mathversion{bold}$0$}, & t>0,
\end{array}\right.
\end{eqnarray*}
where \( \mbox{\mathversion{bold}$G$}  \) is the collection of terms that are lower order in \( \mbox{\mathversion{bold}$V$}\) and depends linearly on 
\( \mbox{\mathversion{bold}$V$}_{0}:= \mbox{\mathversion{bold}$v$}^{\delta ^{\prime}}_{0}-\mbox{\mathversion{bold}$v$}^{\delta}_{0}\). 
By a standard energy method, we have
\begin{eqnarray*}
\frac{1}{2}\frac{{\rm d}}{{\rm d}t}\| \mbox{\mathversion{bold}$V$}\|^{2}_{3}
\leq C\| \mbox{\mathversion{bold}$V$}\|^{2}_{3} + C\big[ (\delta^{\prime} + \delta ) + \| \mbox{\mathversion{bold}$V$}_{0}\|^{2}_{L^{\infty}(\mathbf{R}_{+})}
+ \| \mbox{\mathversion{bold}$V$}_{0s}\|^{2}_{3}\big],
\end{eqnarray*}
where \( C\) depends on \( C_{**}\) defined in Proposition \ref{posest}. Here, we have used identities such as 
\begin{eqnarray*}
\mbox{\mathversion{bold}$v$}^{\delta^{\prime} }_{s}\times \partial ^{4}_{s}(\mbox{\mathversion{bold}$v$}^{\delta ^{\prime}}-\mbox{\mathversion{bold}$v$}^{\delta})
= \mbox{\mathversion{bold}$v$}^{\delta ^{\prime}}_{s}\times \partial ^{4}_{s}\mbox{\mathversion{bold}$v$}^{\delta ^{\prime}}
-\mbox{\mathversion{bold}$v$}^{\delta }_{s}\times \partial ^{4}_{s}\mbox{\mathversion{bold}$v$}^{\delta} - (\mbox{\mathversion{bold}$v$}^{\delta ^{\prime}}-\mbox{\mathversion{bold}$v$}^{\delta})_{s}\times
\partial ^{4}_{s}\mbox{\mathversion{bold}$v$}^{\delta},
\end{eqnarray*}
to obtain the estimate. From this estimate, we see that 
\( \mbox{\mathversion{bold}$v$}^{\delta} \rightarrow \mbox{\mathversion{bold}$v$}\) in \( C\big( [0,T];L^{\infty}(\mathbf{R}_{+})\big) \) and 
\( \mbox{\mathversion{bold}$v$}^{\delta}_{s} \rightarrow \mbox{\mathversion{bold}$v$}_{s}\) in 
\( C\big( [0,T]; H^{2}(\mathbf{R}_{+})\big) \) as \( \delta \rightarrow +0\), and \( \mbox{\mathversion{bold}$v$}\) 
is the solution to (\ref{pos2}).
Combining this with the uniform estimate, we see that
\( \mbox{\mathversion{bold}$v$}_{s}\in \bigcap ^{k}_{j=0} W^{j,\infty}\big( [0,T];H^{2+3(k-j)}(\mathbf{R}_{+})\big) \). 
As before, the uniform estimate is essentially the energy estimate for the limit system, and after an 
approximation argument on the initial datum, the regularity assumption on the initial datum can be relaxed. Thus we have proven Theorem \ref{st1}.

\end{document}